\documentclass[12pt]{article}
\usepackage{amssymb,latexsym}

\title{Chiral polyhedra in ordinary space, II 
}
\author{Egon Schulte
\thanks{Supported, in part, by NSA-grant H98230-04-1-0116. The research was done, in part, while 
the author visited I.H.E.S. in Bures-sur-Yvette, France, for two months in 2002. The author would
like to thank I.H.E.S. for the hospitality. }\\
\small{Northeastern University}\\ 
\small{Boston, MA 02115, USA}\\
\small{schulte@neu.edu}}
\date{Version of \today.}

\advance\voffset by-10mm

\makeatletter
\@addtoreset{equation}{section}
\@addtoreset{figure}{section}
\@addtoreset{table}{section}
\makeatother
\newcommand{\eqref}[1]{(\ref{#1})}

\newenvironment{proof}{\par \medskip \noindent
\emph{Proof.}\quad}{\hspace*{\fill}
$\Box$ \par \vspace*{1ex}}
\newcommand{\bpf}{\begin{proof}}
\newcommand{\epf}{\end{proof}}
\newtheorem{theorem}{Theorem}[section]
\newcommand{\bthm}{\begin{theorem}}
\newcommand{\ethm}{\end{theorem}}
\newtheorem{corollary}[theorem]{Corollary}
\newcommand{\bcor}{\begin{corollary}}
\newcommand{\ecor}{\end{corollary}}
\newtheorem{lemma}[theorem]{Lemma}
\newcommand{\blem}{\begin{lemma}}
\newcommand{\elem}{\end{lemma}}
\newtheorem{proposition}[theorem]{Proposition}
\newcommand{\bprop}{\begin{proposition}}
\newcommand{\eprop}{\end{proposition}}
\newcommand{\bpc}{\begin{picture}}
\newcommand{\epc}{\end{picture}}
\newcommand{\prt}[1]{\\ {#1})\quad}
\newcommand{\prtind}[1]{\indent {#1})\quad}
\newcommand{\beq}{\begin{equation}}
\newcommand{\bqy}{\begin{eqnarray*}}
\newcommand{\bry}{\begin{array}}
\newcommand{\eeq}{\end{equation}}
\newcommand{\eqy}{\end{eqnarray*}}
\newcommand{\ery}{\end{array}}

\newcommand{\RO}{\rm O}

\newcommand{\hole}{{\mkern2mu|\mkern2mu}}

\renewcommand{\geq}{\geqslant}
\renewcommand{\leq}{\leqslant}

\newcommand{\p}{\varphi}
\renewcommand{\r}{\rho}
\newcommand{\sg}{\sigma}

\newcommand{\Ga}{\mathnormal{\Gamma}}
\newcommand{\La}{\mathnormal{\Lambda}}
\newcommand{\Ph}{\mathnormal{\Phi}}
\newcommand{\Ps}{\mathnormal{\Psi}}
\newcommand{\CP}{\mathcal{P}}
\newcommand{\CQ}{\mathcal{Q}}

\newcommand{\s}{{\bf s}}

\newcommand{\BE}{\mathbb{E}}  
\newcommand{\BZ}{\mathbb{Z}}  
\def\mod#1{\allowbreak \mkern9mu ({\rm mod}\mkern6mu #1)}

\newcommand{\seg}{\{\mkern4mu\}}
\def\blend{\mathbin{\rlap={\,\|}}}
\newcommand{\scl}[1]{\langle\mkern1mu{#1}\mkern1mu\rangle}
\newcommand{\tfrac}[2]{{\textstyle\frac{#1}{#2}}}     

\newcommand{\half}{\tfrac{1}{2}}

\begin{document}

\maketitle

\begin{center}
{\em With best wishes to Branko Gr\"unbaum for his 75th birthday.}
\end{center}
\vskip.2in

\begin{abstract}

\noindent
A chiral polyhedron has a geometric symmetry group with two orbits on the flags, such that
adjacent flags are in distinct orbits. Part~I of the paper described the discrete chiral
polyhedra in ordinary Euclidean space $\BE^3$ with finite skew faces and finite skew
vertex-figures; they occur in infinite families and are of types $\{4,6\}$, $\{6,4\}$ and
$\{6,6\}$. Part~II completes the enumeration of all discrete chiral polyhedra in $\BE^3$. There
exist several families of chiral polyhedra of types $\{\infty,3\}$ and $\{\infty,4\}$ with
infinite, helical faces. In particular, there are no discrete chiral polyhedra with finite
faces in addition to those described in Part~I. 
\end{abstract}

\section{Introduction}
\label{intro}

The present paper continues the study of chiral polyhedra in ordinary Euclidean space
$\BE^3$ begun in \cite{schu}. As polyhedra we permit discrete ``polyhedra-like" structures in
$\BE^3$ with finite or infinite, planar or skew, polygonal faces or vertex-figures; more
precisely, a polyhedron in $\BE^3$ is a discrete $3$-dimensional faithful realization of an
abstract regular polytope of rank $3$ (see McMullen \& Schulte~\cite[Ch.5]{msarp}).  A
polyhedron is ({\em geometrically\/}) {\em chiral\/} if its geometric symmetry group has
exactly two orbits on the flags, such that adjacent flags are in distinct orbits. Chiral
polyhedra are nearly regular polyhedra; recall that a (geometrically) regular polyhedron has a
geometric symmetry group with just one orbit on the flags. 

The discrete regular polyhedra in $\BE^3$ were completely enumerated by Gr\"unbaum~\cite{gon}
and Dress~\cite{d1,d2}; see \cite[Section~7E]{msarp} (or \cite{msrpo}) for a quick method of
arriving at the full characterization, as well as for presentations of the symmetry groups.

In Part I, we described a complete classification of the discrete chiral polyhedra with finite
skew faces and finite skew vertex-figures in $\BE^3$ (\cite{schu}). There are three
integer-valued  two-parameter families of chiral polyhedra of this kind for each type
$\{4,6\}$, $\{6,4\}$ or $\{6,6\}$; some infinite families split further into several smaller
subfamilies. Moreover, there are no chiral polyhedra in $\BE^3$ which are finite.

The present Part II completes the enumeration of the discrete chiral polyhedra in $\BE^3$. All
chiral polyhedra not described in Part~I have infinite, helical faces and occur again in
infinite families; their types are $\{\infty,3\}$ or $\{\infty,4\}$, and their faces are
helices over triangles or squares. We discuss their geometry and combinatorics in
detail in Sections~\ref{helfacone}, \ref{helfactwo}, \ref{helfacthree} and~\ref{relships}. The foundations are laid in Section~\ref{spgr}, where we also prove that a chiral polyhedron cannot have finite planar faces or vertex-figures or an affinely reducible symmetry group; in effect, this shows that Part~I
actually enumerated all the chiral polyhedra in $\BE^3$ with finite faces. Basic terminology and
results are reviewed in Section~\ref{bano} (following \cite{msarp,schu}).

Finally, as a birthday greeting, I wish to acknowledge the deep influence Branko Gr\"unbaum's work has had on my research over the years. His pioneering papers \cite{gon, grgcd} on generalized geometric or abstract regular polyhedra or polytopes lie at the heart of the enumeration presented here.

\section{Basic notions}
\label{bano}

Although our main interest is in geometric polyhedra, we begin with a brief review of the 
underlying abstract theory (see \cite[Ch.2]{msarp}). An {\em abstract polyhedron\/} $\CP$ is a
partially ordered set with a strictly monotone {\em rank\/} function with range
$\{-1,0,\ldots,3\}$. The elements of rank $j$ are called the {\em $j$-faces\/} of $\CP$.  For
$j = 0$, $1$ or $2$, we also call $j$-faces {\em vertices}, {\em edges\/} and {\em facets\/},
respectively. When there is no possibility of confusion, we adopt standard terminology for
polyhedra and use the term ``face" to mean ``$2$-face" (facet). The {\em flags\/} of $\CP$
each contain one vertex, one edge and one facet, as well as the unique minimal face $F_{-1}$
and unique maximal face $F_{3}$ of $\CP$ (which usually are omitted from the notation). 
Further, $\CP$ is {\em strongly flag-connected\/}, meaning that any two flags $\Ph$ and $\Ps$
of $\CP$ can be joined by a sequence of flags $\Ph = \Ph_{0},\Ph_{1},\ldots,\Ph_{k} = \Ps$,
where $\Ph_{i-1}$ and $\Ph_{i}$ are {\em adjacent\/} (differ by one face), and $\Ph \cap \Ps
\subseteq \Ph_{i}$ for each $i$.  Finally, if $F$ and $G$ are a $(j-1)$-face and a $(j+1)$-face
with $F < G$ and $0 \leq j \leq 2$, then there are exactly {\em two\/} $j$-faces $H$ such that
$F < H < G$.

An abstract polyhedron $\CP$ is {\em chiral\/} if its ({\em combinatorial automorphism\/}) 
{\em group\/} $\Ga(\CP)$ has two orbits on the flags such that adjacent flags are in distinct
orbits. Let $\Ph := \{F_{0},F_{1},F_{2}\}$ be a fixed flag, or {\em base\/} flag, of $\CP$,
and let $F'_{j}$, with $j=0,1,2$, denote the $j$-face of $\CP$ with $F_{j-1} < F'_j < F_{j+1}$
and $F'_{j} \neq F_j$.  If $\CP$ is chiral, then $\Ga(\CP)$ is generated by {\em distinguished
generators\/} $\sg_{1},\sg_{2}$ ({\em with respect to $\Ph$\/}), where $\sg_{1}$ fixes $F_2$
and cyclically permutes the vertices and edges of $F_2$ such that $F_{1}\sg_{1} = F'_{1}$, and 
$\sg_{2}$ fixes $F_0$ and cyclically permutes the vertices and edges in the vertex-figure at $F_0$ 
such that $F_{2}\sg_{2} = F'_{2}$. (The vertex-figure of $\CP$ at a vertex $F$ consists of the faces 
of $\CP$ with vertex $F$, ordered as in $\CP$. It is isomorphic to a polygon.)  A chiral polyhedron
occurs in two ({\em combinatorially\/}) {\em enantiomorphic forms\/} (see \cite{swc,swp}); an
enantiomorphic form simply is a pair consisting of the underlying abstract polyhedron and an
orbit of flags (specifying a ``combinatorial orientation"). These enantiomorphic forms
correspond to the two (conjugacy classes of) pairs of distinguished generators of $\Ga(\CP)$,
namely the pair $\sg_{1},\sg_{2}$ determined by the original base flag $\Ph$, and the pair
$\sg_{1}\sg_{2}^{2},\sg_{2}^{-1}$ determined by the adjacent flag $\Ph^{2} :=
\{F_{0},F_{1},F_{2}'\}$ (differing from $\Ph$ in the $2$-face).

An abstract polyhedron $\CP$ is {\em regular\/} if $\Ga(\CP)$ is transitive on its flags. If 
$\CP$ is regular, then $\Ga(\CP)$ is generated by {\em distinguished generators\/}
$\r_{0},\r_{1},\r_{2}$ ({\em with respect to $\Ph$\/}), where $\r_{j}$ is the unique
automorphism which keeps all but the $j$-face of a chosen base flag $\Ph$ fixed. The elements
$\sg_{1} := \r_{0}\r_{1}$ and  $\sg_{2} := \r_{1}\r_{2}$ generate the ({\em combinatorial\/})
{\em rotation subgroup\/} $\Ga^{+}(\CP)$ of $\Ga(\CP)$; they have properties similar to those of
the generators for chiral polyhedra. Now the two pairs of generators of $\Ga^{+}(\CP)$ are
conjugate in $\Ga(\CP)$ under $\rho_2$, so the two enantiomorphic forms can be identified.
An abstract polyhedron $\CP$ is {\em directly regular\/} if $\CP$ is regular and $\Ga^{+}(\CP)$ has index $2$ in $\Ga(\CP)$. An abstract regular polyhedron $\CP$ may not be directly regular. For example, if $\CP$ has Petrie polygons of odd length, then $\Ga^{+}(\CP) = \Ga(\CP)$. (Recall that a {\em Petrie polygon\/} of $\CP$ is a path along edges such that any two, but no three, consecutive edges lie in a face.)

We now turn to the geometric theory. Let $\CP$ be an abstract polyhedron, and let $\CP_j$ denote 
its set of $j$-faces.  A {\em realization\/} of $\CP$ is a mapping
$\beta\colon\CP_{0}\rightarrow E$ of the vertex-set $\CP_{0}$ into some euclidean space $E$
(see \cite[Section 5A]{msarp}). In the present context, $E = \BE^3$. Define $\beta_0 := \beta$
and $V_{0} := V := \CP_{0}\beta$, and write $2^X$ for the family of subsets of a set $X$.
Then $\beta$ recursively induces surjections $\beta_{j}\colon \CP_{j}\rightarrow V_{j}$, for
$j=1,2,3$, with $V_{j}\subset 2^{V_{j-1}}$ consisting of the elements
\[ F\beta_{j} := \{G\beta_{j-1} \mid G\in \CP_{j-1} \mbox{ and } G\leq F\} \]
for $F\in\CP_j$; further, $\beta_{-1}$ is given by $F_{-1}\beta_{-1} := \emptyset$. Even
though each $\beta_j$ is determined by $\beta$, it is helpful to think of the realization as
given by all the $\beta_j$. A realization $\beta$ is {\em faithful\/} if each $\beta_j$ is a
bijection; otherwise, $\beta$ is {\em degenerate\/}.  In a {\em symmetric\/} realization $\beta$ of 
$\CP$, each automorphism of $\CP$ induces an isometric permutation of the vertex-set $V$; 
such an isometric permutation extends to an isometry of $E$, uniquely if $E$ is the affine
hull of $V$. Thus associated with a symmetric realization $\beta$ of $\CP$ is a representation of
$\Ga(\CP)$ as a group of euclidean isometries.

We mostly work with discrete and faithful realizations.  In this case the vertices, edges and 
facets of $\CP$ are in one-to-one correspondence with certain points, line segments and simple
(finite or infinite) polygons in $E$, and it is safe to identify a face of $\CP$ and its image
in $E$. The resulting family of points, line segments and polygons is a (discrete) {\em
geometric polyhedron\/} in $E$ and is denoted by $P$; it is understood that $P$ inherits the
partial ordering of $\CP$, and when convenient $P$ will be identified with $\CP$. We let $G(P)$
denote the (geometric) symmetry group of $P$.

Recall that a geometric polyhedron $P$ in $\BE^3$ is {\em geometrically regular\/} if $G(P)$ is 
flag-transitive, and that $P$ is {\em geometrically chiral\/} if $G(P)$ has two orbits on the
flags of $P$ such that adjacent flags are in distinct orbits (see \cite{schu}). In the latter case
the underlying abstract polyhedron $\CP$ must be (combinatorially) chiral or regular. Note that
there are geometric polyhedra which are not chiral but still have a symmetry group with only two
orbits on the flags (see \cite{gruen,wills}). One finite example of type $\{6,6\}$ has $20$
vertices and $20$ faces and has full icosahedral symmetry (its map is listed as $W\#60.57$ in
\cite{wilson}). For a discussion of other classes of highly-symmetric polyhedra see, for example,
\cite{grhol,gr}. 

The structure results for abstract polyhedra carry over to geometric polyhedra as follows. Let $P$ 
be a discrete (geometrically) chiral or regular (geometric) polyhedron in $\BE^3$.  As before,
let $\Ph := \{F_{0},F_{1},F_{2}\}$ be a base flag of $P$, and let $F'_{j}$, with $j=0,1,2$, be
the $j$-face of $P$ with $F_{j-1} < F'_j < F_{j+1}$ and $F'_{j} \neq F_j$.  

If $P$ is chiral, then $G := G(P)$ is generated by {\em distinguished generators\/} $S_{1},S_{2}$ 
({\em with respect to\/} $\Ph$), where $S_{1}$ fixes the base facet $F_2$ and cyclically
permutes its vertices such that $F_{1}S_{1} = F'_{1}$ (and thus $F'_{0}S_{1} = F_{0}$), and
$S_{2}$ fixes the base vertex $F_0$ and cyclically permutes the vertices in the vertex-figure
at $F_0$ such that $F_{2}S_{2} = F'_{2}$ (and thus $F'_{1}S_{2} = F_{1}$). Then 
\beq  
\label{agchi}
S_{1}^p = S_{2}^q = (S_{1}S_{2})^{2}  = I, 
\eeq
the identity mapping, where $p$ and $q$ are determined by the ({\em Schl\"{a}fli\/}) {\em type\/}
$\{p,q\}$ of $P$ (see \cite{cm}); in  general there are also other independent relations. The
involution $T := S_{1}S_{2}$ interchanges the two vertices in $F_{1}$ as well as the two faces
meeting at $F_{1}$. 

If $P$ is regular, then $G(P)$ is generated by {\em distinguished generators\/} $R_{0},R_{1},R_{2}$ 
({\em with respect to $\Ph$\/}), where $R_j$ is the unique symmetry of $P$ which maps all but
the $j$-face of $\Ph$ to itself. Each $R_j$ is a reflection in a point, a line or a plane. If
$S_{1} := R_{0}R_{1}$ and $S_{2} := R_{1}R_{2}$, then the subgroup $G := \scl{S_{1},S_{2}}$ of
$G(P)$ of index $2$ has  properties very similar to those of the group of a chiral polyhedron.
Now $T := S_{1}S_{2} = R_{0}R_{2}$, and $T$ again interchanges the two vertices of $F_{1}$ as
well as the two faces meeting at $F_{1}$. For a regular polyhedron, the {\em dimension
vector\/} 
\[ (\dim R_{0},\dim R_{1},\dim R_{2}) \] 
records the dimensions of the reflection mirrors for $R_{0},R_{1},R_{2}$. See 
\cite[Section 7E]{msarp} for a complete enumeration of the regular polyhedra in $\BE^3$.  

For emphasis, throughout the paper, $G$ will denote the subgroup $\scl{S_{1},S_{2}}$ of $G(P)$, 
irrespective of whether $P$ is chiral or regular.  Then $G=G(P)$ if $P$ is chiral, or $G$ is a
subgroup of index $2$ in $G(P)$ if $P$ is regular.

We require the following lemma, which gives a useful necessary and sufficient condition for a
polyhedron to be regular.

\blem
\label{regcrit}
Let $P$ be a discrete chiral or regular polyhedron in $\BE^3$ with associated group  
$G=\scl{S_{1},S_{2}}$. Let $\Ph = \{F_{0},F_{1},F_{2}\}$ be the corresponding base flag, and
let $F_{j}'$, for $j=0,1,2$, be the $j$-face associated with $\Ph$ as above. Then $P$ is regular 
if and only if there exists an isometry $R$ of $\BE^3$ such that
\[ R^{-1}S_{2}R = S_{2}^{-1}, \quad R^{-1}TR = T, \quad 
F_{0}R=F_{0}, \quad F_{1}R=F_{1}, \quad F_{2}R=F_{2}' . \]
\elem

\bpf
If $P$ is regular, then $R:=R_{2}$ has the required properties. Conversely, suppose that $R$ is
an isometry acting in the way described. Then $R^{-1}G R = G$ because $G$ is generated by
$S_{2}$ and $T$, and 
\[ (F_{i} G)R = (F_{i}R)G = F_{i}G \mbox{ or } F_{i}'G \]
according as $i=0,1$ or $i=2$. But since there is only one orbit on the $i$-faces for each $i$
(even when $i=2$), this implies that $R$ must indeed be a symmetry of $P$, which then takes
$\Phi$ to the adjacent flag $\Phi^{2} = \{F_{0},F_{1},F_{2}'\}$. It follows that $P$ must be
regular, and that $R=R_{2}$.
\epf

Note that Lemma~\ref{regcrit} can be employed to shorten the proof of Theorem~3.1 of
\cite{schu},  which says that there are no finite chiral polyhedra in $\BE^3$.

Geometrically chiral or regular polyhedra can be obtained by a variant of {\em Wythoff's
construction\/} (see \cite{crp} or \cite[Section 5A]{msarp}). Let an abstract polyhedron $\CP$ be
chiral or regular, and let $G := \scl{S_1,S_2}$ be a euclidean representation of its group
$\Ga(\CP)$ or its rotation subgroup $\Ga^{+}(\CP)$, respectively.  Each point $v$ which is fixed by
$S_2$ can serve as the {\em initial\/} vertex of a realization $P$ with vertex-set $V=vG$. Its base
vertex, base edge and base face are $v$, $v\scl{T}$ or $v\scl{S_1}$, respectively, and the
other vertices, edges and faces are their images under $G$. In our applications, $P$ will
generally be a geometric polyhedron. Note that, a priori, an abstract regular polyhedron can have a
realization which is geometrically chiral.

Certain operations that can be applied to regular polyhedra have analogs that also apply to 
abstract chiral polyhedra and frequently to geometrically chiral polyhedra as well. In
particular, we employ the {\em duality operation\/} $\delta$ and the ({\em $2$nd\/}) {\em
facetting operation\/} $\p_{2}$, yielding the generators for the groups of the dual polyhedron
$\CP^{*} = \CP^{\delta}$ or the polyhedron $\CP^{\p_{2}}$, respectively
(see \cite[p.192,194]{msarp}). If an abstract polyhedron $\CP$ is chiral, then $\delta$ and
$\p_{2}$ are given by
\beq
\label{dualop}
\delta\colon (\sg_{1},\sg_{2}) \mapsto (\sg_{2}^{-1},\sg_{1}^{-1}) 
\quad\mbox{and}\quad
\p_{2}\colon (\sg_{1},\sg_{2}) \mapsto (\sg_{1}\sg_{2}^{-1},\sg_{2}^{2}), 
\eeq 
respectively. When applied to the rotation subgroup $\Ga^{+}(\CP)$ of a regular polyhedron $\CP$,
the new generators (on the right) become the distinguished generators for the rotation subgroups of
the polyhedra $\CP^{\delta}$ or $\CP^{\p_{2}}$ of \cite[Section~7B]{msarp}, respectively. Both 
$\delta$ and $\p_{2}$ carry over to geometric polyhedra. The operations in \eqref{dualop} are
examples of {\em mixing operations\/}, in which a new group is derived from a given group by
taking as generators certain suitably chosen products of the generators of the given group, so that
the new group is a subgroup (see \cite[p.83]{msarp}).

Suppose $\CP$ (or $P$) is a chiral polyhedron of type $\{p,q\}$ with $p$ and $q$ finite. For
odd $q$, the geometric effect of the operation $\p_{2}$ on $\CP$ is as described in
\cite[p.195]{msarp} (for regular polyhedra), so we will not repeat it here. However, in the
present context we are only interested in the case when $q\geq 6$ is even, and specifically
when $q=6$. Here it is helpful to think of $\CP$ as a map on a (in our applications,
non-compact) surface $M$, and of $\Ga(\CP)$ as a group of homeomorphisms which preserve this
map and its underlying order complex (``combinatorial barycentric subdivision").  The
triangles in the order complex correspond to the flags of $\CP$, so in particular we have a
base triangle $\Ph$ and a corresponding adjacent triangle $\Ph^2$ as well as the two pairs of
generators $\sg_{1},\sg_{2}$ and $\sg_{1}\sg_{2}^{2},\sg_{2}^{-1}$ of $\Ga(\CP)$ associated
with them. Moreover, the union of the triangles $\Ph$ and $\Ph^2$ is a fundamental region for
the action of $\Ga(\CP)$ on $M$. The map $\CP$ on $M$ can be recovered from this action by
employing (a combinatorial analogue of) Wythoff's construction, with initial vertex the base
vertex in $\Ph$; the latter is fixed by $\sg_{2}$ (and $\sg_{2}^{-1}$).  

Now we can find a model for the new polyhedron $\CP^{\p_{2}}$ as follows. By \eqref{dualop},
the group of $\CP^{\p_{2}}$ is generated by the pair $\sg_{1}\sg_{2}^{-1},\sg_{2}^{2}$, so we can
apply Wythoff's construction with the same initial vertex as before. Then the vertices and edges of
$\CP^{\p_{2}}$ are among those of $\CP$, and a typical face of $\CP^{\p_{2}}$ is a {\it hole\/} (more
exactly, {\em $2$-hole}), which is formed by the edge-path which leaves a vertex by the $2$nd edge
from which it entered, in the same sense (that is, keeping always to the right, say, in some local
orientation of $M$). Generally, since $q$ is even, only every other edge of $\CP$ at a vertex
of $\CP^{\p_{2}}$ belongs to such a face of $\CP^{\p_{2}}$, so a vertex of $\CP^{\p_{2}}$ lies in
$q/2$ faces of $\CP^{\p_{2}}$.  The same polyhedron can also be obtained from Wythoff's construction
applied with the alternative  pair of generators $\sg_{1}\sg_{2}^{3},\sg_{2}^{-2}$ of the (same)
group; this pair is the image under $\p_{2}$ of the alternative pair
$\sg_{1}\sg_{2}^{2},\sg_{2}^{-1}$ of $\Ga(\CP)$. 

With appropriate changes, these considerations extend to regular polyhedra $\CP$ (or $P$) and 
their rotation subgroups $\Ga^{+}(\CP)$.

For a survey about realizations of regular or chiral polytopes in euclidean spaces of small
dimension the reader is also referred to \cite{mssmall}.

\section{The special group}
\label{spgr}

We shall assume from now on that $P$ is a discrete {\sl infinite\/} polyhedron, or {\em
apeirohedron\/}, which is chiral or regular; in fact, there are no finite chiral polyhedra in
$\BE^3$ (see \cite[\S 3]{schu} or \cite{mc}), and the finite regular polyhedra are all known. Then 
$G = \scl{S_{1},S_{2}}$ must be an infinite discrete group of isometries. For now we assume
that $G$ is ({\em affinely\/}) {\em irreducible\/} on $\BE^3$, meaning that there is no
non-trivial linear subspace $L$ of $\BE^{3}$ which is {\em invariant\/} in the sense that $G$
permutes the translates of $L$. (If such an invariant subspace $L$ exists, then its orthogonal
complement $L^{\perp}$ is also invariant in the same sense.)  It follows that $G$ must be a
crystallographic group. 

If $R: x\mapsto xR' + t$ is a general element of $G$, with $R' \in \RO_{3}$, the orthogonal
group, and $t \in \BE^{3}$ a translation vector, then the linear mappings $R'$ form the
{\em special group\/} $G_0$ of $G$. This is a finite subgroup of $\RO_{3}$.  (See \cite{gb} for the 
complete enumeration of the finite subgroups of $\RO_{3}$.)  In particular, if $T(G)$ denotes the
subgroup of all translations in $G$, then $G_{0} \cong G/T(G)$. Note that
$S_{2}' = S_{2}$ and
\[ T' = S_{1}'S_{2}' = S_{1}'S_{2} . \]
Moreover, 
\[ G_{0} = \scl{S_{1}',S_{2}}. \] 

The following refinement of Bieberbach's theorem (see \cite{bi} or \cite[\S 7.4]{ra}) applies to $G$
and was proved in \cite[p.220]{msarp}. (Note that Lemma~7E6 of \cite{msarp} states the lemma in a
slightly weaker form, but its proof actually establishes the stronger version.)

\blem
\label{disrot}
The special group of an irreducible infinite discrete group of isometries in $\BE^{2}$ or
$\BE^{3}$ does not contain rotations of periods other than $2$, $3$, $4$ or $6$.
\elem

The generators $S_{1}$ and $S_{2}$ of $G$ determine the geometry of the faces and
vertex-figures of $P$, respectively. Part I of the paper dealt with the case that $S_{1}$ and
$S_{2}$ are rotatory reflections of finite period. When chiral, the corresponding polyhedra in
$\BE^3$ have finite skew faces and finite skew vertex-figures; however, when regular, either the 
faces or the vertex-figures are planar. Note in this context that an element in $G$ is a
rotatory reflection if and only if its image in $G_0$ is a rotatory reflection (of the same
period).  We now investigate the possibility that a generator $S_{1}'$ or $S_{2}$ of $G_0$ is a
rotation of period at least $3$; the corresponding generator $S_{1}$ or $S_{2}$ of $G$ would
then be a rotation or {\em twist\/} (screw motion of infinite period), and the faces or
vertex-figures of $P$ would be planar or helical. (Of course, helical vertex-figures are ruled out
by our discreteness assumption.)

Before we proceed, note that $S_{1}'$ could only be a rotation of period $2$ if $P$ had (planar) 
zig-zag faces. We shall later exclude this possibility (see Lemma~\ref{chirreds}). For now we just
observe that the group $G$ of such a polyhedron $P$ must necessarily act reducibly on $\BE^3$. In
fact, $S_{1}$ is necessarily a ``half-twist" along the axis of the base (zig-zag) face $F_{2}$
(given by a half-turn about this axis, followed by a translation along this axis), and $T$ is either
the plane reflection in the perpendicular bisector of the base edge $F_{1}$, or a half-turn about a
line through the midpoint of $F_{1}$ and  perpendicular to $F_{1}$. In either case, $G_{0} =
\scl{S_{1}',T'}$ acts reducibly (the axis of $F_{2}$ passes through the midpoint of $F_{1}$),
and so does $G$.

Moreover, we can also eliminate the possibility that $S_{1}' = I$. In fact, in this case $P$ must
have linear apeirogons as faces, forcing $P$ to entirely lie on a line (adjacent faces must lie on
the same line).

\begin{figure}[hbt]
\centering
\begin{center}
\begin{picture}(260,200)
\put(40,5){  
\thinlines
\begin{picture}(180,45)
\multiput(0,0)(120,0){2}{\circle*{4}}
\multiput(50,33.3)(120,0){2}{\circle*{4}}
\multiput(0,0)(50,33.3){2}{\line(1,0){120}}
\multiput(0,0)(120,0){2}{\line(3,2){50}}
\put(0,0){\line(1,0){120}}
\end{picture}}
\put(40,125){  
\begin{picture}(180,45)
\thinlines
\multiput(0,0)(120,0){2}{\circle*{4}}
\multiput(50,33.3)(120,0){2}{\circle*{4}}
\multiput(0,0)(50,33.3){2}{\line(1,0){120}}
\multiput(0,0)(120,0){2}{\line(3,2){50}}
\end{picture}}
\put(40,5){ 
\begin{picture}(180,225)
\thinlines
\multiput(0,0)(120,0){2}{\line(0,1){120}}
\multiput(50,33.3)(120,0){2}{\line(0,1){120}}
\put(170,33.3){\line(0,1){120}}
\end{picture}}
\put(40,5){ 
\begin{picture}(180,225)
\thicklines
\multiput(50,33.3)(70,-33.3){2}{\circle*{7}}
\multiput(0,120)(170,33.3){2}{\circle*{7}}
\multiput(0,120)(1,-1){120}{\circle*{2}}
\multiput(0,120)(0.66,-1.15){75}{\circle*{2}}
\multiput(0,120)(1.42,0.28){120}{\circle*{2}} 
\multiput(120,0)(0.42,1.28){120}{\circle*{2}} 
\multiput(50,33.3)(1,1){120}{\circle*{2}} 
\multiput(50,33.3)(1.16,-0.55){60}{\circle*{2}} 
\end{picture}}
\put(40,5){ 
\begin{picture}(180,225)
\thinlines
\put(85,76.6){\circle*{5}}
\put(85,76.6){\vector(1,0){110}}
\put(199,74){$\xi_2$}
\put(85,76.6){\vector(0,1){100}}
\put(82,182){$\xi_{3}$}
\put(85,76.6){\vector(-3,-2){60}}
\put(15,31){$\xi_{1}$}
\put(60,60){\circle{4}}
\put(145,76.66){\circle{4}}
\put(85,136.66){\circle{4}}
\end{picture}}
\end{picture}
\caption{\it A regular tetrahedron inscribed in a cube.}
\label{tetcub}
\end{center}
\end{figure}
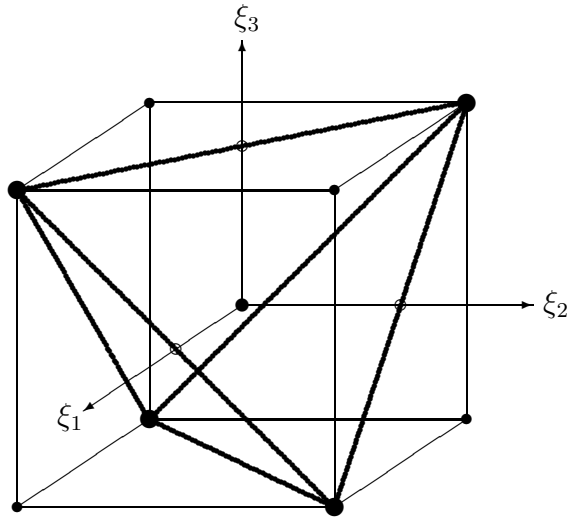

\blem
\label{specgroups}
Let $P$ be a discrete chiral or regular apeirohedron in $\BE^3$ with base vertex $o$, and
let the corresponding group $G = \scl{S_1,S_2}$ act irreducibly, with special group $G_{0} = \scl{S_1',S_2}$. 
\prt{a} Then the generators $S_{1}$ and $S_{2}$ of $G$ cannot both be rotations. 
\prt{b} If one of the generators $S_{1}'$ and $S_{2}$ of $G_{0}$ is a rotatory reflection and
the other a rotation of period at least $3$, then $G_{0}$ can only be one of the following
three groups:\\
\prtind{1}  $\,[3,4]\;(\cong S_{4}\times C_{2})$, the group of a cube $\{4,3\}$ centered at
$o$; \\
\prtind{2}  $[3,3] \;(\cong S_{4})$, the subgroup of $[3,4]$ associated with a regular tetrahedron
$\{3,3\}$ 
\indent inscribed in $\{4,3\}$ (as in Figure~\ref{tetcub}); or \\
\prtind{3}  the subgroup 
\[ [3,3]^{*} : = [3,3]^{+} \cup (-I) [3,3]^{+} \;(\cong A_{4} \times C_{2})\]
\indent of $[3,4]$, obtained by adjoining $-I$ to the rotation subgroup $[3,3]^{+}$ of $\{3,3\}$.
\elem

\bpf
Let $P$ be of type $\{p,q\}$; clearly, since $P$ is discrete, $q$ must be finite. Let $p'$ denote
the period of $S_{1}'$. Most values of $p$, $q$ and $p'$ can be ruled out by appealing to
Lemma~\ref{disrot}, or by observing that the group abstractly defined by the relations in
\eqref{agchi} would be finite. 

First we eliminate the possibility that both generators of $G$ are rotations. If $S_{1}$ and $S_{2}$ 
are rotations, then $G_{0} = \scl{S_{1}',S_{2}}$ must be a finite group in $\BE^3$ generated
by two rotations with distinct axes (recall that $G$ is irreducible).  Hence, $G_{0} =
[3,3]^{+}$ or $[3,4]^{+}$. Since now $p'=p$ and every rotation in these groups is of order $3$ or
$4$, we can only have an infinite group $G$ if $G_{0} = [3,4]^{+}$ with $p=q=4$. However, this
possibility is also excluded because the product of two rotations of period $4$ about distinct axes in
$[3,4]^{+}$ cannot be an involution (on the other hand, $T'= S_{1}'S_{2}$ would have to be an
involution). This settles the first part of the lemma.

We now consider the case that $S_{1}'$ is a rotation of period at least $3$ and that $S_{2}=S_{2}'$ 
is a rotatory reflection. Then $T'$ is a plane reflection; note that we cannot have $T'=-I$
because $G$ is irreducible. We now apply a similar trick to that in the proof of \cite[Theorem
7E4]{msarp} and consider the new mapping 
\[ -S_{2}=(-I)S_{2} = -S_{2}' , \]
which is a rotation about the same axis as the rotational component of $S_{2}$. Then 
\[ \widehat{G}_{0} := \scl{S_{1}',-S_{2}} \] 
is a finite group generated by two rotations with distinct axes, and so we must have
$\widehat{G}_{0} = [3,3]^{+}$ or $[3,4]^{+}$. Note that $G_{0}$ is a subgroup of
$\widehat{G}_{0}\cdot \scl{-I}$ and thus of $[3,3]^{*}$ or $[3,4]$, respectively;  if $-S_{2}$ is 
of period $3$, then $-I = S_{2}^3 \in G_{0}$ and hence $G_{0} = [3,3]^{*}$ or $[3,4]$. Three
cases can occur.

If $\widehat{G}_{0} = [3,3]^{+}$, then $G_{0} = [3,3]^{*}$. In this case, the elements 
$S_{1}',-S_{2}$ are standard generators of the rotation group $[3,3]^+$ of $\{3,3\}$, meaning
that $S_{1}'$ rotates in a face of $\{3,3\}$ and $-S_{2}$ rotates about a vertex of this face,
with orientations such that their product is a half-turn about the midpoint of an edge of this
face containing the vertex. Now $p=3$ or $\infty$, and $q=6$ (already the rotational
component of $S_{2}$ is of period $6$).

If $\widehat{G}_{0} = [3,4]^{+}$, then $S_{1}'$ and $-S_{2}$ must be rotations of period $3$ or $4$;
however, they cannot have the same period. In fact, two rotations of period $3$ generate a
group smaller than $[3,4]^{+}$, and the product of two distinct rotations of period $4$ in
$[3,4]^{+}$ cannot be an involution.  Moreover, we can also rule out the case $p=p'=3$, because
then $G$ would be finite when $q=4$. This leaves the two possibilities that either $(p',q)=(3,4)$
with $p=\infty$ (the case $p=3$ would again yield a finite group), or $(p',q)=(4,6)$ with $p=4$
or $\infty$. 

In the first case, the elements $S_{1}',-S_{2}$ are standard generators (in the above sense) of
$[3,4]^+$, viewed as the rotation group of an octahedron $\{3,4\}$. Now the generators
$S_{1}',S_{2}$ of $G_{0}$ are symmetries of a regular tetrahedron inscribed in $\{3,4\}$ (with
its vertices at face centers of $\{3,4\}$), so that $G_{0}=[3,3]$. 

Similarly, in the second case, the elements $S_{1}',-S_{2}$ are standard generators of
$[3,4]^+$, now viewed as the rotation group of the cube $\{4,3\}$. Now $G_{0} = [3,4]$, because
$-S_{2}$ has period~$3$. 

Finally, a similar analysis goes through in the dual case when $S_{1}'$ is a rotatory
reflection (of period $p'=p$) and $S_{2}$ is a rotation (of period at least $3$). Now we take
\[ \widehat{G}_{0} := \scl{-S_{1}',S_{2}} \] 
as the new group. 

As before, if $\widehat{G}_{0} = [3,3]^{+}$, then $G_{0} = [3,3]^{*}$ and $-S_{1}',S_{2}$
are standard generators of $[3,3]^{+}$. Now $p=6$ (already the rotational component of $S_{1}$
is of period $6$) and $q=3$.  

If $\widehat{G}_{0} = [3,4]^{+}$, then again $-S_{1}'$ and $S_{2}$ are rotations of distinct
periods, $3$ or $4$. We can rule out the case $p=4$ (occurring if $-S_{1}'$ has period $4$), 
because then $G$ would be finite when $q=3$. This only leaves the possibility that
$(p,q)=(6,4)$, obtained when $-S_{1}',S_{2}$ are standard generators of $[3,4]^+$, viewed as
the rotation group of $\{3,4\}$. Now $G_{0} = [3,4]$, because $-S_{1}'$ has period~$3$. 

In summary, the only possible special groups are $[3,3]$, $[3,4]$ and $[3,3]^{*}$, with generators 
$S_{1}',S_{2}$ specified as above.
\epf

The next two lemmas employ the geometry of $P$ to eliminate the possibility of finite faces
in Lemma~\ref{specgroups}(b). In each case we exhibit representations of the putative generators
$S_{1},S_{2}$ of $G$ and prove that these generators must actually determine a regular polyhedron
which either has helical faces or is finite; however, the latter possibility is excluded by our
assumptions on $P$. The results are interesting in their own right.

\blem
\label{twi}
Let $P$ be a discrete chiral or regular apeirohedron in $\BE^3$ with base vertex $o$, and let
the corresponding group $G = \scl{S_1,S_2}$ be irreducible, with special group $G_{0} = \scl{S_1',S_2}$. If $S_{1}'$ is a rotation of period at least $3$ and $S_{2}$ a rotatory reflection, then $P$
is a regular polyhedron with helical faces and $S_{1}$ is a twist. In particular, $P$ is the
Petrial of a Petrie-Coxeter polyhedron $\{4,6\hole 4\}$, $\{6,4\hole 4\}$ or $\{6,6\hole 3\}$;
that is, 
\[ P = \{\infty,6\}_{4,4},\, \{\infty,4\}_{6,4}\, \mbox{ or } \{\infty,6\}_{6,3}, \] 
respectively (notation as in \cite[p.224]{msarp}).
\elem

\bpf
By Lemma~\ref{specgroups}, $G_{0}$ must be one of the groups $[3,3]^{*}$, $[3,3]$ or $[3,4]$.
In each case, $G_{0}$ is associated with standard generators of $\widehat{G}_{0} :=
\scl{S_{1}',-S_{2}}$.  We now investigate the possible generators $S_{1},S_{2}$ of $G$.

We begin with the group $[3,3]^{*}$. Now $S_{1}',-S_{2}$ are standard generators of $[3,3]^+$.
Once $S_{2}$ has been selected, there are three admissible choices for $S_{1}'$ such that 
$S_{1}'S_{2}$ has period $2$; if $S_{1}'$ is one of them, then the two others are 
$S_{2}^{-1}S_{1}'S_{2}$ and $S_{2}^{-2}S_{1}'S_{2}^2$. This follows directly from 
\cite[Lemma~5.1]{schu}, applied to the rotatory reflections $S_{2}$ and $-S_{1}'$. There is
also the further possibility of reversing the orientation of the generator $S_{2}$ and
replacing it by its inverse $S_{2}^{-1}$. However, this simply amounts to replacing the
polyhedron $P$ by its enantiomorphic image; the latter is the same underlying polyhedron, with
the same group, but with the new distinguished generators $S_{1}S_{2}^{2},S_{2}^{-1}$
(determined by a new, adjacent base flag); see Section~\ref{bano} for more details. 

Thus, without loss of generality we may assume that 
\beq 
\label{hesp33st}
\bry{rccl}
S_{1}'\colon & x & \mapsto &  (\xi_{2},-\xi_{3},-\xi_{1}) ,\\
S_{2}\colon  & x & \mapsto & -(\xi_{3},\xi_{1},\xi_{2}),  \\
\ery  
\eeq
described in terms of $x = (\xi_{1},\xi_{2},\xi_{3})$, so that $T'=S_{1}'S_{2}$ is the
reflection in the $\xi_{1}\xi_{3}$-plane. (Note that $-S_{1}',S_{2}$ are the generators for the
special group of the group $G(a,b)$ of \cite[(5.1)]{schu}.) Then $T=S_{1}S_{2}$ must be a
reflection in a plane parallel to the $\xi_{1}\xi_{3}$-plane, and since $T$ interchanges the
vertices in the base edge $F_1$, this plane must coincide with the perpendicular bisector of $F_1$.
In particular, the vertex $v:=oT$ of $F_1$ distinct from $o$ must be of the form $v = (0,a,0)$ for
some non-zero parameter $a$. It follows that $G$ is generated by 
\beq 
\label{he33st}
\bry{rccl}
S_{2}\colon & x & \mapsto & -(\xi_{3},\xi_{1},\xi_{2}),  \\
T\colon & x & \mapsto &  (\xi_{1},-\xi_{2},\xi_{3}) + (0,a,0).  \\
\ery  
\eeq
Since $S_{1} = T S_{2}^{-1}$, we then have
\beq
\label{hesone33st}
\bry{rccl}
S_{1}\colon & x & \mapsto & (\xi_{2},-\xi_{3},-\xi_{1})  + (-a,0,0) ,
\ery 
\eeq
so in particular,
\beq
\label{hson33st}
\bry{rccl}
S_{1}^3\colon & x & \mapsto & (\xi_{1},\xi_{2},\xi_{3})  + a(-1,-1,1) ,
\ery 
\eeq
which is a non-trivial translation. Hence, $S_{1}$ must be a twist and $P$ must have helical
faces (spiraling over an equilateral triangle). Note that $P$ is of type $\{\infty,6\}$.

We now employ Lemma~\ref{regcrit} to show that $P$ must be regular. It is straightforward 
to check that the reflection $R$ (in the plane $\xi_{1}=\xi_{3}$) given by
\beq
\label{refone33st}
\bry{rccl}
R\colon & x & \mapsto & (\xi_{3},\xi_{2},\xi_{1}) 
\ery 
\eeq
satisfies the conditions required in Lemma~\ref{regcrit} with 
\[ F_{0}=o, \;\; F_{1}=\{o,v\}, \;\; F_{2}=o\scl{S_{1}}, \;\; F_{2}'=F_{2}T .\] 
In fact, observe that conjugation by $R$ transforms the pair $S_{1},S_{2}$ to the new generating
pair $S_{1}S_{2}^{2},S_{2}^{-1}$, so in particular $F_{2} = o\scl{S_{1}}$ is mapped to 
$F_{2}' = o\scl{S_{1}S_{2}^{2}}$. Hence, $P$ is regular and $R=R_2$. Inspection of the list of twelve
regular apeirohedra (with an irreducible group) shows that $P = \{\infty,6\}_{4,4}$ (see
\cite[p.224]{msarp}), the Petrial of the Petrie-Coxeter polyhedron $\{4,6\hole 4\}$ (see
\cite{crsp}), with Petrie polygons of length $4$ (indicated by the first subscript in the symbol for
$P$); note here that $S_{1}^{2}S_{2}^{2}\;(= {(R_{0}R_{1}R_{2})}^{2})$ has period $2$.

Next we investigate the possible generators $S_{1},S_{2}$ of $G$ when $G_{0}=[3,4]$. Now
$S_{1}',-S_{2}$ are standard generators of $[3,4]^+$, viewed as the rotation group of a cube
$\{4,3\}$. As in the previous case, once $S_{2}$ has been selected, there are precisely three
admissible choices for $S_{1}'$ such that $S_{1}'S_{2}$ has period $2$; if $S_{1}'$ is one of them, 
then the two others are $S_{2}^{-1}S_{1}'S_{2}$ and $S_{2}^{-2}S_{1}'S_{2}^2$. This
follows from \cite[Lemma~6.1]{schu}, applied to the rotatory reflections $S_{2}$ and $-S_{1}'$.
Moreover, the same remark as above about replacing $S_2$ by $S_{2}^{-1}$ applies in this
context.

Thus we may take as generators
\beq 
\label{hesp34}
\bry{rccl}
S_{1}'\colon & x & \mapsto &  (\xi_{1},-\xi_{3},\xi_{2}), \\
S_{2}\colon  & x & \mapsto & -(\xi_{3},\xi_{1},\xi_{2}),  \\
\ery  
\eeq
so that $T'$ is the reflection in the plane $\xi_{1}+\xi_{2} = 0$. (Now $-S_{1}',S_{2}$ are the
generators for the special group of the group $H(c,d)$ of \cite[(6.1)]{schu}.) The mirror of $T$
once again is the perpendicular bisector of $F_1$ and is parallel to the mirror of $T'$. Hence, if
$v:=oT$, then $v=(a,a,0)$ for some non-zero parameter $a$. It follows that $G$ is generated by 
\beq 
\label{he34}
\bry{rccl}
S_{2}\colon & x & \mapsto & -(\xi_{3},\xi_{1},\xi_{2}),  \\
T\colon     & x & \mapsto &  (-\xi_{2},-\xi_{1},\xi_{3})  + (a,a,0).  \\
\ery  
\eeq
The element $S_{1} = T S_{2}^{-1}$ is given by
\beq
\label{hesone34}
\bry{rccl}
S_{1}\colon & x & \mapsto & (\xi_{1},-\xi_{3},\xi_{2})  + (-a,0,-a) ,
\ery 
\eeq
so in particular, 
\beq
\label{hson34}
\bry{rccl}
S_{1}^4\colon & x & \mapsto & (\xi_{1},\xi_{2},\xi_{3})  + a(-4,0,0) ,
\ery 
\eeq
which is again a non-trivial translation. Hence, $S_{1}$ must be a twist and $P$ must have helical
faces (spiraling over a square).

As in the previous case we can prove that $P$ is regular. The element $R\;( = R_2)$ is the reflection (in the plane $\xi_{1}=\xi_{2}$) given by
\beq
\label{refone34}
\bry{rccl}
R\colon & x & \mapsto & (\xi_{2},\xi_{1},\xi_{3}) .
\ery 
\eeq
In this case, $P = \{\infty,6\}_{6,3}$, the Petrial of the Petrie-Coxeter polyhedron $\{6,6\hole
3\}$ , with Petrie polygons of length $6$ (again indicated by the first subscript in the symbol for
$P$); now $S_{1}^{2}S_{2}^{2}$ has period $3$.

Finally, when $G_{0} = [3,3]\;(\cong S_{4})$, we must take a rotation $S_{1}'$ of period $3$ and
a rotatory reflection $S_{2}$ of period $4$ as generators of $G_{0}$. Once $S_{2}$ has been
selected, there now are exactly four admissible choices for $S_{1}'$; if $S_{1}'$ is one of them,
then the three others are $S_{2}^{-j}S_{1}'S_{2}^{j}$ for $j=1,2,3$. Moreover, the same remark
as before about replacing $S_2$ by $S_{2}^{-1}$ applies.

Thus, with
\beq 
\label{heg33}
\bry{rccl}
S_{1}'\colon & x & \mapsto &  (\xi_{2},\xi_{3},\xi_{1}), \\
S_{2}\colon  & x & \mapsto & (-\xi_{1},-\xi_{3},\xi_{2}),  \\
\ery  
\eeq
we obtain the same $v$, $T'$ and $T$ as in the previous case. In particular we have
\beq
\label{hesone33}
\bry{rccl}
S_{1}  \colon & x & \mapsto & (\xi_{2},\xi_{3},\xi_{1})  + (-a,0,-a) ,\\
S_{1}^3\colon & x & \mapsto & (\xi_{1},\xi_{2},\xi_{3})  + a(-2,-2,-2) .
\ery 
\eeq
Once again, $P$ is regular, by Lemma~\ref{regcrit}. Now $R\;(=R_{2})$ is given by the reflection 
(in the plane $\xi_{3}=0$)
\beq
\label{refone33}
\bry{rccl}
R\colon & x & \mapsto & (\xi_{1},\xi_{2},-\xi_{3}) ,
\ery 
\eeq
and $P = \{\infty,4\}_{6,4}$, the Petrial of $\{6,4\hole 4\}$, with dimension vector $(1,1,2)$
(see \cite[p.224]{msarp}).
\epf

The three polyhedra occurring in Lemma~\ref{twi} are precisely the pure regular polyhedra with 
dimension vector $(1,1,2)$ (see \cite[p.225]{msarp}). Note that the values $1$, $1$, $2$ for
the dimensions of the mirrors for their generating reflections $R_{0},R_{1},R_{2}$ accord with
the geometric type of the generators of $G$, namely a twist for $S_{1}\; ( = R_{0}R_{1})$ and a
rotatory reflection for $S_{2}\; (=R_{1}R_{2})$.

\blem
\label{sonerot}
Let $P$ be a discrete chiral or regular apeirohedron in $\BE^3$, and let the corresponding group 
$G = \scl{S_1,S_2}$ be irreducible. Then the generators cannot be such that $S_{1}$ is a
rotatory reflection and $S_{2}$ a rotation.
\elem

\bpf
Suppose that $S_{1}$ (and hence $S_{1}'$) is a rotatory reflection and $S_{2}$ a rotation. Then
$G_{0} = [3,3]^{*}$ or $[3,4]$, by the proof of Lemma~\ref{specgroups}. We treat these cases
separately and appeal to the discussion of the groups occurring in the proof of Lemma~\ref{twi}. 
Again, $o$ is the base vertex.

If $G_{0} = [3,3]^{*}$, then we may take
\beq 
\label{sp33star}
\bry{rccl}
S_{1}'\colon & x & \mapsto &  -(\xi_{2},\xi_{3},\xi_{1}) ,\\
S_{2}\colon  & x & \mapsto & (-\xi_{3},\xi_{1},-\xi_{2}),  \\
\ery  
\eeq
which are the inverses of the generators $S_{2}$ and $S_{1}'$ of \eqref{hesp33st}, respectively, of 
periods $6$ and $3$. Then $T'$, $T$ and $v$ of \eqref{hesp33st} remain unchanged, because $T$
and $T'$ have period $2$. In particular, $G$ is generated by
\beq 
\label{cu33st}
\bry{rccl}
S_{2}\colon & x & \mapsto & (-\xi_{3},\xi_{1},-\xi_{2}),  \\
T\colon & x & \mapsto &  (\xi_{1},-\xi_{2},\xi_{3}) + (0,a,0),  \\
\ery  
\eeq
and $S_{1} = T S_{2}^{-1}$ is given by
\beq
\label{cu33stpow}
\bry{rccl}
S_{1}\colon & x & \mapsto & (-\xi_{2},-\xi_{3},-\xi_{1})  + (a,0,0) .
\ery 
\eeq
But then the whole group $G$ fixes the point $z := \half (a,a,-a)$, and so $G$ must be a finite 
group, contrary to our assumption that $P$ is infinite. (This finite group consists of the
symmetries of the cube of edge length $a$ centered at $z$, and the associated polyhedron is the
Petrial $\{6,3\}_4$ of $\{4,3\}$; see \cite[p.192,218]{msarp}.)

In a similar fashion we deal with the case $G_{0}=[3,4]$. We take 
\beq 
\label{hesp34a}
\bry{rccl}
S_{1}'\colon & x & \mapsto &  -(\xi_{2},\xi_{3},\xi_{1}),\\
S_{2}\colon  & x & \mapsto & (\xi_{1},\xi_{3},-\xi_{2}),  \\
\ery  
\eeq
which are the inverses of the generators $S_{2}$ and $S_{1}'$ of \eqref{hesp34}, respectively, of 
periods $6$ and $4$. Then $T'$, $T$ and $v$ of \eqref{hesp34} remain unchanged and $G$ is
generated by any two of
\beq 
\label{hesp34oct}
\bry{rccl}
S_{1}\colon & x & \mapsto & -(\xi_{2},\xi_{3},\xi_{1})  + (a,0,a), \\
S_{2}\colon & x & \mapsto & (\xi_{1},\xi_{3},-\xi_{2}),  \\
T\colon     & x & \mapsto &  (-\xi_{2},-\xi_{1},\xi_{3})  + (a,a,0).  \\
\ery  
\eeq
Then $G$ fixes $z' := (a,0,0)$, so it must again be a finite group, contrary to our assumptions on 
$P$. (Now this finite group consists of the symmetries of the octahedron with vertices $z'
\pm a\,e_{1}, z' \pm a\,e_{2}, z' \pm a\,e_{3}$, where $e_{1},e_{2},e_{3}$ is the standard
basis of $\BE^3$; the associated polyhedron is $\{6,4\}_3$, the Petrial of $\{3,4\}$.)
\epf

The above analysis shows that, for a discrete chiral apeirohedron $P$ with an irreducible group
$G$, either the generators $S_{1},S_{2}$ of $G$ must be rotatory reflections, or $S_{1}$ must be a
twist and $S_{2}$ a rotation. 

If $P$ has finite faces, then both generators must be rotatory reflections of finite period, and 
$P$ must be among the polyhedra enumerated in Part~I (\cite{schu}). In particular, this completes
the proof of the following theorem, at least in the irreducible case; the reducible case is settled
by Lemma~\ref{chirreds}.

\bthm
\label{comple}
There are no discrete chiral polyhedra with finite faces in $\BE^3$ other than those (with
skew faces and skew vertex-figures) enumerated in Part~I.
\ethm

For a chiral apeirohedron $P$ with infinite, helical faces, $S_{1}$ must be a twist and $S_{2}$ a 
rotation. In particular, $G_{0}$ is generated by the rotations $S_{1}',S_{2}$ (with distinct
axes), so $G_{0}$ must be a finite group of rotations in $\BE^3$. Bearing in mind that
$S_{1}'S_{2}$ must have period $2$, we immediately arrive at

\blem
\label{chirspecgr}
Let $P$ be a discrete chiral or regular apeirohedron in $\BE^3$ with helical faces, and let
$o$ be its base vertex. Let the corresponding group $G = \scl{S_1,S_2}$ act irreducibly, with 
special group $G_{0} = \scl{S_1',S_2}$, and let $S_{2}$ be a rotation. Then $G_{0} = [3,3]^{+}$
or $[3,4]^{+}$. Moreover, if $p'$ denotes the period of $S_{1}'$ and $P$ is of type
$\{\infty,q\}$, then $(p',q) = (3,3)$, $(3,4)$ or $(4,3)$, and the faces of $P$ are helices
over equilateral triangles or squares according as $p'=3$ or $4$.
\elem

In the next three sections, we shall enumerate the discrete chiral apeirohedra with helical
faces.  They occur in three real-valued two-parameter families, each of type $\{\infty,3\}$
or $\{\infty,4\}$. In each case, $T=S_{1}S_{2}$ is a half-turn whose rotation axis is perpendicular
to the base edge of $P$; in particular this accounts for the two independent parameters describing
the base vertex $v = oT$. By contrast, observe that, in the situations discussed in Lemmas~\ref{twi}
and \ref{sonerot}, the element $T$ is a plane reflection, so that we have only one independent
parameter; this accords with $P$ being regular in these cases, and being determined up to
similarity.

Finally, then, we must settle the reducible case. Recall that there are twelve non-planar regular
polyhedra in $\BE^3$ with a reducible group, namely the blended polyhedra 
\beq
\label{regreds}
\{p,q\} \blend \seg,\;\;  \{p,q\} \blend \{\infty\},\;\;
\{p,q\}^{\pi} \blend \seg,\;\; \{p,q\}^{\pi} \blend \{\infty\} ,
\eeq
where the first component is either a regular tessellation $\{p,q\} = \{4,4\}$, $\{3,6\}$ or
$\{6,3\}$, or the Petrial $\{p,q\}^{\pi}$ of $\{p,q\}$, and where the second component is either a
segment $\seg$ or an apeirogon $\{\infty\}$ (see \cite[p.221-223]{msarp}). The following lemma says
that reducible groups can be eliminated from further discussion.

\blem
\label{chirreds}
A discrete chiral apeirohedron in $\BE^3$ must have an irreducible symmetry group.
\elem

\bpf
It is not possible for a discrete chiral polyhedron in $\BE^3$ to be planar (see \cite[Thm.
3.2]{schu} or \cite{mc}), so it suffices to consider non-planar polyhedra. Let $P$ be a discrete,
non-planar polyhedron, and let $P$ be chiral or regular. Suppose the corresponding group $G =
\scl{S_{1},S_{2}}$ is reducible. We show that $P$ is regular.

Since $G$ is reducible, there exists a plane $L$ through $o$ such that $G$ permutes the translates
of $L$ as well as the translates of $L^{\perp}$, the orthogonal complement of $L$ at $o$. In
particular, since $oS_{2}=o$, the subspaces $L$ and $L^{\perp}$ are invariant under $S_{2}$, and
$S_{2}$ is either a rotation (about $L^{\perp}$) or a rotatory reflection (with a rotation about
$L^{\perp}$ as its rotation component). Moreover, the vertex $v$ of the base edge $F_{1}$ distinct
from $o$ cannot lie in $L$; otherwise, the vertex-figure at $o$, and hence at every vertex, would
lie in $L$, making $P$ planar. We shall see that the geometric type of $S_{2}$ determines the second
component of the blend. 
 
Now suppose $S_{2}$ is a rotation. Let $L'$ denote the plane parallel to $L$ through $v$. Then $P$
has all its vertices in $L$ or $L'$, and the vertex-figure at a vertex in one plane is a flat
(finite or infinite) polygon contained in the other plane; in particular, the edge-graph of $P$ is
bipartite. 

If the faces of $P$ are finite (skew) polygons, then the generator $S_{1}$ is a rotatory
reflection, so $T$ can only be the point reflection in $\tfrac{1}{2}v$; in particular, 
$S_{1}' = T' S_{2}'^{-1} = (-I)S_{2}'^{-1}$. The projection of $P$ on $L$ is discrete, and since
$S_{2}$ and the rotation component of $S_{1}'$ must be rotations of periods $2$, $3$, $4$ or
$6$, it is not hard to see (using basic properties of Wythoff's construction) that it must be a
tessellation $\{p,q\}$; hence $P = \{p,q\} \blend \seg$. 

Similarly, if the faces are infinite (zig-zags), then $S_{1}$ is a half-twist about a line
parallel to $L$ (a half-turn about the axis of the base face, followed by a translation along this
axis), so $T$ must be a half-turn. Now the projection on $L$ is the Petrial $\{p,q\}^{\pi}$, and 
$P = \{p,q\}^{\pi} \blend \seg$. 

Now let $S_{2}$ be a rotatory reflection. Let $L_{k}$ denote the plane parallel to $L$ through
the point $kv$, for $k \in \BZ$. Now the vertices adjacent to $o\;(\in L_{0})$ alternate between
$L_{-1}$ and $L_{1}$. It follows that $P$ has all its vertices in the union of the planes
$L_{k}$, and that the vertices adjacent to a vertex in one of the planes alternate between the
adjacent planes. Hence the faces are helices or zig-zags, each with one vertex on each plane
$L_{k}$. 

If the faces of $P$ are helices, then $S_{1}$ is a twist about an axis (of the base face) parallel
to $L^{\perp}$, which is determined by a rotation of finite period $p'$ (say) in $L$ and a
translation along this axis. The projection of $P$ on $L$ is again discrete. This can be seen as
follows. Let $D$ be a circular disc in $L$ centered at $o$, and let $\widehat{D}$ be a right cylinder
over $D$ with height $|t|$, where $t$ is the translation vector of the translation $S_{1}^{p'}$.
Now, if $x$ is a vertex of $P$ projecting on $D$, then $\widehat{x} := x + nt = xS_{1}^{np'}$, for
some $n\in\BZ$, is a vertex of $P$ in $\widehat{D}$ which projects on the same point as $x$. Thus the
vertices in $\widehat{D}$ represent all the vertices projecting on $D$. However, since $P$ is
discrete, it only has finitely many vertices in $\widehat{D}$, so its projection only has finitely
many points in $D$. This proves that the projection of $P$ on $L$ is discrete. As before we then
establish that this projection is a tessellation $\{p,q\}$ in $L$. Hence $P = \{p,q\} \blend \{
\infty \}$.

Finally, if the faces of $P$ are zig-zags, then $S_{1}$ must be a glide reflection whose
respective reflection mirror is perpendicular to $L$ and contains the axis of the base face.
Moreover, $T$ is a half-turn about a line parallel to $L$. The projection of $P$ on $L$ is again
discrete, for similar reasons as in the previous case (now $S_{1}^{2}$ is a translation). Now this
projection is a Petrial $\{p,q\}^{\pi}$, and $P = \{p,q\}^{\pi} \blend \{\infty\}$. This concludes
the proof.
\epf

\section{Type $\{\infty,3\}$, with helical faces over triangles}
\label{helfacone}

In this section we enumerate the chiral apeirohedra of type $\{\infty,3\}$ whose faces are
helices over triangles. They occur when $G_{0}=[3,3]^{+}$ in Lemma~\ref{chirspecgr}. Throughout we
may assume that the symmetry group $G$ acts irreducibly; we saw in Lemma~\ref{chirreds} that the
reducible case adds no further polyhedra to the list.

Each polyhedron is obtained from a suitable group $G$ by Wythoff's construction with initial (or 
base) vertex $F_{0}:=o$. For a polyhedron $P$ of type $\{\infty,3\}$ with helical faces over
triangles, we must necessarily begin with a group $G=\scl{S_{1},S_{2}}$ whose special group is
$G_{0}=[3,3]^{+}$. Then $S_{1}',S_{2}$ are standard generators of $[3,3]^{+}$, and each is a
rotation of period $3$. As in similar situations discussed earlier, once $S_{2}$ has been selected,
there are three equivalent choices for $S_{1}'$ such that $S_{1}'S_{2}$ has period $2$; if $S_{1}'$
is one of them, then the other two are $S_{2}^{-1}S_{1}'S_{2}$ and $S_{2}^{-2}S_{1}'S_{2}^{2}$.
Moreover, as mentioned before, if we substitute $S_{2}$ by its inverse, we arrive at (the
enantiomorphic image of) the same underlying polyhedron, so we need not investigate this case
separately.

Thus we may take 
\beq 
\label{heltrisg}
\bry{rccl}
S_{1}'\colon & x & \mapsto &  (-\xi_{3},-\xi_{1},\xi_{2}) ,\\
S_{2}\colon  & x & \mapsto &  (\xi_{2},\xi_{3},\xi_{1}).  \\
\ery  
\eeq
Then $T'=S_{1}'S_{2}$ is the half-turn about the $\xi_{2}$-axis. Since $T=S_{1}S_{2}$ must be a
half-turn interchanging the two vertices in the base edge $F_{1}$, its rotation axis must
necessarily be perpendicular to $F_{1}$ and pass through the midpoint of $F_{1}$. But this axis
is parallel to the $\xi_{2}$-axis, so $F_{1}$ must lie in the $\xi_{1}\xi_{3}$-plane.
Hence, $v:=oT = (a,0,b)$ for some real parameters $a$ and $b$, not both zero, so that we have
$G=G_{1}(a,b)$, the group generated by
\beq 
\label{heltrigr}
\bry{rccl}
S_{2}\colon  & x & \mapsto &  (\xi_{2},\xi_{3},\xi_{1}),  \\
T    \colon  & x & \mapsto &  (-\xi_{1},\xi_{2},-\xi_{3}) + (a,0,b).
\ery  
\eeq
From $S_{1}=TS_{2}^{-1}$ we then obtain
\beq 
\label{heltrisone}
\bry{rccl}
S_{1}    \colon & x & \mapsto &  (-\xi_{3},-\xi_{1},\xi_{2}) + (b,a,0), \\
S_{1}^{3}\colon & x & \mapsto &  (\xi_{1},\xi_{2},\xi_{3}) + (b-a)(1,-1,-1). \\ 
\ery  
\eeq
In particular, $S_{1}^{3}$ is the translation by $(b-a)(1,-1,-1)$, and is trivial if $b=a$. Hence
the case $b=a$ is special and yields the finite group
\[ G_{1}(a,a) \cong [3,3]^{+}  \]
(with fixed point $\half (a,a,a)$); the corresponding polyhedron is a regular tetrahedron. It is
often convenient to exclude this case from the discussion.

When necessary, we indicate the parameters $a$ and $b$ more explicitly and write $T=T(a,b)$, 
$S_{1} = S_{1}(a,b)$, and so on; the generator $S_{2}$ does not actually depend on $a,b$.

The polyhedron $P = P_{1}(a,b)$ is obtained from $G = G_{1}(a,b)$ by Wythoff's construction with
base vertex $F_{0}=o$. The base edge $F_{1}$ has vertices $o$ and $v$. If we identify faces of $P$
with their vertex-sets, then a helical face $F$ (over a triangle) is denoted as follows:  if
$x,y,z$ are consecutive vertices of $F$ and $t$ is the translation vector associated with $F$,
then 
\[  F = \{x,y,z\} + {\Bbb Z}\cdot t , \]
with the understanding that the order in which the vertices of $F$ occur is given by
\[ \ldots, x\!-\!2t,y\!-\!2t,z\!-\!2t,x\!-\!t,y\!-\!t,z\!-\!t,x,y,z,x\!+
\!t,y\!+\!t,z\!+\!t,x\!+\!2t,y\!+\!2t,z\!+\!2t, \ldots \;.\]
The same notation, with $t=o$, is also adopted in the special case $b=a$; then $F = \{x,y,z\}$, so
we usually drop the term ${\Bbb Z}\cdot t$. In particular, for the base face we have
\beq
\label{facetspone}
F_{2} = o\scl{S_{1}} = \{(a,0,b),(0,0,0),(b,a,0)\} + {\Bbb Z}\!\cdot\!c(1,-1,-1) ,
\eeq
where $c:=b-a$. The remaining vertices, edges and faces of $P$ are the images of $F_{0}$,
$F_{1}$ and $F_{2}$ under $G$. The vertices adjacent to $o$ are 
\[ v = (a,0,b),\; vS_{2} = (0,b,a),\; vS_{2}^{2} = (b,a,0) . \]

If the polyhedra $P_{1}(a,b)$ are only considered up to similarity, then the two parameters $a,b$ 
can be reduced to a single parameter. This is similar to what we observed for the polyhedra in
\cite{schu}. In fact, if $s$ is a non-zero scalar and $R_{s}:=s I$, then conjugation of the
generators $S_{2},T(a,b)$ for $G_{1}(a,b)$ by $R_s$ gives the generators $S_{2},T(sa,sb)$ for
$G_{1}(sa,sb)$.  It follows that 
\beq
\label{simpone}
P_{1}(sa,sb) = P_{1}(a,b)R_{s} =: s\,P_{1}(a,b) ,
\eeq
so $P_{1}(sa,sb)$ and $P_{1}(a,b)$ are similar. In particular, $P_{1}(a,b)$ is similar to 
$P_{1}(1,\tfrac{b}{a})$ or $P_{1}(0,1)$ according as $a\neq 0$ or $a=0$. 

Before we move on, observe that the polyhedra $P_{1}(a,b)$ of type $\{\infty,3\}$ are related
to the chiral (or regular) polyhedra $P(a,b)$ of type $\{6,6\}$ described in \cite[\S 5]{schu}
(their notation has no suffix). These polyhedra have finite skew hexagonal faces and vertex-figures.
The relationship is established by the facetting operation $\p_{2}$ defined in \eqref{dualop}.
Recall from \cite[(5.1),(5.2)]{schu} that the generators $\widehat{S}_{1},\widehat{S}_{2}$ (say) and
the  corresponding element $\widehat{T} := \widehat{S}_{1}\widehat{S}_{2}$ of the symmetry group
$G(a,b)$ of $P(a,b)$ are given by
\beq 
\label{dhgroup}
\bry{rccl}
\widehat{S}_{1}\colon & x & \mapsto & (-\xi_{2},\xi_{3},\xi_{1}) + (0,-b,-a),\\
\widehat{S}_{2}\colon & x & \mapsto & -(\xi_{3},\xi_{1},\xi_{2}),  \\
\widehat{T}\colon & x & \mapsto &  (-\xi_{1},\xi_{2},-\xi_{3}) + (a,0,b).  
\ery  
\eeq
We now have the following 

\blem
\label{p2pab}
$P_{1}(a,b) = P(a,b)^{\p_2}$, for all real parameters $a$ and $b$. 
\elem

\bpf
In the present context, $\p_{2}$ takes the form
\beq
\label{geomp2es}
\bry{rccl}
\p_{2}\colon & (\widehat{S}_{1},\widehat{S}_{2}) & \mapsto &
(\widehat{S}_{1}\widehat{S}_{2}^{-1},\widehat{S}_{2}^{2}).
\ery 
\eeq
Now observe that, when $\p_{2}$ is expressed in terms of the generators 
$\widehat{S}_{2},\widehat{T}$ of $G(a,b)$, we obtain
\beq
\label{geomp2st}
\bry{rccl}
\p_{2}\colon & (\widehat{S}_{2},\widehat{T}) & \mapsto &
(\widehat{S}_{2}^{2},\widehat{T}) = (S_{2},T). 
\ery 
\eeq
Hence $P_{1}(a,b) = P(a,b)^{\p_2}$ and $G_{1}(a,b)$ is a subgroup of $G(a,b)$.
\epf

Note that the above relationship between $P(a,b)$ and $P_{1}(a,b)$ holds without any 
restriction on the parameters.  On the other hand, $P(a,b)$ is known to be discrete only when
$a$ and $b$ can be scaled to a pair of relatively prime integers, whereas there is no such condition
for $P_{1}(a,b)$.  Thus we may (in most cases, do) have the interesting phenomenon that the derived
polyhedron $P_{1}(a,b)$ is discrete but the original polyhedron $P(a,b)$ is non-discrete. We can conclude from this that then there are infinitely many copies of $P_{1}(a,b)$ inscribed in $P(a,b)$.

We already note here two special cases. If $b=-a$, then the polyhedra are similar to those obtained 
for $(a,b)=(1,-1)$, and are regular. In particular, since $P(1,-1) = \{6,6\}_{4}$ (see
\cite[Thm.5.17]{schu}) and $\{6,6\}_{4}^{\p_2} =  \{\infty,3\}^{(a)}$ (notation as in
\cite[p.224]{msarp}), we observe that
\beq 
\label{infa}
P_{1}(1,-1) = \{\infty,3\}^{(a)} .
\eeq
Similarly, the case $b=a$ reduces to $(a,b)=(1,1)$ and yields regular polyhedra. Now 
$P(1,1) = \{6,6 \hole 3\}$ and $\{6,6 \hole 3\}^{\p_2} = \{3,3\}$ (again, see \cite[Thm.5.17]{schu}
and \cite[p.224]{msarp}), so in particular
\beq 
\label{threeones}
P_{1}(1,1) = \{3,3\}. 
\eeq

It follows from our discussion in Section~\ref{bano} that the faces of the polyhedron $P_{1}(a,b)$ 
are among the holes of $P(a,b)$. These holes are helices over triangles, or triangles, according as 
$b \neq a$ or $b=a$.

We now determine the translation subgroup $T(G_{1}(a,b))$ of $G_{1}(a,b)$. When convenient, we 
identify a translation with its translation vector. It is convenient to use the following
notation for lattices. Let $s$ be a non-zero real number, let $k=1$, $2$ or $3$, and let $\s :
= (s^k,0^{3-k})$, the vector with $k$ components $s$ and $3-k$ components $0$. Following
\cite[p.166]{msarp}, we write $\La_{\s}$ for the sublattice of the cubic lattice $\BZ^3$
generated by $\s$ and its images under permutation and changes of sign of coordinates. Then
\[  \La_{\s} = s\La_{(1^{k},0^{3-k})}.  \]
Note that $\La_{(1,0,0)} = \BZ^{3}$. The lattice $\La_{(1,1,0)}$ is the {\em face-centered cubic 
lattice\/} (the root lattice $D_{3}$) consisting of all integral vectors whose coordinate
sum is even, and $\La_{(1,1,1)}$ is the {\em body-centered cubic\/} lattice (see \cite{cosl}). 

We know from \eqref{heltrisone} that $S_{1}^3$ is the translation by $c(1,-1,-1)$, where again $c :=
b-a$. Hence, if $R$ is any element of $G=G_{1}(a,b)$ and $R'$ its image in the special group
$G_{0}$, then $R^{-1}S_{1}^{3}R$ is the translation by $c(1,-1,-1)R'$. Since $G_{0} \cong
[3,3]^{+}$, we then have translations by the vectors $c(1,1,1)$, $c(1,-1,-1)$, $c(-1,1,-1)$,
$c(-1,-1,1)$ in $G$, and hence also by their integral linear combinations. It follows that
\[ c \La_{(1,1,1)} \leq T(G) .\]
(The left hand side is trivial if $b=a$.)  In fact, we have the following stronger result. 

\blem
\label{ponetrans}
The translation subgroup of $G_{1}(a,b)$ is given by $T(G_{1}(a,b)) = c \La_{(1,1,1)}$, with 
$c = b-a$.
\elem

\bpf
Suppose $b\neq a$. Since $G = G_{1}(a,b) = \scl{S_{2},T}$, we have $G = N \cdot \scl{S_{2}}$ (as a
product of subgroups), where $N := \scl{T_{1},T_{2},T_{3}}$ with $T_{2}:=T$,
$T_{1} := S_{2}^{-1}TS_{2}$ and $T_{3} := S_{2}^{-2}TS_{2}^{2}$. Bearing in mind that the generator
$S_{2}$ for $P_{1}(a,b)$ is just the square of the generator $\widehat{S}_{2}$ for $P(a,b)$ (see
\eqref{geomp2st}), we see that $T_{1},T_{2},T_{3}$ are the half-turns occurring in
\cite[(5.1),(5.6)]{schu} with the same notation, namely
\beq 
\label{t123}
\bry{rccl}
T_{1}\colon & x & \mapsto &  (\xi_{1},-\xi_{2},-\xi_{3}) + (0,b,a),  \\
T_{2}\colon & x & \mapsto &  (-\xi_{1},\xi_{2},-\xi_{3}) + (a,0,b), \\
T_{3}\colon & x & \mapsto &  (-\xi_{1},-\xi_{2},\xi_{3}) + (b,a,0). 
\ery  
\eeq
Observe that a translation in $G$ must necessarily belong to $N$; in fact, in the special 
group, the images of a translation in $G$ or of an element in $N$ must involve an even
number of sign changes but no permutation of coordinates, whereas a non-trivial element in
$\scl{S_2}$ involves a non-trivial permutation of coordinates.

It is straightforward to check that each product $T_{i}T_{j}T_{k}$, with $i,j,k$ distinct, is a 
translation in $G$. The translation vectors of $T_{1}T_{3}T_{2}$, $T_{2}T_{1}T_{3}$ and
$T_{3}T_{2}T_{1}$ are $c(-1,-1,1)$, $c(1,-1,-1)$ or $c(-1,1,-1)$, respectively; the three
other products yield the negatives of these vectors. Hence the six products are translations in
$c\La_{(1,1,1)}$.

Now suppose we have an element $R := T_{j_{1}}T_{j_{2}} \ldots T_{j_{n}}$ of $N$, with $j_{m}=1$, 
$2$ or $3$ for each $m$. For a non-trivial translation we must have $n \geq 3$. If $n \geq 3$
and $j_{n-2},j_{n-1},j_{n}$ are distinct, then we can reduce $R$ modulo $c\La_{(1,1,1)}$ to a
product with $n-3$ terms $T_{j_m}$ (by splitting off the product of the last three
generators), and then proceed by induction. If $n \geq 3$ and $j_{n-2},j_{n-1},j_{n}$ are not
mutually distinct, then two cases are possible. If $j_{n-1}=j_{n}$, then we directly
eliminate the last two generators. If $j_{n-1} \neq j_{n}$, we first insert the trivial
element $T_{k}T_{k}$ with $k\neq j_{n-1},j_{n}$ into the product between $T_{j_{n-2}}$ and
$T_{j_{n-1}}$, and then split off the product of the last three terms as before. In any
case, by induction we can reduce $R$ modulo $c\La_{(1,1,1)}$ to a product with at most two
terms, and hence to the empty product if $R$ is a translation. 

Now it follows that each translation in $G$ necessarily belongs to $c \La_{(1,1,1)}$. 
\epf

Note that Lemma~\ref{ponetrans} can be rephrased as
\beq
\label{specgone}
G_{1}(a,b)\slash c\La_{(1,1,1)}  \cong [3,3]^{+} \cong A_{4} ,
\eeq
where each term represents the special group. Observe also that the translation subgroup only 
depends on the single parameter $c = b-a$. 

\begin{figure}[hbt]
\centering
\begin{center}
\begin{picture}(200,190)
\put(10,0){
\begin{picture}(150,150)
\multiput(40,60)(40,0){2}{\circle*{6}}
\put(80,60){\circle{10}}
\multiput(90,93.33)(40,0){2}{\circle*{6}}
\multiput(25,55.5)(0,40){2}{\circle*{6}}
\multiput(145,54.5)(0,42){2}{\circle*{6}}
\put(145,96.5){\circle{10}}
\multiput(76,10.5)(0,120){2}{\circle*{6}}
\put(76,130.5){\circle{10}}
\multiput(90.5,22.2)(0,120){2}{\circle*{6}}
\end{picture}}
\multiput(10,0)(0,60){3}{
\begin{picture}(120,30)
\thinlines
\multiput(0,0)(60,0){3}{\circle*{2}}
\multiput(25,16.66)(60,0){3}{\circle*{2}}
\multiput(50,33.3)(60,0){3}{\circle*{2}}
\multiput(0,0)(25,16.66){3}{\line(1,0){120}}
\multiput(0,0)(60,0){3}{\line(3,2){50}}
\end{picture}}
\put(10,0){
\begin{picture}(180,140)
\thinlines
\multiput(0,0)(60,0){3}{\line(0,1){120}}
\multiput(25,16.66)(60,0){3}{\line(0,1){120}}
\multiput(50,33.33)(60,0){3}{\line(0,1){120}}
\end{picture}}
\put(10,0){ 
\begin{picture}(180,140)
\thinlines
\put(85,76.66){\vector(1,0){100}}
\put(189,74.66){$\xi_{2}$}
\put(85,76.66){\vector(0,1){90}}
\put(82.5,172){$\xi_{3}$}
\put(85,76.66){\vector(-3,-2){50}}
\put(26,37){$\xi_{1}$}
\put(85,76.66){\circle*{2}}
\end{picture}}

\end{picture}
\caption{\it The set $V$ for $a=1$ and $b=3$. The points in the vertex-star $V_0$ for
$P_{1}(a,b)$ are circled. The cubes are drawn in for reference.}
\label{figvab}
\end{center}
\end{figure}
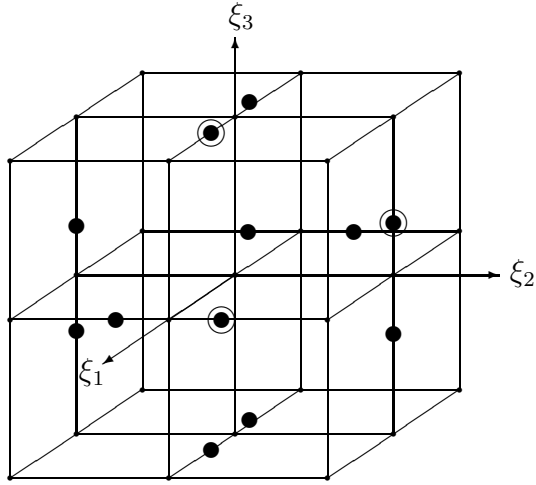

Next we determine the vertex-stars of $P_{1}(a,b)$. The {\em vertex-star\/} at a vertex $x$ of the 
polyhedron is a (cyclically ordered) set consisting of the vectors $y-x$, where $y$ runs
through the (cyclically ordered) vertices of the vertex-figure at $x$. In particular, the vertex-star
at
$o$ is given by
\beq
\label{v0star}
V_{0} := \{(a,0,b),(0,b,a),(b,a,0)\}. 
\eeq
If $R$ is any element of $G=G_{1}(a,b)$ and again $R'$ its image in the special group $G_{0}$, then
the vertex-star at the vertex $oR$ of $P_{1}(a,b)$ is given by 
\[ V_{0}R' = \{(a,0,b)R',(0,b,a)R',(b,a,0)R'\}. \]
It is easy to see that the union of all these vertex-stars comprises the set of (generally twelve) 
vectors
\beq
V := \{ (\pm a,0,\pm b), (\pm b,\pm a,0), (0,\pm b,\pm a) \},
\eeq
which is the same set as for the polyhedron $P(a,b)$ from which $P_{1}(a,b)$ is derived by 
$\p_{2}$ (see \cite[(5.9)]{schu} and Figure~\ref{figvab}). The points in $V$ are the vertices
of a convex $3$-polytope, which is combinatorially equivalent to an icosahedron if $a,b\neq 0$
and $a \neq \pm b$, or is a cuboctahedron if $a = \pm b \neq 0$, or an octahedron if $a=0$ or
$b=0$. 

Now we can find a full set of vertex-stars for $P_{1}(a,b)$. In fact, since $G$ acts transitively
on the vertices of $P_{1}(a,b)$, the vertex-stars of $P_{1}(a,b)$ are images of $V_0$ under the
special group $G_{0}\;(\cong A_{4})$, so $G_{0}$ must act transitively on them. But since $G_{0}$
contains $S_2$, the stabilizer of $V_0$ must be at least of order $3$, and hence the number of
vertex-stars cannot exceed $4$. In fact, there are precisely four vertex-stars (even if $a=0$ or
$b=0$), namely $V_{0}$ of \eqref{v0star} and
\beq 
\label{vonestars}
\bry{rccl}
V_{1} := V_{0}T_{1}' = \{(a,0,-b),(0,-b,-a),(b,-a,0)\}, \\
V_{2} := V_{0}T_{2}' = \{(-a,0,-b),(0,b,-a),(-b,a,0)\}, \\
V_{3} := V_{0}T_{3}' = \{(-a,0,b),(0,-b,a),(-b,-a,0)\}.
\ery  
\eeq
They occur as vertex-stars at the vertices $o$, $oT_{1} = (0,b,a)$, $oT_{2} = (a,0,b)$ and 
$oT_{3} = (b,a,0)$, respectively. Each point in $V$ lies in exactly one or exactly two vertex-stars according as $a,b\neq 0$, or $a=0$ or $b=0$.

Recall that the {\em edge-module\/} $\La$ of a polyhedron $P$, with $o$ among its vertices, is
defined by
\beq
\label{edmodule}
\La : = {\langle x-y \mid x,y \mbox{ adjacent vertices of } P \rangle}_{\BZ}
\eeq
(see \cite{schu}). Then $\La$ is the $\BZ$-module generated by the ``oriented" edges of $P$; 
alternatively, $\La$ is the $\BZ$-module generated by all the vertex-stars of $P$. If $V(P)$
denotes the vertex-set of $P$, then 
\[ V(P) \subset \La ,\]
because $P$ is connected. Note that the special group of $G(P)$ leaves $\La$ invariant. 

Now let $\La := \BZ[V]$ denote the $\BZ$-module spanned by the vectors in $V$. This is the
edge-module of $P_{1}(a,b)$ (and $P(a,b)$) and has rank at most $12$.  We proved in 
\cite[Lemma 5.2]{schu} that $\La$ is a (discrete) lattice if and only if $a$ or $b$ is zero or
$a$ and $b$ are rational multiples of each other, or, equivalently, if $a,b$ can be scaled to a
pair of integers which are relatively prime. On the other hand, if $\La$ is not a lattice, then
$\La$ is dense in $\BE^3$. This follows from general structure results for $\BZ$-modules
embedded in euclidean spaces (see \cite[p.265]{sen}); bear in mind that $\La$ cannot have a
discrete component, because such a component would have to be invariant under the
(irreducible) special group $[3,3]^{+}$. For example, if $a=1$ and $b=\tau$, the golden ratio,
then $\La$ is the (dense) $\BZ$-module spanned by the vertex-set $V$ of a regular icosahedron;
similar $\BZ$-modules have been studied in the context of icosahedral quasicrystals (see
\cite{chen}).

We now determine the vertex-set of $P_{1}(a,b)$. Since $P_{1}(a,b) = P(a,b)^{\p_2}$ and 
$P_{1}(a,b)$ and $P(a,b)$ share the same base vertex and base edge, the vertices and edges,
respectively, of $P_{1}(a,b)$ are among the vertices and edges of $P(a,b)$. However, by
\eqref{geomp2st}, the vertex-stars for $P_{1}(a,b)$ consist only of every other vector in the
corresponding vertex-star for $P(a,b)$. This indicates that the vertex-set and edge-set of
$P_{1}(a,b)$ are generally only a ``thin" subset of those of
$P(a,b)$. In fact, we have

\blem
\label{vertpone}
The vertex-set of $P_{1}(a,b)$ is given by
\[ V(P_{1}(a,b)) = \{(0,0,0), (a,0,b), (0,b,a), (b,a,0) \} + c\La_{(1,1,1)}, \]
with $c:=b-a$, where the four cosets of $c\La_{(1,1,1)}$ occurring are distinct. In particular, the vertex-set
is discrete (it is finite when $b=a$).
\elem

\bpf
Set $T_{0}:= I$. It follows from \eqref{specgone} and the proof of Lemma~\ref{ponetrans} that the
twelve elements
\[ S_{2}^{j}T_{k} \quad (j=0,1,2; k=0,1,2,3) \]
are a system of representatives of $G$ modulo $c\La_{(1,1,1)}$. Hence the vertex-set is 
the union of the twelve cosets $oS_{2}^{j}T_{k} + c\La_{(1,1,1)}$.  However, since $S_{2}$
fixes the base vertex $o$, there are only four cosets, namely $oT_{k} + c\La_{(1,1,1)}$ with $k =
0,1,2,3$.  Thus the vertex-set consists of the cosets of $o$ and its adjacent vertices.
Moreover, these four cosets are distinct. In fact, if $(a,0,b) \equiv (0,0,0)$ (say) modulo
$c\La_{(1,1,1)}$, then $(\tfrac{a}{c},0,\tfrac{b}{c}) = (\tfrac{a}{c},0,\tfrac{a}{c}+1) \in
\La_{(1,1,1)}$, so $\tfrac{a}{c}$ and $\tfrac{a}{c}+1$ must both be even integers; however, this is impossible. This proves the lemma.
\epf

Note that the vertex-set of $P_{1}(a,b)$ is discrete, irrespective of the values of $a$ and $b$. By 
contrast, if $a$ and $b$ are non-zero and are not rational multiples of each other, then the
vertex-set of $P(a,b)$ is dense in $\BE^3$; in fact, since $2\La$ is a subgroup of the translation group of $P(a,b)$, the vertex-set must contain the points in $2\La$ (see \cite[(5.13)]{schu}). 

Lemma~\ref{vertpone} also enables us to find the vertex-star at a given vertex $x$ of $P_{1}(a,b)$. 
In fact, since $c\La_{(1,1,1)}$ is the translation subgroup, the vertex-star at $x$ is the same as the vertex-star at the representative vertex of $x + c\La_{(1,1,1)}$ listed in the lemma. In other words, first determine which vertex among $(0,0,0)$, $(a,0,b)$, $(0,b,a)$, $(b,a,0)$ is equivalent to $x$ modulo $c\La_{(1,1,1)}$, and then take its vertex-star, $V_{j}$ (say). The vertices of $P_{1}(a,b)$ adjacent to $x$ then are the points in $x + V_{j}$. Recall that the vertex-stars at $(0,0,0)$, $(a,0,b)$, $(0,b,a)$ or $(b,a,0)$ are given by $V_{0}$, $V_{2}$, $V_{1}$ or $V_{3}$, respectively, so we have $j=0$, $2$, $1$ or $3$.

\blem
\label{mulverone}
$P_{1}(a,b)$ is a geometric polyhedron.
\elem

\bpf
We have altogether four translation classes of vertices and four vertex-stars, and the
vertex-stars at vertices of distinct translation classes are also distinct. Hence each
translation class of vertices is characterized by its unique vertex-star. It follows that
$P_{1}(a,b)$ cannot have multiple vertices; in other words, $P_{1}(a,b)$ must be a faithful
realization of the underlying abstract polyhedron. We can also see this more directly as follows.
Let $x$ be a vertex of $P_{1}(a,b)$ with vertex-star $V_{j}$, and let $y$ be a vertex adjacent to $x$
such that $y-x \in V_{j}$. We need to show that $x-y$ belongs to the vertex-star $V_{k}$ (say) at
$y$, so that $x \in y + V_{k}$. Since we can reduce $x$ modulo $c\La_{(1,1,1)}$, we may assume
that $x$ is among $(0,0,0)$, $(a,0,b)$, $(0,b,a)$ or $(b,a,0)$. However, now it is straightforward
to verify that the corresponding adjacent vertices and their vertex-stars have the desired property. 
\epf  

Notice that our methods also provide an alternative description of the edge-graph of $P_{1}(a,b)$
which does not appeal to the group. In fact, take the set of points described in
Lemma~\ref{vertpone} as vertex-set of the graph and place, at each vertex $x$, the corresponding
vertex-star from the list found by reducing $x$ modulo $c\La_{(1,1,1)}$. Then the resulting graph
coincides with the edge-graph of $P_{1}(a,b)$.

The vertex-figures of $P_{1}(a,b)$ are triangles. However, the vertex-stars of $P_{1}(a,b)$ are 
generally not planar. In fact, we have the following

\blem
\label{poneplanar}
The vertex-stars of $P_{1}(a,b)$ are planar if and only if $b=-a$. 
\elem

\bpf
The special group $G_{0}$ acts transitively on the vertex-stars. The three vectors in $V_0$ have
determinant $\pm (a^{3}+b^{3})$, so they lie in a plane if and only if $b=-a$.  
\epf

Next we determine which polyhedra are regular. 

\blem
\label{regpone}
The polyhedron $P_{1}(a,b)$ is geometrically chiral if $b \neq \pm a$, or geometrically regular if 
$b=\pm a$. In particular, $P_{1}(a,-a)$ is similar to $P_{1}(1,-1) = \{\infty,3\}^{(a)}$, and
$P_{1}(a,a)$ is similar to $P_{1}(1,1) = \{3,3\}$.
\elem

\bpf
The case $b = \pm a$ was already settled (see \eqref{infa} and \eqref{threeones}), so we may assume
that $b \neq \pm a$. Suppose $P_{1}(a,b)$ is regular. Then $P_{1}(a,b)$ must have an involutory
symmetry $R \;(=R_{2})$ which fixes $o$ and $v$ and interchanges the neighbors $(0,b,a)$ and
$(b,a,0)$ of $o$ (see also Lemma~\ref{regcrit}). Since $V_0$ is not planar, $R$ must necessarily be
the (linear) reflection in the plane through $o,v$ with normal vector 
\[ n := (0,b,a)-(b,a,0) = (-b,b-a,a). \]  
Moreover, since $R^{-1}TR = T$, the rotation axis of $T$ must be invariant under $R$. However,
the direction vector $e_{2} = (0,1,0)$ of this axis is not a scalar multiple of $n$, so $e_{2}$
must lie in the mirror of $R$. Hence $b-a = n \cdot e_{2} = 0$ (where $\cdot$ denotes the scalar
product), which contradicts our assumption on $a,b$. (Recall that $R_{2}$ is a half-turn for 
$\{\infty,3\}^{(a)}$, so a contradiction can only concern the case $b=a$.)  Now the lemma follows.
\epf

The two enantiomorphic forms of a chiral polyhedron are represented by different pairs of 
generators of its group. If $S_{1}(a,b),S_{2}$ is the pair associated with the base flag 
$\Ph(a,b) := \{F_{0}(a,b),F_{1}(a,b),F_{2}(a,b)\}$ of $P_{1}(a,b)$, then $S_{1}(a,b)S_{2}^{2},S_{2}^{-1}$ is the pair associated with the adjacent flag 
$\Ph^{2}(a,b) := \{F_{0}(a,b),F_{1}(a,b),F'_{2}(a,b)\}$ of $P_{1}(a,b)$. Since the product of
each pair is $T(a,b)$, we can also represent the two forms by the pairs of generators
$T(a,b),S_{2}$ and $T(a,b),S_{2}^{-1}$, respectively. When Wythoff's construction is applied
to $G_{1}(a,b)$ with the new pair of generators $T(a,b),S_{2}^{-1}$ and with the same initial
vertex $o$, which is also fixed by $S_{2}^{-1}$, then we again obtain the same polyhedron
$P_{1}(a,b)$, but now with a new base flag adjacent to the original base flag. We can see this
as follows.

First we claim that
\beq
\label{poneabba}
P_{1}(b,a) = P_{1}(a,b)R , 
\eeq
where $R$ is the reflection (in the plane $\xi_{1}=\xi_{3}$) given by
\beq
\label{err}
R\colon\;  (\xi_{1},\xi_{2},\xi_{3}) \; \mapsto \;  (\xi_{3},\xi_{2},\xi_{1}) .
\eeq
In fact, conjugation by $R$ transforms the generators $T(a,b),S_{2}$ for $G_{1}(a,b)$ to the 
generators $T(b,a),S_{2}^{-1}$ for $G_{1}(b,a)$, so in particular $G_{1}(b,a) = R^{-1}G_{1}(a,b)R$.
Then
\beq
\label{vpone}
V(P_{1}(b,a)) = o G_{1}(b,a) = o R^{-1}G_{1}(a,b)R = o G_{1}(a,b)R = V(P_{1}(a,b))R , 
\eeq
so $R$ maps the vertex-sets of the polyhedra onto each other.  Similarly, since
\beq
\label{epone}
F_{1}(b,a) G_{1}(b,a) = F_{1}(b,a) R^{-1}G_{1}(a,b)R = (F_{1}(a,b) G_{1}(a,b))R ,
\eeq
the same holds for the edge-sets. Finally, $R$ takes $F_{2}(b,a)$ to the face $F_{2}(a,b)' \;
(=F_{2}(a,b)T(a,b))$ of $\Ph^{2}(a,b)$ (see \eqref{facetspone}), so the corresponding
statement for face-sets follows from
\beq
\label{fapone}
F_{2}(b,a) G_{1}(b,a) = F_{2}(b,a) R^{-1}G_{1}(a,b)R = (F_{2}(a,b)' G_{1}(a,b))R . 
\eeq
This establishes \eqref{poneabba}.  Moreover, notice that 
\[ \Ph(b,a)R = \Ph^{2}(a,b), \quad \Ph^{2}(b,a)R =  \Ph(a,b) . \]

Now we can complete our argument about enantiomorphism. In fact, bearing in mind how conjugation by 
$R$ affects the generators we observe that the polyhedron $P$ (say) obtained by Wythoff's
construction from $G_{1}(a,b)$ with generators $T(a,b),S_{2}^{-1}$, is mapped by $R$ to the
polyhedron obtained from $G_{1}(b,a)$ with generators $T(b,a),S_{2}$. (Note that we
interchanged $a$ and $b$.)  However, the latter is just $P_{1}(b,a)$, so by \eqref{poneabba}
its preimage under $R$ is $P_{1}(a,b)$ itself. Hence, $P = P_{1}(a,b)$.

These considerations also justify our initial hypothesis that it suffices to concentrate on only 
one orientation for the generator $S_2$. The opposite orientation (given by $S_{2}^{-1}$) then
is implied by enantiomorphism.

Alternatively we could have appealed to similar such results for the polyhedra $P(a,b)$, by 
observing that conjugation by $R$ commutes with the facetting operation $\p_{2}$ (see
\cite[(5.18)]{schu}). In particular, 
\[ (P(a,b)R)^{\p_{2}} = (P(a,b)^{\p_{2}})R . \]

We now discuss the question when two polyhedra $P_{1}(a,b)$ and $P_{1}(c,d)$ with parameter sets 
$a,b$ and $c,d$, respectively, are affine images of each other. The result is consistent with
\cite[Lemma 5.15]{schu}, which describes the answer for the polyhedra $P(a,b)$ and $P(c,d)$.

\blem
\label{affeq}
Let $a,b,c,d$ be real numbers, and let $(a,b) \neq (0,0) \neq (c,d)$. Then the polyhedra
$P_{1}(a,b)$ and $P_{1}(c,d)$ are affinely equivalent if and only if $(c,d) = s(a,b)$ or $(c,d) =
t(b,a)$ for some non-zero scalars $s$ or $t$. Moreover, $P_{1}(a,b)$ and $P_{1}(c,d)$ are congruent
if and only of $(c,d) = \pm (a,b)$ or $(c,d) = \pm (b,a)$.
\elem

\bpf
We already know that the two polyhedra are affinely equivalent if $(c,d) = s(a,b)$ or $t(b,a)$ 
(see \eqref{simpone} and \eqref{poneabba}). 

Suppose $S$ is an affine transformation of $\BE^3$ such that $P_{1}(c,d)= P_{1}(a,b)S$. Since 
$G_{1}(c,d)$ acts transitively on the vertices, and the stabilizer of the vertex $o$ acts
transitively on the vertices adjacent of $o$, we can assume that $oS=o$ (so $S$ is linear) and
$(a,0,b)S=(c,0,d)$. First observe that $S$ and $S^{-1}$ preserve planarity of vertex-stars, so
certainly $b=-a$ if and only if $d=-c$; in this case the assertion is obvious. Let $b\neq -a$,
so now the vertex-stars are not planar. The affine transformation $S^{-1}T(a,b)S$ interchanges
the vertices in the base edge of $P_{1}(c,d)$, as well as the two faces meeting at this base
edge. Since this is also true for the symmetry $T(c,d)$ of $P_{1}(c,d)$, we must have
$S^{-1}T(a,b)S = T(c,d)$; recall that an affine transformation is uniquely determined by its
effect on four independent points. Similarly, $S^{-1}S_{2}S$ cyclically permutes the vertices
in the vertex-figure of $P_{1}(c,d)$ at $o$, so we must also have $S^{-1}S_{2}S = S_{2}$ or
$S_{2}^{-1}$.  These properties for the generators translate into conditions for the matrix entries
of $S$ and show that, up to similarity, there are only the following two possibilities. If
$S^{-1}S_{2}S = S_{2}$, then $S=sI$ for some scalar $s$, and hence $(c,d)=s(a,b)$. If $S^{-1}S_{2}S
= S_{2}^{-1}$, then $S = tR$ with $R$ as in \eqref{err}, for some scalar $t$; hence
$(c,d)=t(b,a)$. Moreover, if $S$ is an isometry, then necessarily $s=\pm 1$ or $t = \pm 1$. 
\epf

Figure~\ref{affclasspone} illustrates the situation described in Lemma~\ref{affeq}. Each affine
equivalence class of polyhedra is represented by a polyhedron $P_{1}(1,b)$ corresponding to the
point $(1,b)$ on the right side of the square. The equivalence classes of the regular polyhedra
$\{3,3\}$ and $\{\infty,3\}^{(a)}$ then occur at the upper right corner and lower right corner,
respectively. The other three sides of the square give alternative representations of the
equivalence classes, as indicated by the four symmetrically related nodes. 

\begin{figure}[hbt]
\centering
\begin{center}
\begin{picture}(180,190)
\put(80,90){
\begin{picture}(100,100)
\put(-90,0){\line(1,0){180}}
\put(90,0){\vector(1,0){2}}
\put(96,-3){$a$}
\put(0,-90){\line(0,1){180}}
\put(0,90){\vector(0,1){2}}
\put(-2,96){$b$}
\multiput(59.4,-60)(0.3,0){5}{\line(0,1){120}}
\put(-60,60){\line(1,0){120}}
\put(-60,-60){\line(1,0){120}}
\put(-60,-60){\line(0,1){120}}
\put(60,-60){\line(0,1){120}}
\put(85,85){\vector(-1,1){10}}
\put(85,85){\vector(1,-1){10}}
\put(80,95){$R$}
\put(60,60){\circle*{5}}
\put(-60,-60){\circle*{5}}
\put(-60,60){\circle*{5}}
\put(60,-60){\circle*{5}}
\put(47,67){$\{3,3\}$}
\put(40,-74){$\{\infty,3\}^{(a)}$}
\put(-80,67){$\{\infty,3\}^{(a)}$}
\put(-73,-74){$\{3,3\}$}
\put(20,60){\circle*{5}}
\put(5,48){\footnotesize ${P_{1}(b,1)}$}
\put(60,20){\circle*{5}}
\put(23,17){\footnotesize $P_{1}(1,b)$}
\put(61.5,-7){\scriptsize $1$}
\put(1,62){\scriptsize $1$}
\put(-60,-20){\circle*{5}}
\put(-20,-60){\circle*{5}}
\end{picture}}
\end{picture}
\caption{\it The affine classes of polyhedra $P_{1}(a,b)$.}
\label{affclasspone}
\end{center}
\end{figure}

The group $G_{1}(a,b)$ associated with $P_{1}(a,b)$ is generated by rotations and thus contains 
only direct isometries. This implies that the faces of $P_{1}(a,b)$ consist of either all
right-hand helices or all left-hand helices. Note that the faces of $P_{1}(a,b)$ are
right-hand helices if and only if the faces of $P_{1}(b,a)$ are left-hand helices; in fact,
the plane reflection $R$ of \eqref{err} maps $P_{1}(a,b)$ to $P_{1}(b,a)$, and vice versa.  When
$b=-a$, these polyhedra are regular and congruent, but they do not coincide (so $R$ is not a
symmetry); in fact, every symmetry of $\{\infty,3\}^{(a)}$ is a direct isometry. On the other
hand, if $b=a$, then $R$ is indeed a symmetry.

\begin{figure}[hbt]
\centering
\begin{center}
\begin{picture}(200,240)
\put(80,120){
\begin{picture}(100,100)
\put(0,0){\circle*{3}}
\put(0,-10){$o$}
\put(0,0){\oval(30,30)[t]}
\put(0,0){\oval(30,30)[r]}
\put(10.5,10.5){\vector(1,-1){2}}
\put(13,13){$S_{2}$}
\put(0,0){\line(1,0){80}}
\put(80,0){\line(3,5){25}}
\put(80,0){\line(3,-5){25}}
\put(80,0){\circle*{3}}
\put(85,-3){$v$}
\put(35,45){$F_{2}$}
\put(30,-50){$F_{2}S_{2}$}
\put(-65,-3){$F_{2}S_{2}^{2}$}
\put(110,-3){$F_{2}S_{2}^{2}T$}
\put(0,0){\line(-3,5){40}}
\put(-40,-66){\circle*{3}}
\put(-40,-66){\line(-1,0){50}}
\put(-40,-66){\line(3,-5){25}}
\put(-35,-69){$vS_{2}$}
\put(-85,-100){$F_{2}S_{2}^{2}TS_{2}$}
\put(0,0){\line(-3,-5){40}}
\put(-40,66){\circle*{3}}
\put(-40,66){\line(-1,0){50}}
\put(-40,66){\line(3,5){25}}
\put(-35,63){$vS_{2}^{2}$}
\put(-85,96){$F_{2}S_{2}^{2}TS_{2}^{2}$}
\end{picture}}
\end{picture}
\caption{\it The faces near $o$.}
\label{fanear}
\end{center}
\end{figure}
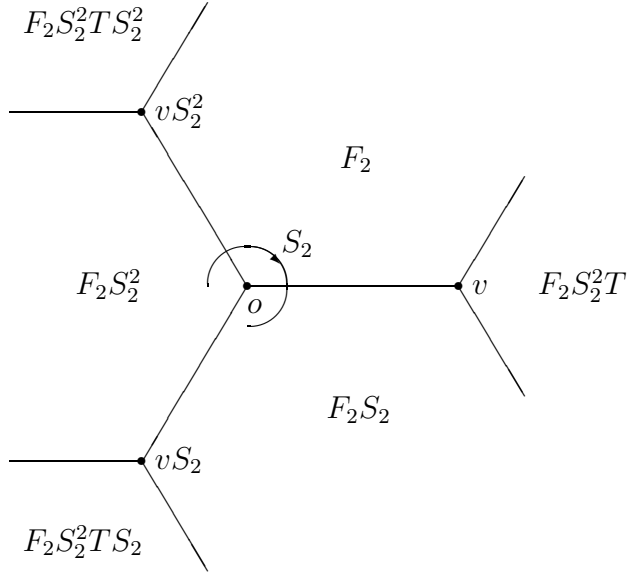

We conclude this section with a discussion of the {\em translation classes\/} of faces of
$P_{1}(a,b)$, by which we mean the transitivity classes of faces under the translation subgroup
$c\La_{(1,1,1)}$ of $G=G_{1}(a,b)$, where again $c:=b-a$. When $P_{1}(a,b)$ is regular, no
translational symmetries in addition to those in $G_{1}(a,b)$ occur, so the translation classes
with respect to the full symmetry group of $P_{1}(a,b)$ are the same. 

\blem 
\label{trclpone}
$P_{1}(a,b)$ has precisely four translation classes of faces. They are represented by the four 
faces $F_{2}$, $F_{2}S_{2}$, $F_{2}S_{2}^{2}$ and $F_{2}S_{2}^{2}T$, where $F_{2}$ is the base
face given by \eqref{facetspone}.
\elem

\bpf
By Lemma~\ref{vertpone} each vertex of $P_{1}(a,b)$ is equivalent to $o$, $v$, $vS_{2}$ or 
$vS_{2}^{2}$ modulo $c\La_{(1,1,1)}$. Hence, modulo $c\La_{(1,1,1)}$, every face of $P_{1}(a,b)$ is
equivalent to a face containing one of these vertices. There are six such faces, including the
three which contain $o$, namely $F_{2}$, $F_{2}S_{2}\;(= F_{2}')$ and $F_{2}S_{2}^{2}$. The other
three faces, $F_{2}S_{2}^{2}T$, $F_{2}S_{2}^{2}TS_{2}$ and $F_{2}S_{2}^{2}TS_{2}^{2}$, contain
the vertices $v$, $vS_{2}$ and $vS_{2}^{2}$, respectively. The combinatorial picture is illustrated
in Figure~\ref{fanear}. In the first set of three, the helical faces have an axis with direction
vector $(1,-1,-1)$, $(-1,-1,1)$ or $(-1,1,-1)$, respectively. However, in the second set, the faces
are all equivalent modulo $c\La_{(1,1,1)}$ and have an axis with direction vector $(1,1,1)$.  In
particular,
\[ F_{2}S_{2}^{2}TS_{2}    = F_{2}S_{2}^{2}T + (c,c,-c), \quad 
F_{2}S_{2}^{2}TS_{2}^{2} = F_{2}S_{2}^{2}T + (2c,0,0). \]
This concludes the proof.
\epf

\section{Type $\{\infty,3\}$, with helical faces over squares}
\label{helfactwo}

In this section we describe another family of chiral apeirohedra of type $\{\infty,3\}$, now
with  helical faces over squares. They occur when $G_{0}=[3,4]^{+}$ in Lemma~\ref{chirspecgr}.
As in the previous section, each polyhedron is obtained from an irreducible group $G$ by Wythoff's
construction with initial vertex $F_{0}:=o$. 

For a polyhedron $P$ of type $\{\infty,3\}$ with helical faces over squares, we must begin with a 
group $G=\scl{S_{1},S_{2}}$ with special group $G_{0} = [3,4]^{+}$, where $S_{1}',S_{2}$ are
standard generators of $G_0$ considered as the group of the cube $\{4,3\}$, so that $S_{1}'$
and $S_{2}$ have periods $4$ and $3$, respectively, and their product has period $2$. (Recall
that $R'$ denotes the image in $G_0$ of an element $R$ of $G$.)  As already mentioned in the
proof of Lemma~\ref{twi}, once $S_{2}$ has been selected, there are three equivalent choices
for $S_{1}'$ such that $S_{1}'S_{2}$ has period $2$; if $S_{1}'$ is one of them, then the other
two are $S_{2}^{-1}S_{1}'S_{2}$ and $S_{2}^{-2}S_{1}'S_{2}^{2}$. Moreover, substituting
$S_{2}$ by its inverse would lead to (the enantiomorphic image of) the same underlying
polyhedron, so we need not discuss this case separately.

As in the previous section, we may confine ourselves to some very specific choices for the 
generators. Thus, with
\beq 
\label{helsqasg}
\bry{rccl}
S_{1}'\colon & x & \mapsto &  (-\xi_{3},\xi_{2},\xi_{1}) ,\\
S_{2}\colon  & x & \mapsto &  (\xi_{2},\xi_{3},\xi_{1}),  \\
\ery  
\eeq
we observe that $T'=S_{1}'S_{2}$ is the half-turn about the line $\xi_{2}=\xi_{1}$ in the 
$\xi_{1}\xi_{2}$-plane. Then $T=S_{1}S_{2}$ is a half-turn about a (parallel) line
perpendicular to the plane $\xi_{2}=-\xi_{1}$, and since $T$ must interchange the two vertices
in the base edge $F_{1}$, this base edge must lie in this plane and the rotation axis of $T$
must pass through its midpoint.  In particular, $v:=oT =  (c,-c,d)$ for some real parameters
$c$ and $d$, not both zero. The resulting group $G = G_{2}(c,d)$ is generated by
\beq 
\label{helsqagr}
\bry{rccl}
S_{2}\colon  & x & \mapsto &  (\xi_{2},\xi_{3},\xi_{1}),  \\
T    \colon  & x & \mapsto &  (\xi_{2},\xi_{1},-\xi_{3}) + (c,-c,d).
\ery  
\eeq
For $S_{1} = TS_{2}^{-1}$ we have
\beq 
\label{helsqatwo}
\bry{rccl}
S_{1}    \colon & x & \mapsto &  (-\xi_{3},\xi_{2},\xi_{1}) + (d,c,-c), \\
S_{1}^{4}\colon & x & \mapsto &  (\xi_{1},\xi_{2},\xi_{3}) + 4c(0,1,0). \\ 
\ery  
\eeq
The latter is the translation by $4c(0,1,0)$, and is trivial if $c=0$. The special case $c=0$
yields the finite group 
\[ G_{2}(0,d) \cong [3,4]^+  \]
(with fixed point $\half (d,d,d)$); the corresponding polyhedron is a cube.

The polyhedron $P = P_{2}(c,d)$ is obtained from $G=G_{2}(c,d)$ by Wythoff's construction with base
vertex $F_{0}=o$. Its base edge $F_{1}$ has vertices $o$ and $v = (c,-c,d)$; its base face is
given by
\beq
\label{facetsptwo}
F_{2} = o\scl{S_{1}} = \{(c,-c,d),(0,0,0),(d,c,-c),(c+d,2c,-c+d)\} + {\Bbb Z}\!\cdot\!4c(0,1,0) 
\eeq
and consists of a helix over a square. (The notation for helical faces is as before, but now there
are four consecutive vertices that are  translated.)  As usual, the remaining vertices, edges and
faces of $P$ are the images of $F_{0}$, $F_{1}$ and $F_{2}$ under $G$. The vertices adjacent to $o$
are 
\[ v = (c,-c,d),\; vS_{2} = (-c,d,c),\; vS_{2}^{2} = (d,c,-c) . \]
When the polyhedra are only considered up to similarity, the two parameters $c,d$ reduce to a single
parameter. In fact, \eqref{simpone} carries over and shows that $P_{2}(c,d)$ is similar to
$P_{2}(1,\tfrac{d}{c})$ or $P_{2}(0,1)$ according as $c\neq 0$ or $c=0$. 

As in the previous section, the facetting operation $\p_{2}$ establishes a link to polyhedra with 
finite faces. Now the polyhedra $P_{2}(c,d)$ are related to the chiral (or regular) polyhedra
$Q(c,d)$ of type $\{4,6\}$ described in \cite[\S 6]{schu}; the latter generally have skew square
faces and skew hexagonal vertex-figures. Recall that the generators
$\widehat{S}_{1},\widehat{S}_{2}$ (say) and corresponding element $\widehat{T} :=
\widehat{S}_{1}\widehat{S}_{2}$ of the symmetry group $H(c,d)$ of $Q(c,d)$ are given by
\beq 
\label{quadhexgroup}
\bry{rccl}
\widehat{S}_{1}\colon & x & \mapsto & (-\xi_{1},\xi_{3},-\xi_{2}) + (c,-d,-c),\\
\widehat{S}_{2}\colon & x & \mapsto & -(\xi_{3},\xi_{1},\xi_{2}),  \\
\widehat{T}\colon & x & \mapsto &  (\xi_{2},\xi_{1},-\xi_{3}) + (c,-c,d).  
\ery  
\eeq
Now we have

\blem
\label{p2pcd} 
$P_{2}(c,d) = Q(c,d)^{\p_2}$, for all real parameters $c$ and $d$. 
\elem

\bpf 
In terms of the generators $\widehat{S}_{2},\widehat{T}$ of $H(c,d)$, the operation $\p_{2}$ takes 
the same form as in \eqref{geomp2st}. Now $P_{2}(c,d) = Q(c,d)^{\p_2}$ and $G_{2}(c,d)$ is a
subgroup of $H(c,d)$.
\epf

This relationship holds again without any restriction on the parameters, so in particular we may
(and often do) have the situation that $P_{2}(c,d)$ is discrete but $Q(c,d)$ is not. Moreover,
the lemma implies that the faces of $P_{2}(c,d)$ are among the holes of $Q(c,d)$, which now are
helices over squares, except when $c=0$ (then they are squares).

The facetting operation $\p_{2}$ preserves regularity, so again there are two special cases.
From $Q(1,0)= \{4,6\}_{6}$ we obtain the regular polyhedron
\beq
\label{inf3b}
P_{2}(1,0) =  \{\infty,3\}^{(b)} 
\eeq
(notation as in \cite[p.224]{msarp}); all other polyhedra $P_{2}(c,0)$ are similar to
$P_{2}(1,0)$.  Similarly, the case $c=0$ reduces to $(c,d)=(0,1)$ and again yields regular
polyhedra. Now $Q(0,1) = \{4,6 \hole 4\}$ and $\{4,6 \hole 4\}^{\p_2} = \{4,3\}$ (see
\cite[Thm.6.12]{schu} and \cite[p.224]{msarp}), so in particular
\beq
\label{p201}
P_{2}(0,1) =  \{4,3\} .
\eeq

Next we determine the translation subgroup $T(G)$ of $G=G_{2}(c,d)$. Since $S_{1}^4$ is the 
translation by $4c(0,1,0)$, and its conjugates under the special group $G_{0}\;( \cong
[3,4]^{+})$ consist of the vectors $\pm 4c(1,0,0)$, $\pm 4c(0,1,0)$ and $\pm 4c(0,0,1)$, we
immediately see that
\[ 4c \BZ^{3} = 4c\La_{(1,0,0)} \leq T(G) .\]
In fact, we have 

\blem
\label{ptwotrans}
The translation subgroup of $G_{2}(c,d)$ is given by $T(G_{2}(c,d)) = 4c \BZ^{3}$.
\elem

\bpf
The proof is more involved than the proof of Lemma~\ref{ponetrans} and is interesting in its
own right. As before, since $G = \scl{S_{2},T}$, we have $G = N \cdot \scl{S_{2}}$ (as a product of
subgroups), where $N := \scl{T_{1},T_{2},T_{3}}$ with $T_{2}:=T$, $T_{1} := S_{2}^{-1}TS_{2}$ and
$T_{3} :=  S_{2}^{-2}TS_{2}^{2}$. Now $S_{2}$ is the square of the generator $\widehat{S}_{2}$ for
the group $H(c,d)$ of $Q(c,d)$ (see \eqref{geomp2st} and \eqref{quadhexgroup}), so
$T_{1},T_{2},T_{3}$ are just the half-turns occurring in \cite[\S 6]{schu} with the same notation,
namely
\beq 
\label{t123ptwo}
\bry{rccl}
T_{1}\colon & x & \mapsto &  (\xi_{3},-\xi_{2},\xi_{1}) + (-c,d,c),  \\
T_{2}\colon & x & \mapsto &  (\xi_{2},\xi_{1},-\xi_{3}) + (c,-c,d), \\
T_{3}\colon & x & \mapsto &  (-\xi_{1},\xi_{3},\xi_{2}) + (d,c,-c). 
\ery  
\eeq
It is straightforward to check that $(T_{i}T_{j})^{3}$ and $(T_{i}T_{j}T_{k})^{4}$, but not 
$T_{i}T_{j}$ and $T_{i}T_{j}T_{k}$, are translations in $4c \BZ^{3}$, for all mutually distinct
$i,j,k$.

We now consider the six elements $T_{i}T_{j}T_{k}T_{j}$ of $N$ with $i,j,k$ distinct. They are
\beq
\label{theui}
\begin{array}{l}
U_{1}:=T_{1}T_{2}T_{3}T_{2},\quad U_{2}:=T_{2}T_{1}T_{3}T_{1},\quad U_{3}:=T_{3}T_{1}T_{2}T_{1},\\
\widehat{U}_{1}:=T_{1}T_{3}T_{2}T_{3},\quad  \widehat{U}_{2} :=T_{2}T_{3}T_{1}T_{3},\quad\,
\widehat{U}_{3}:=T_{3}T_{2}T_{1}T_{2}.
\end{array}
\eeq
In particular, we have
\beq 
\label{u123}
\bry{rccl}
U_{1}\colon & x & \mapsto &  (-\xi_{1},\xi_{2},-\xi_{3}) + (c+d,-2c,3c+d),  \\
U_{2}\colon & x & \mapsto &  (-\xi_{1},-\xi_{2},\xi_{3}) + (-c+d,-3c+d,2c), \\
U_{3}\colon & x & \mapsto &  (\xi_{1},-\xi_{2},-\xi_{3}) + (-2c,3c+d,c+d). 
\ery  
\eeq
Then $U_{i}^{2}$, $\widehat{U}_{i}^{2}$ and $U_{i}\widehat{U}_{i}$ are translations in
$4c\BZ^{3}$, for each $i$. Now define the subgroups
\beq
\label{kk2}
K := \scl{U_{1},U_{2},U_{3}},\quad K^{2} := \scl{U_{1}^{2},U_{2}^{2},U_{3}^{2}}
\eeq
of $N$. Then it is easy to see that $4c \BZ^{3} = K^{2} \leq K$, and hence also
$\widehat{U}_{1},\widehat{U}_{2},\widehat{U}_{3} \in K$. In particular, modulo $K^{2}$, each 
element $U_{i}$ has period $2$ and coincides with $\widehat{U}_{i}$. Moreover, each product
$U_{i}U_{j}U_{k}$, with $i,j,k$ distinct, is a translation in $K^2$. 

Next observe that $S_{2}\in N$, so in particular $N = G$. In fact, $S_{1}T_{1}U_{3}$  is a
translation in $K^{2}$, so 
\[ S_{1}\in K^{2}U_{3}T_{1}  \subseteq K T_{1} \subseteq N .\] 
Hence also $S_{2} = S_{1}^{-1}T_{2}\in N$. 

Now, modulo $K$, each element of $G$ can be represented by an element (a right coset
representative) which is a product  of at most three generators $T_{i}$. In fact, if $\varphi =
T_{i_{1}}T_{i_{2}}\ldots T_{i_{n}} \in G$, then we can split off elements of the form
$T_{i}T_{j}T_{k}T_{j}$ and proceed inductively; for example, if $i_{1},i_{2},i_{3}$ are distinct,
then insert the trivial product $T_{i_{2}}T_{i_{2}}$ between the third and fourth term of $\varphi$,
split off $T_{i_{1}}T_{i_{2}}T_{i_{3}}T_{i_{2}}$, and proceed with $T_{i_{2}}T_{i_{4}}\ldots
T_{i_{n}}$. 

Moreover, modulo $K^{2}$, each element of $K$ can be represented by $I$, $U_{1}$, $U_{2}$ or
$U_{3}$ (as right coset representatives), so in particular the subgroup $K^{2}$ has index $4$ in
$K$. In fact, suppose $\psi = U_{i_{1}}U_{i_{2}}\ldots U_{i_{m}} \in K$; modulo $K^{2}$, each
generator has period $2$, so we need not consider any product involving their inverses. In this
case we can split off elements of the form $U_{i}U_{j}U_{k}$ and again proceed inductively; for
example, if $i_{1},i_{2}$ are distinct and $j\neq i_{1},i_{2}$, then insert the (trivial, modulo
$K^2$) product $U_{j}U_{j}$ between the second and third term of $\psi$, split off
$U_{i_{1}}U_{i_{2}}U_{j}$, and proceed with $U_{j}T_{i_{3}}\ldots U_{i_{m}}$. This, then,
reduces $\psi$ modulo $K^{2}$ to an element which is a product of at most two generators
$U_{i}$. If in fact this product involves two generators and is given by $U_{i_{1}}U_{i_{2}}$
(say), then we can further reduce it to $U_{j}$ with $j\neq i_{1},i_{2}$, by multiplying by
$U_{j}U_{j}$ on the right and again splitting off $U_{i_{1}}U_{i_{2}}U_{j}$. 

Now it follows that each translation in $K$ already lies in $K^{2}$; that is, $T(G) \cap K =
K^{2}$. In fact, none of the elements $U_{1}$, $U_{2}$, $U_{3}$ is a translation, so any
translation in $K$ must be equivalent to $I$ modulo $K^{2}$. Note that this also shows that
$K^{2}$ is normal in $K$. In particular, $K\slash K^{2} \cong C_{2}\times C_{2}$, with $C_{2}\times
C_{2}$ generated by any two of the generators $U_{i}$.

We now complete the proof by considering images of elements in the special group $G_{0}$. 
First observe that, since this is true for the generators, each element of $K$, considered in
the special group, takes the form $x \mapsto (\pm\xi_{1},\pm\xi_{2},\pm\xi_{3})$, with exactly
two minus-signs or no minus-sign (for the trivial element); in particular, no permutation of
coordinates occurs. On the other hand, modulo $K$, each element of $G$ can be represented by an
element which is a product of at most three generators $T_{i}$, and it is not hard to see that in
fact the six elements 
\[ I,T_{1},T_{2},T_{3},T_{1}T_{2},T_{2}T_{1}\]
suffice and give a complete system of right coset representatives of $G$ modulo $K$. However, the
elements in each  coset of $K$ distinct from $K$ all involve a (in fact, the same) non-trivial
permutation of the coordinates, because this is true for the corresponding representatives. Hence
no such coset can contain a translation, because any translation is trivial in the special group.
Now it follows that $T(G) = 4c \BZ^{3}$.

An alternative system of representatives of $G$ modulo $K$ consists of the six elements
\beq
\label{altsys}
I,T_{1},S_{2},T_{1}S_{2},S_{2}^{2},T_{1}S_{2}^{2} , 
\eeq 
which turns out to be more useful for the vertex-set computation. From this we also obtain a system 
of right coset representatives for $G$ modulo $4c \BZ^{3}$, by multiplying each of these elements
on the left by $I$, $U_{1}$, $U_{2}$ or $U_{3}$.
\epf

The translation subgroup depends only on a single parameter, $c$. Now Lemma~\ref{ponetrans} can be 
rephrased as
\beq
\label{specgtwo}
G_{2}(c,d)\slash 4c\BZ^{3}  \cong [3,4]^{+} \cong S_{4} ,
\eeq
with each term representing the special group. 

We now compute the vertex-set of $P_{2}(c,d)$. Again the vertex-set is discrete, irrespective of 
the values of the parameters $c$ and $d$. 

\blem
\label{vertptwo}
The vertex-set of $P_{2}(c,d)$ is given by
$V(P_{2}(c,d)) = {\displaystyle \bigcup_{z \in Z} (z + 4c\BZ^{3})}$, 
where $Z$ consists of the points
\[\left\{ \begin{array}{l}
(0,0,0);\\
(c+d,2c,-c+d), \mbox{ and its cyclic shifts}; \\
(2c+d,2c+d,2c+d);\\
(c,-c,d), \mbox{ and its cyclic shifts}.
\end{array} \right.\]
\elem

\bpf
Let $U_{1},U_{2},U_{3}$ and $K$ be as in the proof of the previous lemma (see \eqref{u123} and 
\eqref{kk2}). First observe that $oK$, the orbit of the base vertex $o$ under $K$, is the union
of the translation classes of the four vertices $o$, $oU_{1}$, $oU_{2}$, $oU_{3}$ modulo
$4c\BZ^{3}$, or, equivalently, of the four vertices
\beq
\label{fvtwo1}
(0,0,0),  (c+d,2c,-c+d), (-c+d,c+d,2c), (2c,-c+d,c+d) ,
\eeq
respectively. Then $oKT_{1}$, the image of $oK$ under $T_{1}$, is the union of the translations
classes of the four vertices
\beq
\label{fvtwo2}
(-c,d,c), (2c+d,2c+d,2c+d), (c,-c,d), (d,c,-c) . 
\eeq
These two sets of four vertices are invariant under $S_{2}$, so their union is also invariant under 
$S_{2}$. It follows that no additional vertices arise from the images of $oK$ under the elements
$S_{2}$, $T_{1}S_{2}$, $S_{2}^{2}$ or $T_{1}S_{2}^{2}$, that is, from the remaining cosets
$KS_{2}$, $KT_{1}S_{2}$, $KS_{2}^{2}$ and $KT_{1}S_{2}^{2}$ of $K$ (see \eqref{altsys}). This proves
the lemma.
\epf

Observe that, when $c=0$, the polyhedron is the cube $P_{2}(0,d)$ with vertex-set $Z$.

For certain parameter values the eight translation classes of vertices may not all be distinct. 
However, bearing in mind that the set of eight translation classes is invariant under $S_{2}$,
inspection shows that two classes can only coincide if $d=kc$ with $k\in\BZ$ and $k \equiv 2 \mod
4$. This occurs precisely when $(2c+d,2c+d,2c+d)$ is equivalent to $(0,0,0)$ modulo $4c\BZ^{3}$,
and then $(c,-c,d)$ and $(-c+d,c+d,2c)$, as well as their cyclic shifts, are also equivalent modulo
$4c\BZ^{3}$.  However, no further collapses occur, because $(c,-c,d)$ and its cyclic shifts are not
equivalent to $(0,0,0)$ modulo $4c\BZ^{3}$ for any parameter values. This proves 

\blem
\label{verttransptwo}
$P_{2}(c,d)$ has eight translation classes of vertices modulo $4c\BZ^{3}$, except when $d=kc$ and 
$k$ is an integer with $k \equiv 2 \mod 4$; they are represented by the vertices listed in
Lemma~\ref{vertptwo}. If $d=kc$ and $k$ is an integer with $k \equiv 2 \mod 4$, then
$P_{2}(c,d)$ has four translation classes of vertices, represented by 
\[\left\{ \begin{array}{l}
(0,0,0);\\
(c,-c,d), \mbox{ and its cyclic shifts}.
\end{array} \right.\]
\elem

Recall from \cite[\S 6]{schu} that the polyhedron $Q(c,d)$, from which $P_{2}(c,d)$ is obtained by 
the facetting operation, has vertices of multiplicity $2$ (every point in $\BE^3$ taken by a
vertex is occupied by exactly $2$ vertices), and thus is not a faithful realization of the underlying abstract polyhedron, if $c$ and $d$ are integers with $c$ odd and $d \equiv 2 \mod 4$; for other integral parameters it has single vertices. This explains why the number of translation classes of vertices of $P_{2}(c,d)$ collapses to $4$ in certain cases. In fact, up to similarity, the polyhedra with four translation classes are all of the form $P_{2}(1,k)$ with $k\in\BZ$ and $k \equiv 2 \mod 4$, so the corresponding polyhedron $Q(1,k)$ does have vertices of multiplicity $2$. However, there are also cases with
$c$ odd and $d \equiv 2 \mod 4$ where $P_{2}(c,d)$ still has eight translation classes of vertices.
This happens because $P_{2}(c,d)$ has fewer translational symmetries than $Q(c,d)$; in fact,
the translation group of $Q(c,d)$ is $2\La_{(1,1,0)}$.

We now determine the vertex-stars of $P_{2}(c,d)$. First observe that, since $G=G_{2}(c,d)$ acts 
transitively on the vertices of $P_{2}(c,d)$, the vertex-stars of $P_{2}(c,d)$ are images of
the vertex-star $W_0$ at $o$ under the special group $G_{0}\;(\cong S_{4})$, so $G_{0}$ certainly
acts transitively on them. The stabilizer of $W_0$ contains $\scl{S_{2}}$, and hence the number of
vertex-stars cannot exceed $8$. This is consistent with the numbers $8$ or $4$ of transitivity
classes of vertices; clearly, vertices in the same translation class must have the same vertex-star.

The vertex-stars at the eight vertices occurring in Lemma~\ref{vertptwo} are given by
\beq 
\label{vtwostars}
\bry{lllll}
W_{0} &:=&             & &\{ (c,-c,d),(-c,d,c),(d,c,-c) \},\\ 
W_{1} &:=& W_{0}U_{1}' &=& \{ (-c,-c,-d),(c,d,-c),(-d,c,c) \},\\ 
W_{2} &:=& W_{0}U_{2}' &=& \{ (-c,c,d),(c,-d,c),(-d,-c,-c) \},\\ 
W_{3} &:=& W_{0}U_{3}' &=& \{ (c,c,-d),(-c,-d,-c),(d,-c,c) \},\\ 
W_{4} &:=& W_{0}T_{1}' &=& \{ (-c,-c,d),(c,-d,-c),(d,c,c) \},\\ 
W_{5} &:=& W_{0}U_{1}'T_{1}' &=& \{ (c,-c,-d),(-c,-d,c),(-d,c,-c) \},\\ 
W_{6} &:=& W_{0}U_{2}'T_{1}' &=& \{ (-c,c,-d),(c,d,c),(d,-c,-c) \},\\  
W_{7} &:=& W_{0}U_{3}'T_{1}' &=& \{ (c,c,d),(-c,d,-c),(-d,-c,c) \}, 
\ery  
\eeq
where the labelling accords with the order in which the vertices are listed in \eqref{fvtwo1} and
\eqref{fvtwo2}; that is, $W_{0},\ldots,W_{3}$ for the vertices in \eqref{fvtwo1}, and 
$W_{4},\ldots,W_{7}$ for the vertices in \eqref{fvtwo2}, in each case with the orders
matching. 

Each vertex-star consists of three of six vectors in a vertex-star of the corresponding polyhedron 
$Q(c,d)$. (In the notation of \cite[\S 6]{schu}, these vertex-stars are $W_{0}$, $W_{7}$,
$W_{6}$, $W_{5}$, $W_{1}$, $W_{4}$, $W_{2}$ and $W_{3}$, respectively, so in particular the
present labelling does not accord with the old one. The reflection $R$ of \cite[(6.8)]{schu} is not
an element of $G$, so it cannot be employed here.)  The vertex-stars at the vertices
$(c,-c,d)$, $(-c,d,c)$, $(d,c,-c)$ adjacent to $o$ are $W_{6}$, $W_{4}$, $W_{7}$, respectively. 

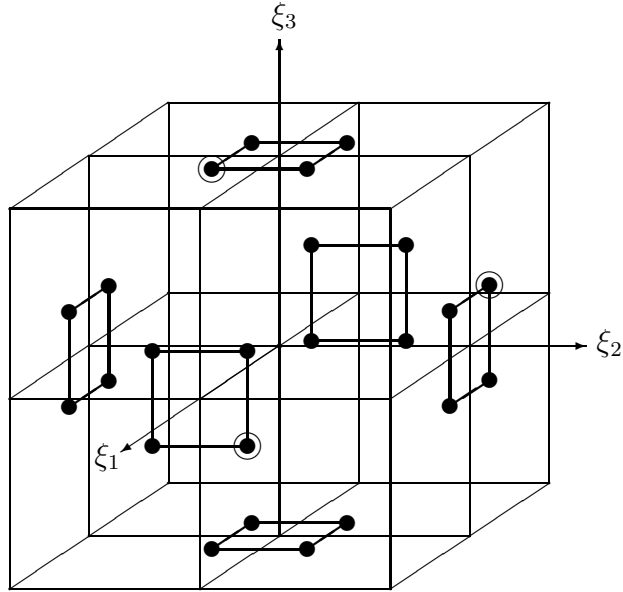
\begin{figure}[hbt]
\centering
\begin{center}
\begin{picture}(220,250)
\multiput(10,20)(0,72){3}{ 
\begin{picture}(150,40)
\thinlines
\multiput(0,0)(72,0){3}{\circle*{1}}
\multiput(30,20)(72,0){3}{\circle*{1}}
\multiput(60,40)(72,0){3}{\circle*{1}}
\multiput(0,0)(30,20){3}{\line(1,0){144}}
\multiput(0,0)(72,0){3}{\line(3,2){60}}
\end{picture}}
\put(10,20){ 
\begin{picture}(150,200)
\thinlines
\multiput(0,0)(72,0){3}{\line(0,1){144}}
\multiput(30,20)(72,0){3}{\line(0,1){144}}
\multiput(60,40)(72,0){3}{\line(0,1){144}}
\end{picture}}
\put(10,20){ 
\begin{picture}(150,200)
\multiput(54,54)(36,0){2}{\circle*{6}}
\put(90,54){\circle{10}}
\multiput(54,90)(36,0){2}{\circle*{6}}
\multiput(22.5,68.83)(0,36){2}{\circle*{6}}
\multiput(37.5,78.83)(0,36){2}{\circle*{6}}
\multiput(76.5,15)(36,0){2}{\circle*{6}}
\multiput(91.5,25)(36,0){2}{\circle*{6}}
\end{picture}}
\put(214,204){ 
\begin{picture}(150,200)
\multiput(-54,-54)(-36,0){2}{\circle*{6}}
\multiput(-54,-90)(-36,0){2}{\circle*{6}}
\multiput(-22.5,-68.83)(0,-36){2}{\circle*{6}}
\put(-22.5,-68.83){\circle{10}}
\multiput(-37.5,-78.83)(0,-36){2}{\circle*{6}}
\multiput(-76.5,-15)(-36,0){2}{\circle*{6}}
\multiput(-91.5,-25)(-36,0){2}{\circle*{6}}
\put(-127.5,-25){\circle{10}}
\end{picture}}
\put(10,20){ 
\begin{picture}(150,200)
\thicklines
\multiput(54,54)(0,36){2}{\line(1,0){36}}
\multiput(54,54)(36,0){2}{\line(0,1){36}}
\multiput(22.5,68.83)(0,36){2}{\line(3,2){15}}
\multiput(22.5,68.83)(15,10){2}{\line(0,1){36}}
\multiput(76.5,15)(15,10){2}{\line(1,0){36}}
\multiput(76.5,15)(36,0){2}{\line(3,2){15}}
\end{picture}}
\put(214,204){ 
\begin{picture}(150,200)
\thicklines
\multiput(-54,-54)(0,-36){2}{\line(-1,0){36}}
\multiput(-54,-54)(-36,0){2}{\line(0,-1){36}}
\multiput(-22.5,-68.83)(0,-36){2}{\line(-3,-2){15}}
\multiput(-22.5,-68.83)(-15,-10){2}{\line(0,-1){36}}
\multiput(-76.5,-15)(-15,-10){2}{\line(-1,0){36}}
\multiput(-76.5,-15)(-36,0){2}{\line(-3,-2){15}}
\end{picture}}
\put(10,20){ 
\begin{picture}(150,200)
\thinlines
\put(102,92){\vector(1,0){116}}
\put(222,90){$\xi_{2}$}
\put(102,92){\vector(0,1){116}}
\put(99,214){$\xi_{3}$}
\put(102,92){\vector(-3,-2){60}}
\put(32,47){$\xi_{1}$}
\put(102,92){\circle*{2}}
\end{picture}}
\end{picture}
\caption{\it The set $W$ for $c=1$ and $d=4$. The points in the vertex-star $W_{0}$ for
$P_{2}(c,d)$ are circled. The cubes are drawn in for reference.}
\label{ptwostarw}
\end{center}
\end{figure}

Just as in similar situations before, the union of all vertex-stars of $P_{2}(c,d)$ comprises the 
set of (generally twenty-four) vectors
\beq
\label{doubu}
W  :=  \bigcup_{i=0}^{7} W_{i}  
= \{ (\pm c, \pm c,\pm d),  (\pm c, \pm d, \pm c), (\pm d, \pm c, \pm c) \},      
\eeq
which is the same as the corresponding set for $Q(c,d)$ (see \cite[(6.13)]{schu}). This is also the 
orbit of $v$ under $G_{0}$. Note that $W$ is the vertex-set of one of the following convex
polyhedra:\ a suitably truncated octahedron or cube, with $24$ vertices, if $0<|c|<|d|$ or
$0<|d|<|c|$, respectively; a cube if $c = \pm d$; an octahedron if $c=0$; or a cuboctahedron if
$d=0$. Figure~\ref{ptwostarw} shows the vertex-set of the truncated octahedron obtained for $c=1$
and $d=4$; the fat lines indicate the square faces. 

The edge-module of $P_{2}(c,d)$ is given by $\La := \BZ[W]$, the $\BZ$-module spanned by 
the vectors in $W$. As before, $V(P_{2}(c,d)) \subset \La$. Recall from \cite[Lemma~6.2]{schu}
that $\La$ is a lattice if and only if $c,d$ can be scaled to a pair of relatively prime integers.
Hence we observe that, while the vertex-set of $P_{2}(c,d)$ is always discrete, its
edge-module is generally not discrete. On the other hand, the vertex-set of $Q(c,d)$ is discrete
if and only if its edge-module is discrete; in fact, in this case, $2\La$ is a subset of the
vertex-set.

The eight vertex-stars in \eqref{vtwostars} are mutually disjoint (and hence distinct), unless 
$c =\pm d$ or $c=0$ or $d=0$; the circled points in Figure~\ref{ptwostarw} indicate the
vertex-star $W_{0}$ for $c=1$ and $d=4$. The vertex-stars are still mutually distinct (but not
disjoint) if $c = \pm d$ or $c=0$; in these cases they correspond to the triangular faces of a
stella octangula (a pair of regular tetrahedra inscribed in a cube as in Figure~\ref{tetcub}) or
an octahedron, respectively. Finally, if $d=0$, the vertex-stars coincide in pairs, namely we have
\beq
\label{coinc}
W_{0}=W_{5},\; W_{1}=W_{4},\; W_{2}=W_{6},\; W_{3}=W_{7};
\eeq
now the vertex-stars consist of alternating vertices of the four equatorial hexagons in the 
cuboctahedron.

We now have

\blem
\label{geomp2}
$P_{2}(c,d)$ is a geometric polyhedron, except when $d=kc$ and $k$ is an integer with $k \equiv
2 \mod 4$.
\elem

\bpf
The case $c=0$ is obvious. Suppose $c\neq 0$. We have altogether eight translation classes of
vertices and eight vertex-stars, unless $d=0$ or $d=kc$ with $k\in\BZ$ and $k \equiv 2 \mod 4$.
Moreover, the vertex-stars at vertices of distinct translation classes are also distinct, so each
translation class of vertices is characterized by its unique vertex-star. Hence $P_{2}(c,d)$ cannot
have multiple vertices. To see this more explicitly, let $x$ be a vertex of $P_{2}(c,d)$ with
vertex-star $W_{j}$ (say), and let $y$ be a vertex adjacent to $x$ such that $y-x \in W_{j}$. We
must prove that $x-y$ belongs to the vertex-star $W_{k}$ (say) at $y$, that is, $x \in y + W_{k}$.
Reducing $x$ modulo $4c\BZ^3$ (if need be), we may further assume that $x$ is among the eight
vertices listed in Lemmas~\ref{vertptwo}. Now inspection of the corresponding adjacent vertices and
their vertex-stars shows that the required property holds in each case.

If $d=kc$ with $k\in\BZ$ and $k \equiv 2 \mod 4$, then we have four translation classes of
vertices but still eight vertex-stars. Hence each point occupied by a vertex is occupied by two
vertices, so we have a double vertex associated with two (disjoint) vertex-stars; then this
point has ``valency" $6$. Finally, when $d=0$, there are eight translation classes of vertices and
four vertex-stars, so each vertex-star is associated with two translation classes. Now we have a
regular geometric polyhedron (see \eqref{inf3b}).
\epf

The vertex-stars of $P_{2}(c,d)$ are generally not planar, although the vertex-figures are
triangles. In fact, we have 

\blem
\label{ptwoplanar}
The vertex-stars of $P_{2}(c,d)$ are planar if and only if $d=0$. 
\elem

\bpf
The three vectors in $W_0$ have determinant $\pm d(3c^{2}+d^{2})$, so they lie in a plane if
and only if $d=0$.  
\epf

Next we determine the regular polyhedra. 

\blem
\label{regptwo}
The polyhedron $P_{2}(c,d)$ is geometrically chiral if $c,d\neq 0$, or geometrically regular if
$c=0$ or $d=0$. In particular, $P_{2}(c,0)$ is similar to $P_{2}(1,0) = \{\infty,3\}^{(b)}$, and 
$P_{2}(0,d)$ is similar to $P_{2}(0,1) = \{4,3\}$.
\elem

\bpf
The second claim was already settled (see \eqref{inf3b} and \eqref{p201}). Now suppose $c,d \neq
0$.  As in the proof of Lemma~\ref{regpone} we must refute the assumption that there exists an
involutory symmetry $R$ which fixes $o$ and $v$ and interchanges the neighbors
$(-c,d,c)$ and $(d,c,-c)$ of $o$. Since $W_0$ is not planar when $d\neq 0$, this symmetry $R$
must necessarily be the reflection in the plane through $o,v$ with normal vector 
\[ n := (-c,d,c) - (d,c,-c) = (-c-d,-c+d,2c) . \]  
Moreover, since $R^{-1}TR = T$, the rotation axis of $T$, with direction vector $e_{1}+e_{2}$, must
necessarily be invariant under $R$. But $e_{1}+e_{2}$ is not a scalar multiple of $n$, so
$e_{1}+e_{2}$ must lie in the mirror of $R$. Hence $-2c = n \cdot (e_{1}+e_{2}) = 0$, which is a
contradiction.
\epf

The discussion of enantiomorphism is similar to that in the previous section. Now we have
\beq
\label{ptwocddc}
P_{2}(-c,d) = P_{2}(c,d)R ,
\eeq
where $R$ is the reflection (in the plane $\xi_{1}=\xi_{2}$) given by 
\beq
\label{errtwo}
R\colon\;  (\xi_{1},\xi_{2},\xi_{3}) \; \mapsto \;  (\xi_{2},\xi_{1},\xi_{3}) .
\eeq
In this case conjugation by $R$ transforms the pair of generators $T(c,d),S_{2}$ for $G_{2}(c,d)$ 
to the pair of generators $T(-c,d),S_{2}^{-1}$ for $G_{2}(-c,d)$, so in particular we have
$G_{2}(-c,d) = R^{-1}G_{2}(c,d)R$. However, $R\not\in G_{2}(c,d)$, because $G_{2}(c,d)$ consists
only of direct isometries. As before, if $\Ph(c,d) := \{F_{0}(c,d),F_{1}(c,d),F_{2}(c,d)\}$ is the
base flag of $P_{2}(c,d)$, and if $\Ph^{2}(c,d) := \{F_{0}(c,d),F_{1}(c,d),F'_{2}(c,d)\}$ denotes
the adjacent flag differing from $\Ph(c,d)$ in its $2$-face, then we can use arguments for
$G_{2}(c,d)$ similar to those in \eqref{vpone}, \eqref{epone} and \eqref{fapone} for $G_{1}(a,b)$ to
establish \eqref{ptwocddc} as well as
\[ \Ph(-c,d)R = \Ph^{2}(c,d), \quad \Ph^{2}(-c,d)R = \Ph(c,d) . \] 
Here the equations take the form
\beq
\label{eqforg2}
\bry{c}
V(P_{2}(-c,d)) = o G_{2}(-c,d) = o R^{-1}G_{2}(c,d)R = o G_{2}(c,d)R = V(P_{2}(c,d))R ,\\[.05in]
\bry{lclcl}
F_{1}(-c,d) G_{2}(-c,d) &=& F_{1}(-c,d) R^{-1}G_{2}(c,d)R &=& (F_{1}(c,d) G_{2}(c,d))R ,\\[.02in]
F_{2}(-c,d) G_{2}(-c,d) &=& F_{2}(-c,d) R^{-1}G_{2}(c,d)R &=& (F_{2}(c,d)' G_{2}(c,d))R .
\ery
\ery
\eeq
One implication is again that Wythoff's construction applied to $G_{2}(c,d)$ with the alternative 
pair of generators $T(c,d),S_{2}^{-1}$ and the same initial vertex, yields the same polyhedron
$P_{2}(c,d)$ as before, but now with a different base flag, namely $\Ph^{2}(c,d)$. 

Alternatively we could have appealed to analogous results for the polyhedra $Q(c,d)$, again 
observing that conjugation by $R$ commutes with the facetting operation $\p_{2}$ (see
\cite[\S 6]{schu}); in particular, 
\[ (Q(c,d)R)^{\p_{2}} = (Q(c,d)^{\p_{2}})R . \]

We can also tell when two polyhedra are affinely equivalent.  The result is consistent with 
\cite[Lemma 6.10]{schu} for the polyhedra $Q(c,d)$.

\blem
\label{affeqtwo}
Let $c,d,e,f$ be real numbers, and let $(c,d) \neq (0,0) \neq (e,f)$. Then the polyhedra
$P_{2}(c,d)$ and $P_{2}(e,f)$ are affinely equivalent if and only if $(e,f) = s(\pm c,d)$ for some
non-zero scalar $s$. Moreover, $P_{2}(c,d)$ and $P_{2}(e,f)$ are congruent if and only of $(e,f) =
(\pm c,\pm d)$.
\elem

\bpf
Suppose we have an affine transformation $S$ with $P_{2}(e,f)= P_{2}(c,d)S$. As in the proof of 
Lemma~\ref{affeq} we can assume that $oS=o$ and $(c,-c,d)S=(e,-e,f)$.  First note that, since $S$ 
and $S^{-1}$ preserve planarity of vertex-stars, we certainly have $d=0$ if and only if $f=0$,
so in this case we are done. Now let $d\neq 0$, so that the vertex-stars are not planar. Then
we can argue as in Lemma~\ref{affeq} and conclude that $S^{-1}T(c,d)S = T(e,f)$ and
$S^{-1}S_{2}S = S_{2}$ or $S_{2}^{-1}$.  In particular, we must have $S=sI$ or $S = sR$, with $R$ as
in \eqref{errtwo}, for some scalar $s$. It follows that $(e,f) = s(\pm c,d)$. Moreover, if $S$ is
an isometry, then $s=\pm 1$. The converse is obvious (by \eqref{ptwocddc}).
\epf

\begin{figure}[hbt]
\centering
\begin{center}
\begin{picture}(180,235)
\put(80,110){
\begin{picture}(100,100)
\put(-115,0){\line(1,0){235}}
\put(120,0){\vector(1,0){2}}
\put(126,-3){$c$}
\put(0,-115){\line(0,1){230}}
\put(0,115){\vector(0,1){2}}
\put(-2,121){$d$}
\put(-80,0){\line(1,1){80}}
\put(-80,0){\line(1,-1){80}}
\put(80,0){\line(-1,1){80}}
\put(80,0){\line(-1,-1){80}}
\multiput(79.4,0)(0.3,0.3){5}{\line(-1,1){80}}
\put(80,0){\circle*{5}}
\put(-80,0){\circle*{5}}
\put(0,80){\circle*{5}}
\put(0,-80){\circle*{5}}
\put(0,-100){\vector(1,0){10}}
\put(0,-100){\vector(-1,0){10}}
\put(-15,-113){$R$}
\put(77,-15){$\{\infty,3\}^{(b)}$}
\put(3,82){$\{4,3\}$}
\put(3,-86){$\{4,3\}$}
\put(-125,-15){$\{\infty,3\}^{(b)}$}
\put(30,50){\circle*{5}}
\put(30,-50){\circle*{5}}
\put(37,48){\footnotesize ${P_{2}(c,d)}$}
\put(-30,50){\circle*{5}}
\put(-30,-50){\circle*{5}}
\put(-78,48){\footnotesize ${P_{2}(-c,d)}$}
\put(78,5){\scriptsize $1$}
\put(-8,78){\scriptsize $1$}
\end{picture}}
\end{picture}
\caption{\it The affine classes of polyhedra $P_{2}(c,d)$.}
\label{affclassptwo}
\end{center}
\end{figure}
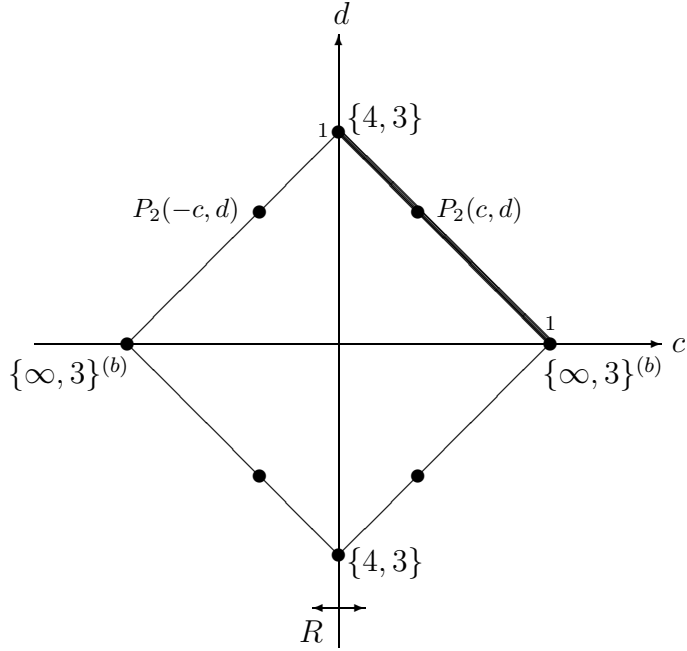

The affine equivalence classes of polyhedra $P_{2}(c,d)$ are illustrated in 
Figure~\ref{affclassptwo}. At the corners of the square we find the regular polyhedra $P_{2}(1,0)
= \{\infty,3\}^{(b)}$ and $P_{2}(0,1) = \{4,3\}$. Each class is represented by one polyhedron whose
parameter point $(c,d)$ is located on the upper right side of the square; the other three sides of
the square give alternative representations, as indicated by the four symmetrically related nodes.

As $G_{2}(c,d)$ contains only direct isometries, the faces of $P_{2}(c,d)$ consist of either
all right-hand helices or all left-hand helices. Each face has its axis parallel to a coordinate
axis. The reflection $R$ of \eqref{errtwo} changes orientation, so the faces of $P_{2}(c,d)$ are
right-hand helices if and only if the faces of $P_{2}(-c,d)$ are left-hand helices. This is
similar to what we observed in the previous section. The following lemma describes the
translation classes of faces modulo $4c\BZ^3$; in the regular case they coincide with the
translation classes relative to the full symmetry group. There are six classes, even when there are
only four translation classes of vertices.

\blem 
\label{trclptwo}
$P_{2}(c,d)$ has precisely six translation classes of faces, namely those represented by the
faces $F_{2}$, $F_{2}T_{1}$, $F_{2}S_{2}$, $F_{2}T_{1}S_{2}$, $F_{2}S_{2}^{2}$ and
$F_{2}T_{1}S_{2}^{2}$, where $F_{2}$ is the base face given by \eqref{facetsptwo}.
\elem

\bpf
The details are tedious, so we just outline the argument. Suppose $c\neq 0$. Each vertex of
$P_{2}(c,d)$ is equivalent modulo $4c\BZ^{3}$ to one of the eight vertices $o$, $oU_{1}$, $oU_{2}$,
$oU_{3}$, $oT_{1}$, $oU_{1}T_{1}$, $oU_{2}T_{1}$ or $oU_{3}T_{1}$ (see the proof of
Lemma~\ref{vertptwo}), so each face is equivalent to a face containing one of these vertices. This
latter set of faces consists of the three faces $F_{2}$, $F_{2}S_{2}$, $F_{2}S_{2}^{2}$ containing
$o$, as well as their images under $U_{1}$, $U_{2}$, $U_{3}$, $T_{1}$, $U_{1}T_{1}$,
$U_{2}T_{1}$ and $U_{3}T_{1}$. Inspection of these faces then reduces the number of classes to
$6$, with $2$ classes for each possible direction of axes. The representative faces include
$F_{2}$ and $F_{2}T_{1}$, which are parallel to the $y$-axis, as well as their images under
$S_{2}$ and $S_{2}^{2}$ (see also \eqref{altsys}). The six faces are those depicted in
Figure~\ref{fanear}, with $F_{2}T_{1} = F_{2}S_{2}^{2}TS_{2}$.
\epf

\section{Type $\{\infty,4\}$, with helical faces over triangles}
\label{helfacthree}

In this section we describe the third family of chiral apeirohedra with helical faces. In this case
the polyhedra are of type $\{\infty,4\}$ and have faces consisting of helices over triangles.
Again we employ Wythoff's construction.

Now we must begin with an irreducible group $G=\scl{S_{1},S_{2}}$ with special group $G_{0} =
[3,4]^{+}$,  where $S_{1}',S_{2}$ are standard generators of $G_0$ considered as the group of the
octahedron $\{3,4\}$, so that $S_{1}'$ and $S_{2}$ have periods $3$ and $4$, respectively, and
their product has period $2$. As in similar situations discussed before, once $S_{2}$ has been
selected, there are four equivalent choices for $S_{1}'$ such that $T'=S_{1}'S_{2}$ has period $2$;
if $S_{1}'$ is one of them, then the other three are $S_{2}^{-1}S_{1}'S_{2}$,
$S_{2}^{-2}S_{1}'S_{2}^{2}$ and $S_{2}^{-3}S_{1}'S_{2}^{3}$. Moreover, as before, replacing $S_{2}$
by its inverse would give the same polyhedron.

Again we consider specific representations for the generators. Now we take
\beq 
\label{helthi}
\bry{rccl}
S_{1}'\colon & x & \mapsto &  (\xi_{3},\xi_{1},\xi_{2}) ,\\
S_{2}\colon  & x & \mapsto &  (\xi_{3},\xi_{2},-\xi_{1})  \\
\ery  
\eeq
as generators of $G_{0}$; these are the inverses of the elements $S_{2}$ and $S_{1}'$ of
\eqref{helsqasg}. Then
$T'$ and
$T$ are the same transformations as in the previous section. However, the resulting group $G =
G_{3}(c,d)$ is not the same as $G_{2}(c,d)$; in fact, the two groups have different translation
subgroups. In particular, $G_{3}(c,d)$ is generated by
\beq 
\label{helthir}
\bry{rccl}
S_{2}\colon  & x & \mapsto &  (\xi_{3},\xi_{2},-\xi_{1}),  \\
T    \colon  & x & \mapsto &  (\xi_{2},\xi_{1},-\xi_{3}) + (c,-c,d),
\ery
\eeq
with real parameters $c$ and $d$, not both zero. Then
\beq 
\label{helthird}
\bry{rccl}
S_{1}    \colon & x & \mapsto &  (\xi_{3},\xi_{1},\xi_{2}) + (-d,-c,c), \\
S_{1}^{3}\colon & x & \mapsto &  (\xi_{1},\xi_{2},\xi_{3}) + d(-1,-1,-1). \\ 
\ery  
\eeq
Then $S_{1}^{3}$ is the translation by $d(-1,-1,-1)$, and is trivial if $d=0$. The special case
$d=0$ gives the finite group 
\[ G_{3}(c,0) \cong [3,4]^+  \]
(with fixed point $(0,-c,0)$); the polyhedron is an octahedron. 

The polyhedron $P = P_{3}(c,d)$ is obtained from $G=G_{3}(c,d)$ by Wythoff's construction with base
vertex $F_{0}=o$. Its base edge $F_{1}$ has vertices $o$ and $v:=oT = (c,-c,d)$, and hence is
the same as for $P_{2}(c,d)$. The base face $F_{2}$ is given by
\beq
\label{facetspthree}
F_{2} = o\scl{S_{1}} = \{(-d,-c,c),(0,0,0),(c,-c,d)\} + {\Bbb Z}\!\cdot\!d(1,1,1) 
\eeq
and consists of a helix over a triangle. All other vertices, edges and faces of $P$ are the
images of $F_{0}$, $F_{1}$ and $F_{2}$ under $G$. The vertices adjacent to $o$ are 
\[ v = (c,-c,d),\;\, vS_{2} = (d,-c,-c),\;\, 
vS_{2}^{2} = (-c,-c,-d),\;\, vS_{2}^{3} = (-d,-c,c) .\] 
The polyhedron $P_{3}(c,d)$ is similar to $P_{3}(1,\tfrac{d}{c})$ or $P_{3}(0,1)$ according as
$c\neq 0$ or $c=0$, so again there is only a single parameter for the similarity classes. 

There does not seem to be any obvious relationship between $P_{3}(c,d)$ and any of the chiral
polyhedra with finite faces, so in particular the facetting operation $\p_{2}$ does not apply
here. We shall elaborate on this in Section~\ref{relships}.

Next we find the translation subgroup $T(G)$ of $G=G_{3}(c,d)$. Since $S_{1}^3$ is the translation
by $d(-1,-1,-1)$, and its conjugates under the special group $G_{0} \;(\cong [3,4]^{+})$ consists
of the vectors $d(\pm 1,\pm 1,\pm 1)$, we see that
\[ d\La_{(1,1,1)} \leq T(G) .\]
In fact, we have 

\blem
\label{pthreetrans}
The translation subgroup of $G_{3}(c,d)$ is given by $T(G_{3}(c,d)) = d\La_{(1,1,1)}$.
\elem

\bpf
Now $G = N \cdot \scl{S_{2}}$, where $N := \scl{T_{0},T_{1},T_{2},T_{3}}$ with $T_{j} := 
S_{2}^{-j}TS_{2}^{j}$ for $j=0,1,2,3$. (The labelling of the half-turns $T_j$ differs from that in
previous sections.)  In particular,
\beq 
\label{t123pthree}
\bry{rccl}
T_{0} = T\colon & x & \mapsto &  (\xi_{2},\xi_{1},-\xi_{3}) + (c,-c,d),  \\
T_{1}\colon & x & \mapsto &  (-\xi_{1},-\xi_{3},-\xi_{2}) + (d,-c,-c), \\
T_{2}\colon & x & \mapsto &  (-\xi_{2},-\xi_{1},-\xi_{3}) + (-c,-c,-d),\\
T_{3}\colon & x & \mapsto &  (-\xi_{1},\xi_{3},\xi_{2}) + (-d,-c,c). 
\ery  
\eeq
Note that $N$ is normal in $G$.

The element $T_{0}T_{1}T_{2}T_{3}T_{2}T_{1}$ of $N$ is the translation by $(0,0,-2d)$. Its 
conjugates by $S_{1}$ and $S_{1}^{2}$ are also in $N$ and are the translations by $(-2d,0,0)$
and $(0,-2d,0)$, respectively. Hence $2d\BZ^{3} \leq N$; in fact, $2d\BZ^{3}$ is normal in $N$.
Moreover, $S_{2}^{2}T_{0}T_{2}T_{1}T_{3}$ is the translation by $(-2d,0,-2d)$ contained in $N$, so 
$S_{2}^{2}\in N$ and $N$ has index at most $2$ in $G$.

Now consider the quotient $M := N\slash 2d\BZ^{3}$. The elements $(T_{j}T_{j+1})^{3}$ with
$j=0,1,2,3$ (considered modulo $4$), as well as $(T_{0}T_{2})^{2}$, are translations in
$2d\BZ^{3}$, hence each is trivial when considered in $M$. Moreover, since
$T_{0}T_{1}T_{2}T_{3}T_{2}T_{1}$ is also in $2d\BZ^{3}$, we have $T_{3} =
T_{2}T_{1}T_{0}T_{1}T_{2}$ in $M$, so $M$ is generated by (the images of) $T_{0},T_{1},T_{2}$.
These generators satisfy the standard Coxeter-type relations for $[3,3]\;(\cong S_{4})$, and since
$T_{0}T_{1}$, $T_{1}T_{2}$ and $T_{0}T_{2}$ themselves are not translations, these relations give a
presentation for $M$, so in particular $M\cong S_{4}$. On the other hand, the normal subgroup
$d\La_{(1,1,1)}$ of $G$ (contained in $T(G)$) is not a (normal) subgroup of $N$; otherwise,
since $2d\BZ^{3}$ has index $2$ in $d\La_{(1,1,1)}$, the group $N\slash d\La_{(1,1,1)}$ would
be a quotient of $N\slash 2d\BZ^{3}$ of order $12$ (however, $S_4$ does not have such a
quotient). Hence $N$ must be a proper subgroup of $G$ (that is, $S_{2} \not\in N$), and $G =
N\cdot T(G)$. Now it follows that
\[ N\slash 2d\BZ^{3} = M \cong S_{4} \cong [3,4]^{+} \cong G_{0} \cong 
G\slash T(G) \cong N\slash (T(G) \cap N), \] 
so in particular $T(G) \cap N = 2d\BZ^{3}$ and $2d\BZ^{3}$ has index $2$ in $T(G)$. Hence $T(G) =
d\La_{(1,1,1)}$, as claimed. 
\epf

Again we observe that the translation group depends only on a single parameter, $d$ in this
case. Now Lemma~\ref{pthreetrans} can be rephrased as
\beq
\label{specgthree}
G_{3}(c,d)\slash d\La_{(1,1,1)}  \cong [3,4]^{+} \cong S_{4} ,
\eeq
where each term represents the special group. 

The next lemma describes the vertex-set of $P_{3}(c,d)$, which again is discrete for all
parameter values.

\blem
\label{vertpthree}
The vertex-set of $P_{3}(c,d)$ is given by
\[ V(P_{3}(c,d)) = 
\{(0,0,0), (c,-c,d), (-c,-c,d), (0,-2c,0), (d,-c,c), (d,-c,-c) \} + d\La_{(1,1,1)}. \]
The six cosets of $d\La_{(1,1,1)}$ occurring are distinct, except when $c/d$ is an integer. If
$c/d$ is an odd integer, then $V(P_{3}(c,d)) = d\La_{(1,1,1)}$ (and the cosets all
coincide). If $c/d$ is an even integer, then 
\[ V(P_{3}(c,d)) = \{(0,0,0), (c,-c,d), (d,-c,c)\} + d\La_{(1,1,1)} \]
(and the cosets coincide in pairs).
\elem

\bpf
We appeal to the proof of Lemma~\ref{pthreetrans}. First observe that, modulo $d\La_{(1,1,1)}$,
each element of $G$ can be written as 
\beq
\label{elrep}
S_{2}^{i} (T_{0}T_{1}T_{2})^{j} (T_{0}T_{1})^{k} T_{0}^{m} 
\eeq
with $i=0,1$, $j=0,1,2,3$, $k=0,1,2$ and $m=0,1$. In fact, the two elements $I,S_{2}$ give a
system of left coset representatives of $G$ modulo its subgroup $N$, and the four elements
$(T_{0}T_{1}T_{2})^{j}$ with $j=0,1,2,3$ give a system of left coset representatives of 
\[ N\slash 2d\BZ^{3} = M =\scl{T_{0},T_{1},T_{2}} \cong [3,3] \] 
modulo its subgroup $\scl{T_{0},T_{1}}\;(\cong S_{3})$; moreover, since  $2d\BZ^{3}\leq
d\La_{(1,1,1)}$, representations of elements modulo $2d\BZ^{3}$ also yield representations modulo
$d\La_{(1,1,1)}$.

The vertex-set of $P_{3}(c,d)$ is the orbit of $o$ under $G$ and consists of cosets of
$d\La_{(1,1,1)}$, namely those of the images of $o$ under the elements in \eqref{elrep}. The first
term in \eqref{elrep} fixes $o$, because $S_{2}$ fixes $o$. Next, it is easily verified that
$T_{0}T_{1}T_{2} = S_{2}^{-1}$ in $G\slash d\La_{(1,1,1)}$, so the second term of \eqref{elrep}
also does not contribute a non-trivial coset. Thus contributions can only arise from the images of
$o$ under the six elements 
\beq 
\label{elrepleft}
(T_{0}T_{1})^{k} T_{0}^{m},
\eeq 
with $k=0,1,2$ and $m=0,1$. Modulo $d\La_{(1,1,1)}$, these images are $(0,0,0)$, $(-c,-c,d)$,
$(d,-c,c)$ if $k=0,1,2$ and $m=0$, and $(c,-c,d)$, $(0,-2c,0)$, $(d,-c,-c)$  if $k=0,1,2$ and
$m=1$. They give a set of six vertices modulo $d\La_{(1,1,1)}$.

The corresponding six cosets are distinct, unless two representing vertices become equivalent
modulo $d\La_{(1,1,1)}$. Inspection shows that such equivalencies occur precisely when ${c/d}
\in \BZ$. If $c/d$ is an odd integer, then the cosets all coincide and hence
$V(P_{3}(c,d)) = d\La_{(1,1,1)}$. However, if $c/d$ is an even integer, then there are three
distinct cosets represented by $(0,0,0)$, $(c,-c,d)$ or $(d,-c,c)$. This completes the proof.
\epf

Note that, when $d=0$, the polyhedron is the octahedron $P_{3}(c,0)$ whose vertices are the six
points listed in Lemma~\ref{vertpthree}.

Next we determine the vertex-stars of $P_{3}(c,d)$. They are the images of the vertex-star
$W_0$ at $o$ under the special group $G_{0}\;(\cong S_{4})$, so the latter certainly acts
transitively on them. On the other hand, since $W_0$ is stabilized by $\scl{S_{2}}$, the number of
vertex-stars cannot exceed $6$. This bound is consistent with the numbers $6$, $3$ or $1$ of
transitivity classes of vertices.

The vertex-stars at the six representing vertices of Lemma~\ref{vertpthree} are given by
\beq 
\label{vthreestars}
\bry{lllll}
W_{0} &:=&                          & & \{ (c,-c,d),(d,-c,-c),(-c,-c,-d),(-d,-c,c) \},\\
W_{1} &:=& W_{0}T_{0}'              &=& \{ (-c,c,-d),(-c,d,c),(-c,-c,d), (-c,-d,-c) \},\\
W_{2} &:=& W_{0}T_{0}'T_{1}'        &=& \{ (c,c,d),(c,d,-c),(c,-c,-d),(c,-d,c) \},\\
W_{3} &:=& W_{0}T_{0}'T_{1}'T_{0}'  &=& \{ (d,c,c),(-c,c,d),(-d,c,-c),(c,c,-d) \}, \\
W_{4} &:=& W_{0}T_{1}'T_{0}'        &=& \{ (d,c,-c),(-c,d,-c),(-d,-c,-c),(c,-d,-c) \},\\ 
W_{5} &:=& W_{0}T_{1}'              &=& \{ (-d,c,c),(c,d,c),(d,-c,c),(-c,-d,c) \} 
\ery  
\eeq
(see \eqref{elrepleft}), where the labelling matches the order in which the vertices are listed in
the lemma. Note here that
$(T_{0}'T_{1}')^{3}=I$. The vertex-stars at the vertices $(c,-c,d)$, $(d,-c,-c)$,
$(-c,-c,-d)$ and $(-d,-c,c)$ adjacent to $o$ are $W_{1}$, $W_{5}$, $W_{2}$ and $W_{4}$,
respectively. The vertex-star at an arbitrary vertex $x$ is found by first reducing $x$ modulo
$d\La_{(1,1,1)}$ to a representing vertex and then taking the vertex-star at that vertex. 

\begin{figure}[hbt]
\centering
\begin{center}
\begin{picture}(220,250)
\multiput(10,20)(0,72){3}{ 
\begin{picture}(150,40)
\thinlines
\multiput(0,0)(72,0){3}{\circle*{1}}
\multiput(30,20)(72,0){3}{\circle*{1}}
\multiput(60,40)(72,0){3}{\circle*{1}}
\multiput(0,0)(30,20){3}{\line(1,0){144}}
\multiput(0,0)(72,0){3}{\line(3,2){60}}
\end{picture}}
\put(10,20){ 
\begin{picture}(150,200)
\thinlines
\multiput(0,0)(72,0){3}{\line(0,1){144}}
\multiput(30,20)(72,0){3}{\line(0,1){144}}
\multiput(60,40)(72,0){3}{\line(0,1){144}}
\end{picture}}
\put(10,20){ 
\begin{picture}(150,200)
\multiput(54,54)(36,0){2}{\circle*{6}}
\put(54,54){\circle{10}}
\multiput(54,90)(36,0){2}{\circle*{6}}
\multiput(22.5,68.83)(0,36){2}{\circle*{6}}
\multiput(37.5,78.83)(0,36){2}{\circle*{6}}
\multiput(76.5,15)(36,0){2}{\circle*{6}}
\multiput(91.5,25)(36,0){2}{\circle*{6}}
\put(91.5,25){\circle{10}}
\end{picture}}
\put(214,204){ 
\begin{picture}(150,200)
\multiput(-54,-54)(-36,0){2}{\circle*{6}}
\put(-90,-54){\circle{10}}
\multiput(-54,-90)(-36,0){2}{\circle*{6}}
\multiput(-22.5,-68.83)(0,-36){2}{\circle*{6}}
\multiput(-37.5,-78.83)(0,-36){2}{\circle*{6}}
\multiput(-76.5,-15)(-36,0){2}{\circle*{6}}
\multiput(-91.5,-25)(-36,0){2}{\circle*{6}}  
\put(-127.5,-25){\circle{10}}

\end{picture}}
\put(10,20){ 
\begin{picture}(150,200)
\thicklines
\multiput(54,54)(0,36){2}{\line(1,0){36}}
\multiput(54,54)(36,0){2}{\line(0,1){36}}
\multiput(22.5,68.83)(0,36){2}{\line(3,2){15}}
\multiput(22.5,68.83)(15,10){2}{\line(0,1){36}}
\multiput(76.5,15)(15,10){2}{\line(1,0){36}}
\multiput(76.5,15)(36,0){2}{\line(3,2){15}}
\end{picture}}
\put(214,204){ 
\begin{picture}(150,200)
\thicklines
\multiput(-54,-54)(0,-36){2}{\line(-1,0){36}}
\multiput(-54,-54)(-36,0){2}{\line(0,-1){36}}
\multiput(-22.5,-68.83)(0,-36){2}{\line(-3,-2){15}}
\multiput(-22.5,-68.83)(-15,-10){2}{\line(0,-1){36}}
\multiput(-76.5,-15)(-15,-10){2}{\line(-1,0){36}}
\multiput(-76.5,-15)(-36,0){2}{\line(-3,-2){15}}
\end{picture}}
\put(10,20){ 
\begin{picture}(150,200)
\thinlines
\put(102,92){\vector(1,0){116}}
\put(222,90){$\xi_{2}$}
\put(102,92){\vector(0,1){116}}
\put(99,214){$\xi_{3}$}
\put(102,92){\vector(-3,-2){60}}
\put(32,47){$\xi_{1}$}
\put(102,92){\circle*{2}}
\end{picture}}
\end{picture}
\caption{\it The set $W$ for $c=1$ and $d=4$. Now the points in the vertex-star $W_{0}$ for
$P_{3}(c,d)$ are circled.}
\label{pthreestarw}
\end{center}
\end{figure}
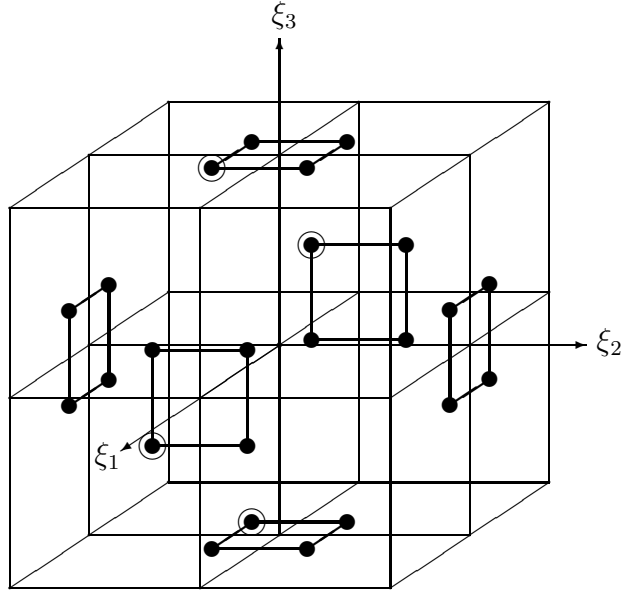

The union $W$ of the vertex-stars of $P_{3}(c,d)$ is the same set as in the previous section,
that is,
\beq
\label{doubuthr}
W  :=  \bigcup_{i=0}^{5} W_{i}  
= \{ (\pm c, \pm c,\pm d),  (\pm c, \pm d, \pm c), (\pm d, \pm c, \pm c) \} .     
\eeq 
Hence the edge-module $\La := \BZ[W]$ of $P_{3}(c,d)$ is also the same as before.

The six vertex-stars in \eqref{vthreestars} are mutually disjoint (and hence distinct), unless
$c =\pm d$ or $c=0$ or $d=0$; Figure~\ref{pthreestarw} illustrates the vertex-star $W_{0}$.
The vertex-stars are still mutually distinct (but not disjoint) if $c = \pm d$ or $d=0$; in these
cases they correspond to the square faces of a cube or cuboctahedron, respectively. Finally, if
$c=0$, the vertex-stars coincide in pairs and consist of the equatorial squares of an octahedron; in
particular, 
\beq
\label{coincpairs}
W_{0}=W_{3},\; W_{1}=W_{2},\; W_{4}=W_{5}. 
\eeq

Next we prove

\blem
\label{geomp3}
$P_{3}(c,d)$ is a geometric polyhedron, except when $c/d$ is a non-zero integer.
\elem

\bpf
The case $d=0$ is trivial, so let $d \neq 0$. Then we have six translation classes of vertices and
six vertex-stars, except when $c/d$ is an integer (possibly zero). Moreover, the vertex-stars at
vertices of distinct translation classes are distinct, so again each translation class of vertices
is characterized by its unique vertex-star. Hence $P_{3}(c,d)$ cannot  have multiple vertices. Again
we can see this more explicitly.  Let $x$ be a vertex of $P_{3}(c,d)$ with vertex-star $W_{j}$, and
let $y$ be a vertex adjacent to $x$ such that $y-x \in W_{j}$. We must verify that $x-y$ belongs to
the vertex-star $W_{k}$ at $y$, that is, $x \in y + W_{k}$. Reducing $x$ modulo $d\La_{(1,1,1)}$ we
can achieve that $x$ is among the six vertices listed in Lemmas~\ref{vertpthree}. Inspection of the
corresponding adjacent vertices and their vertex-stars now establishes the required property in
each case.

If $c/d$ is a non-zero even integer, then we have three translation classes of vertices and six 
vertex-stars, so again each vertex has multiplicity $2$ and is associated with two disjoint 
vertex-stars.  

If $c/d$ is an odd integer, then there is just one translation classes of vertices but still six
vertex-stars. Now each point $x$ taken by a vertex is in fact occupied by $6$ vertices, so each
vertex has multiplicity $6$ and all six vertex-stars occur at it. Hence the set of points
connected to $x$ by an edge of $P_{3}(c,d)$ is given by $x+W$. This case is highly degenerate.

Finally, when $c=0$, there are three translation classes of vertices and three vertex-stars, so
we have a geometric polyhedron. As we shall see, the polyhedron is regular.
\epf

Since $S_{2}$ is a rotation about the $\xi_{2}$-axis, it is clear that the four vertices of 
$P_{3}(c,d)$ adjacent to $o$ all lie in a plane perpendicular to the $\xi_{2}$-axis, so the
vertex-figures certainly are planar. Moreover, we have

\blem
\label{pthreeplanar}
The vertex-stars of $P_{3}(c,d)$ are planar if and only if $c=0$. 
\elem

\bpf
The vertex-stars are planar if and only if the affine hull of the vertex-figure at $o$ contains $o$
itself. This occurs if and only if $v = (c,-c,d)$ lies in the $\xi_{1}\xi_{3}$-plane. Hence the
condition is $c=0$.
\epf

Next we determine the regular polyhedra. Recall from \cite[p.224,230]{msarp} that 
$\{\infty,4\}_{\cdot,* 3}$ denotes the (self-Petrie) regular polyhedron of type $\{\infty,4\}$, 
whose full symmetry group is specified by the single extra relation 
\[ (R_{2}(R_{1}R_{0})^{2})^{3} = I . \]
(The ``$\cdot$" in the suffix means that the length of the Petrie polygons is unspecified, and the
$3$ with $*$ prefix indicates the length of the $2$-zigzags of the combinatorial dual.)

\blem
\label{regpthree}
The polyhedron $P_{3}(c,d)$ is geometrically chiral if $c,d\neq 0$, or geometrically regular if 
$c=0$ or $d=0$. In particular, $P_{3}(0,d)$ is similar to $P_{3}(0,1) = \{\infty,4\}_{\cdot,*3}$,
and $P_{3}(c,0)$ is similar to $P_{3}(1,0) = \{3,4\}$.
\elem

\bpf
Clearly, by construction, $P_{3}(c,d)$ is at least chiral. Now suppose $P_{3}(c,d)$ is regular. 
Then there exists an involutory symmetry $R$ which fixes $o$ and $v$ (and $vS_{2}^{2}$) and
interchanges the vertices $vS_{2}$ and $vS_{2}^{-1}$ in the vertex-figure at $o$; moreover,
$RS_{2}R=S_{2}^{-1}$ and $RTR=T$. If $c\neq 0$, the vectors $v = (c,-c,d)$ and
$vS_{2}^{2}=(-c,-c,-d)$ span an $R$-invariant plane with normal vector
\[ n := vS_{2}-vS_{2}^{-1} = (2d,0,-2c) , \]
so $R$ must be the reflection in this plane. The rotation axis of $T$ is the line through
$\tfrac{1}{2}v$ with direction vector $e_{1}+e_{2}$.  (As before, $e_{1},e_{2},e_{3}$ denotes the
standard basis of $\BE^3$.)  Since $R$ must also leave this line invariant, it must lie in the
mirror of $R$; in fact, it cannot be perpendicular to the mirror because $e_{1}+e_{2}$ and $n$ are
not collinear. But then $d = n\cdot(e_{1}+e_{2}) = 0$. Hence, if $P_{3}(c,d)$ is regular, then
either $c=0$ or $d=0$. The latter case yields the regular octahedron  $P_{3}(c,0)$.

Now let $c=0$. Suppose $R$ is an involutory symmetry of $P_{3}(0,d)$ with the properties mentioned.
Then $v = de_{3} =-vS_{2}^{2}$ and $vS_{2} = de_{1} =-vS_{2}^{-1}$, so necessarily $e_{3}R=e_{3}$
and $e_{1}R = -e_{1}$. Moreover, since $RS_{2}R=S_{2}^{-1}$, the symmetry $R$ must leave the
rotation axis of $S_{2}$ (the $y$-axis) invariant, so we also have $e_{2}R = \pm e_{2}$.  We can
further eliminate the possibility that $R$ is a plane reflection; in fact, then necessarily
$e_{2}R = e_{2}$, so $R$ would not leave the rotation axis of the half-turn $T$ invariant. This,
then, leaves only the possibility that $R$ is the half-turn about the $z$-axis (spanned by $v$). 
On the other hand, it is straightforward to check that this half-turn, denoted again by $R$,
indeed is a symmetry of $P_{3}(0,d)$ with the required properties. In particular, it fixes $o$ and
$v$ and maps the base face $F_{2}$ of \eqref{facetspthree} to the face
\[ \begin{array}{rclcl}
F_{2}R &=& \{(d,0,0),(0,0,0),(0,0,d)\} + {\Bbb Z}\!\cdot\!d(-1,-1,1) \\
       &=& \{(0,-d,d),(0,0,d),(0,0,0)\} + {\Bbb Z}\!\cdot\!d(1,1,-1) &=& F_{2}T
\end{array} \]
adjacent to $F_{2}$ along the base edge. It follows that $P_{3}(0,d)$ is regular (see
Lemma~\ref{regcrit}). Moreover, since $R_{2}:=R$ and $S_{2}$ are rotations and $S_{1}$ is a
twist, the symmetries $R_{1}:=S_{2}R_{2}^{-1}$ and $R_{0}:=S_{1}R_{1}^{-1}$ must also be
half-turns. Hence $P_{3}(0,d)$ has dimension vector 
\[ (\dim R_{0},\dim R_{1},\dim R_{2})=(1,1,1) . \]
Comparison with the regular polyhedra in $\BE^3$ now establishes that $P_{3}(0,d)$ is similar to 
$\{\infty,4\}_{\cdot,*3}$ (see \cite[p.225]{msarp}). 
\epf

The discussion of enantiomorphism follows the same pattern as in previous sections. Now 
\beq
\label{ptthreecddc}
P_{3}(c,-d) = P_{3}(c,d)R ,
\eeq
where $R$ is the reflection (in the plane $\xi_{3}=0$) given by 
\beq
\label{errthree}
R\colon\;  (\xi_{1},\xi_{2},\xi_{3}) \; \mapsto \;  (\xi_{1},\xi_{2},-\xi_{3}) .
\eeq
In this case conjugation by $R$ transforms the generators $T(c,d),S_{2}$ for $G_{3}(c,d)$ to
the generators $T(c,-d),S_{2}^{-1}$ for $G_{3}(c,-d)$, so in particular
$G_{3}(c,-d) = R^{-1}G_{3}(c,d)R$. However, $R\not\in G_{3}(c,d)$, again because $G_{3}(c,d)$ only
contains direct isometries. If $\Ph(c,d) := \{F_{0}(c,d),F_{1}(c,d),F_{2}(c,d)\}$ is the base
flag of $P_{3}(c,d)$, and if 
\[ \Ph^{2}(c,d) := \{F_{0}(c,d),F_{1}(c,d),F'_{2}(c,d)\} \] 
is the adjacent flag differing from $\Ph(c,d)$ in its $2$-face, then equations for $G_{3}(c,d)$
similar to those in \eqref{eqforg2} for $G_{2}(c,d)$ establish \eqref{ptthreecddc} and 
\[ \Ph(c,-d)R = \Ph^{2}(c,d), \quad \Ph^{2}(c,-d)R = \Ph(c,d) . \] 
In particular this implies that Wythoff's construction applied to $G_{3}(c,d)$ with the alternative 
pair of generators $T(c,d),S_{2}^{-1}$ and the same initial vertex, yields the same polyhedron
$P_{3}(c,d)$ as before, but now with $\Ph^{2}(c,d)$ as the base flag. 
\smallskip

Next we investigate affine equivalence. Here we have

\blem
\label{affeqthree}
Let $c,d,e,f$ be real numbers, and let $(c,d) \neq (0,0) \neq (e,f)$. Then the polyhedra
$P_{3}(c,d)$ and $P_{3}(e,f)$ are affinely equivalent if and only if $(e,f)=t(c,\pm d)$ for some
non-zero scalar $t$. Moreover, $P_{3}(c,d)$ and $P_{3}(e,f)$ are congruent if and only of
$(e,f) = (\pm c,\pm d)$.
\elem

\bpf
The proof is very similar to that of Lemma~\ref{affeqtwo}. Suppose $S$ is an affine transformation
with $P_{3}(e,f)= P_{3}(c,d)S$, $oS=o$ and $(c,-c,d)S=(e,-e,f)$. Since $S$ and $S^{-1}$ preserve
planarity of vertex-stars, we certainly have $c=0$ if and only if $e=0$. If $c\neq 0$, we have
non-planar vertex-stars and hence can conclude as in Lemma~\ref{affeq} that $S^{-1}T(c,d)S = T(e,f)$
and $S^{-1}S_{2}S = S_{2}$ or $S_{2}^{-1}$. Then either $S=tI$ or $S=tR$, with $R$ as in
\eqref{errthree}, for some scalar $t$. Hence $(e,f)=t(c,\pm d)$. Moreover, if $S$ is an isometry,
then $t=\pm 1$. The converse is clear (by \eqref{ptthreecddc}).
\epf

\begin{figure}[hbt]
\centering
\begin{center}
\begin{picture}(180,235)
\put(80,110){
\begin{picture}(100,100)
\put(-115,0){\line(1,0){235}}
\put(122,0){\vector(1,0){2}}
\put(128,-3){$c$}
\put(0,-115){\line(0,1){230}}
\put(0,115){\vector(0,1){2}}
\put(-2,121){$d$}
\put(-80,0){\line(1,1){80}}
\put(-80,0){\line(1,-1){80}}
\put(80,0){\line(-1,1){80}}
\put(80,0){\line(-1,-1){80}}
\multiput(79.4,0)(0.3,0.3){5}{\line(-1,1){80}}
\put(80,0){\circle*{5}}
\put(-80,0){\circle*{5}}
\put(0,80){\circle*{5}}
\put(0,-80){\circle*{5}}
\put(112,0){\vector(0,1){10}}
\put(112,0){\vector(0,-1){10}}
\put(114,6){$R$}
\put(77,-15){$\{3,4\}$}
\put(3,82){$\{\infty,4\}_{\cdot,*3}$}
\put(3,-86){$\{\infty,4\}_{\cdot,*3}$}
\put(-106,-15){$\{3,4\}$}
\put(30,50){\circle*{5}}
\put(-30,50){\circle*{5}}
\put(37,48){\footnotesize ${P_{3}(c,d)}$}
\put(30,-50){\circle*{5}}
\put(-30,-50){\circle*{5}}
\put(37,-52){\footnotesize ${P_{3}(c,-d)}$}
\put(78,5){\scriptsize $1$}
\put(-8,78){\scriptsize $1$}
\end{picture}}
\end{picture}
\caption{\it The affine classes of polyhedra $P_{3}(c,d)$.}
\label{affpthree}
\end{center}
\end{figure}

Notice that the conditions of Lemmas~\ref{affeqtwo} and \ref{affeqthree} are the same, that is,
$P_{3}(c,d)$ and $P_{3}(e,f)$ are affinely equivalent if and only if $P_{2}(c,d)$ and $P_{2}(e,f)$
are affinely equivalent. Thus the affine equivalence classes of polyhedra $P_{3}(c,d)$ can be
illustrated in a similar fashion as those of $P_{2}(c,d)$ (see Figure~\ref{affpthree}). At the
corners of the square we now find the regular polyhedra $P_{3}(0,1) = \{\infty,4\}_{\cdot,*3}$ and
$P_{3}(1,0) = \{3,4\}$. As before, each equivalence class is represented by a polyhedron with its
parameter point $(c,d)$ located on the upper right side of the square; the other three sides of the
square give alternative representations.

The faces of $P_{3}(c,d)$ consist of either all right-hand helices or all left-hand helices, and
the helices for $P_{3}(c,d)$ are right-hand if and only if those for $P_{3}(c,-d)\;(= P_{3}(c,d)R)$
are left-hand. Each face has its axis parallel to a body diagonal of $\BZ^3$. The next lemma
describes the translation classes of faces modulo $d\La_{(1,1,1)}$; in the regular case they also
are the translation classes relative to the full symmetry group. 

\blem 
\label{trclpthree}
$P_{3}(c,d)$ has eight translation classes of faces, except when $c=0$ or $c=\pm d$; they are
represented by the faces $F_{2}T_{1}^{i}S_{2}^{j}$ with $i=0,1$ and $j=0,1,2,3$, where $F_{2}$ is
the base face given by \eqref{facetspthree}. If $c=0$ or $c=\pm d$, then $P_{3}(c,d)$ has four
translation classes of faces, represented by $F_{2}S_{2}^{j}$ with $j=0,1,2,3$.
\elem

\bpf
We only outline the argument. Suppose $d \neq 0$. Each vertex of $P_{3}(c,d)$ is equivalent modulo
$d\La_{(1,1,1)}$ to one of the (generally six) vertices $o (T_{0}T_{1})^{k}T_{0}^{m}$ with $k=0,1,2$
and $m=0,1$ (see the proof of Lemma~\ref{vertpthree}, particularly \eqref{elrepleft}), so each
face must be equivalent to a face containing one of these vertices. This latter set of faces
consists of the four faces $F_{2}S_{2}^{j}$ with $j=0,1,2,3$, each of which contains $o$ and is
parallel to one of the four body diagonals of $\BZ^3$, as well as their images under the six elements
$U:=(T_{0}T_{1})^{k}T_{0}^{m}$ with $k,m$ as above. The face $\widehat{F}_{2} :=  F_{2}S_{2}^{3}$
has its axis parallel to the diagonal with direction vector $(-1,1,1)$, and since $T_{0}'$ and
$T_{1}'$ (the images of $T_{0}$ and $T_{1}$ in the special group) reverse this direction, this also
remains true for each face $\widehat{F}_{2}U$ with $U$ as above. Inspection shows that, for all
parameter values of $c$ and $d$, these six faces fall into two sets of three, where each set
consists of mutually equivalent faces modulo $d\La_{(1,1,1)}$ (all $U$ with $m=0$, or all $U$ with
$m=1$). Moreover, faces from different sets can only be equivalent if $c=0$ or $c=\pm d$; in fact,
in these cases any two faces of the six are equivalent modulo $d\La_{(1,1,1)}$. Hence, the faces
parallel to the direction vector $(-1,1,1)$ fall into two classes or one class, respectively, and
are represented by either $\widehat{F}_{2},\widehat{F}_{2}T_{0}$ or $\widehat{F}_{2}$.
This takes care of one diagonal direction. However, since the special group $G_{0}$ (in fact,
already its subgroup $\scl{S_{2}}$) acts transitively on the four body diagonals, this also holds
for other body diagonals, so we must have a total of eight or four translation classes, respectively.
Finally, since 
\[ \widehat{F}_{2}T_{0}S_{2}^{j+1} = F_{2} (S_{2}^{-1}T_{0}S_{2})S_{2}^{j} = 
F_{2}T_{1}S_{2}^{j} \]
(see \eqref{t123pthree}), the images of $\widehat{F}_{2}T_{0}$ under $\scl{S_{2}}$ are just those
of $F_{2}T_{1}$ under $\scl{S_{2}}$, so the eight or four classes, respectively, are represented by
the images under $\scl{S_{2}}$ of either $F_{2},F_{2}T_{1}$ or $F_{2}$. This completes the proof.
\epf

\section{Coverings and relationships}
\label{relships}

In this final section, we briefly investigate coverings and relationships for the families of chiral (or
regular) polyhedra described in the previous sections and \cite{schu}. Recall that we have three
large families of chiral (or regular) polyhedra with finite faces, consisting of the polyhedra
$P(a,b)$ of type $\{6,6\}$, $Q(c,d)$ of type $\{4,6\}$, or $Q(c,d)^{*}$ (the dual of $Q(c,d)$) of
type $\{6,4\}$, respectively. In addition we have the three large families of chiral (or regular)
helix-faced polyhedra $P_{1}(a,b)$, $P_{2}(c,d)$ or $P_{3}(c,d)$ discussed in this paper.

We begin with the observation that the helix-faced polyhedra $P_{1}(a,b)$, $P_{2}(c,d)$ and
$P_{3}(c,d)$ can be thought of as figures ``unraveling" a tetrahedron $\{3,3\}$, a cube
$\{4,3\}$ or an octahedron $\{3,4\}$, respectively; that is, in a way, the helical faces can
be pushed down to their bases, like springs, to become faces of $\{3,3\}$, $\{4,3\}$ or $\{3,4\}$. 
The underlying covering relationship of polyhedra is induced by the canonical homomorphism 
\[ \bry{ccl}
G & \rightarrow & G_{0} \;(\cong G \slash T(G))\\ 
S    & \mapsto & S',
\ery \]
which maps the distinguished generators $S_{1},S_{2}$ of $G$ to those of $G_{0}$. For example,
for $P_{1}(a,b)$ with group $G = G_{1}(a,b) = \scl{S_{1},S_{2}}$, the images
$S_{1}',S_{2}'\;(=S_{2})$ are distinguished generators of $G_{0} = [3,3]^{+}$. Thus we have a
covering $P_{1}(a,b) \searrow \{3,3\}$ (see \cite[Sect.2D]{msarp}).

\bthm 
\label{rav}
There are the following coverings of polyhedra:
\[ P_{1}(a,b) \searrow \{3,3\}, \quad
P_{2}(c,d) \searrow \{4,3\}, \quad
P_{3}(c,d) \searrow \{3,4\}.\]
In each case, the polyhedron on the right is the only finite polyhedron in its family and is 
obtained for $b=a$, $c=0$ or $d=0$, respectively. 
\ethm

In a similar fashion we can also derive the coverings
\beq 
\label{covq}
Q(c,d) \searrow \{4,6\hole \cdot, 2\}, \quad Q(c,d)^{*} \searrow \{4,6\hole \cdot, 2\}^{*} ,
\eeq
where $\{4,6\hole \cdot, 2\}$ denotes the regular map of type $\{4,6\}$ determined by the length,
$2$, of its $3$-holes (the $2$-holes are not needed for the specification). Recall that a $3$-hole is
an edge-path which leaves a vertex by the $3$rd edge from which it entered, always keeping to the
right or always to the left (see \cite[p.196]{msarp}).  In fact, for $Q(c,d)$ and its dual, the
special group $G_{0}\;(=[3,4])$ admits an involutory group automorphism with $S_{1}' \mapsto
S_{1}'{S_{2}'}^{2}$ and $S_{2}'\mapsto {S_{2}'}^{-1}$, so the abstract polyhedron $\CQ^{0}$
associated with $G_{0}$ and its generators $S_{1}',S_{2}'$ must actually be directly regular, with full
automorphism group $G_{0} \rtimes C_{2}$ (see \cite[Thm.1(c)]{swc}).  (Recall that an abstract 
polyhedron is directly regular if it is regular and its combinatorial rotation subgroup has index
$2$ in the full automorphism group. The corresponding surface is necessarily orientable.)   For
$Q(c,d)$, this group automorphism can be realized by (inner) conjugation in $G_{0}$ with the
reflection in the plane $\xi_{1}=\xi_{2}$ (see \cite[(6.17)]{schu}). In particular, since $\CQ^{0}$
is orientable and of genus $3$, we must have $\CQ^{0} = \{4,6\hole \cdot, 2\}$ (see
\cite[p.470]{sherk}).

Note, however, that there is no such covering associated with the polyhedra $P(a,b)$. The
generators $S_{1}',S_{2}'$ of $G_{0}\;(=[3,3]^{*})$ do not determine a (directly) regular or chiral
map in this case; in fact, $-I \in \scl{S_{1}'} \cap \scl{S_{2}'}$, so the intersection property
fails.

Next we study operations between polyhedra or their groups. Recall that a mixing operation
derives a new group $H$ from a given group $G = \scl{S_{1},S_{2}}$ by taking as generators
$\widehat{S}_{1},\widehat{S}_{2}$ (say) certain suitably chosen products of $S_{1},S_{2}$,
so that $H$ is a subgroup of $G$ (see \cite[p.183]{msarp}). Mixing operations do not always produce
groups of chiral or regular polyhedra, but in certain special cases they do. Besides the duality
operation $\delta$ and the facetting operation $\p_{2}$ described in \eqref{dualop}, we also employ
the {\em halving operation\/} 
\beq
\label{eta}
\eta:\ (S_{1},S_{2}) \mapsto
(S_{1}^{2}S_{2},S_{2}^{-1}) = : (\widehat{S}_{1},\widehat{S}_{2}) 
\eeq
(see \cite[\S 7]{schu}), which is an analog of the halving operation for regular polyhedra (see
\cite[p.197]{msarp}). Recall from \cite[\S 7]{schu} that $\eta$ links the two families of
finite-faced polyhedra $Q(c,d)$ and $P(a,b)$ through
\[ Q(c,d)^{\eta} \cong P(c-d,c+d) , \] 
where $Q(c,d)^{\eta}$ denotes the image of $Q(c,d)$ under $\eta$ (and $\cong$ means congruence). We
also observed in Lemmas~\ref{p2pab} and \ref{p2pcd} that
\[ P_{1}(a,b) = P(a,b)^{\p_2}, \quad P_{2}(c,d) = Q(c,d)^{\p_2} . \]

Generally it is difficult to decide whether or not two classes of polyhedra are related by a
mixing operation. In some cases obstructions already arise from the structure of the special
groups. In fact, we have the following lemma.

\blem
\label{obstr}
Let $G$ and $H$ be crystallographic groups with special groups $G_{0}$ and $H_{0}$, respectively.
If $H$ is a subgroup of $G$, then $H_{0}$ is isomorphic to a subgroup of $G_{0}$. 
\elem

\bpf
If $H\leq G$, then the mapping $H \rightarrow G/T(G)$ defined by $S \mapsto T(G)S$ has kernel 
$H \cap T(G) = T(H)$, so $H\slash T(H)$ is isomorphic to a subgroup of $G\slash T(G)$. Now the
lemma follows.
\epf

For example, the special groups tell us immediately that $P_{3}(c,d)$ cannot be derived from a
mixing operation applied to a polyhedron $P(a,b)$ (see \cite[\S 5]{schu}). In fact, by the previous
lemma, $G_{3}(c,d)$ cannot be a subgroup of $G(a,b)$ (the group for $P(a,b)$), because clearly
$(G_{3}(c,d))_{0}\,(= [3,4]^{+} \cong S_{4})$ is not isomorphic to a subgroup of $(G(a,b))_{0}\,(=
[3,3]^{*} \cong A_{4}\times C_{2})$. Similarly we can eliminate the possibility that $P_{3}(c,d)$
is obtained from a mixing operation applied to a helix-faced polyhedron $P_{1}(a,b)$ (with special
group $[3.3]^{+}\cong A_{4}$). 

On the other hand, we cannot, a priori, rule out a mixing operation between $P_{3}(c,d)$ and a
polyhedron $Q(c',d')$ or its dual; in fact, $Q(c',d')$ and its dual have special group
$[3,4]\;(\cong S_{4} \times C_{2})$. The same remark also applies to $P_{3}(c,d)$ and a polyhedron
$P_{2}(c',d')$ (with special group $[3,4]^{+}$). In these cases it is open whether or not mixing
operations exist. On the other hand, we do know that the regular polyhedron $P_{3}(0,1) =
\{\infty,4\}_{\cdot,*3}$ is obtained from the Petrie-Coxeter polyhedron $Q(0,1)^{*} = \{6,4 \hole
4\}$ by the {\em skewing\/} (mixing) {\em operation\/} $\sg$ defined in \cite[p.199,224]{msarp};
that is, $P_{3}(0,1) = (Q(0,1)^{*})^{\sg}$. However, $\sg$ is only defined for regular polyhedra,
and it is not clear how to generalize it meaningfully to chiral polyhedra.

We now list the families of polyhedra, with the various known relationships between them. Our
first diagram emphasizes operations between families rather than between individual polyhedra.
In particular, we drop the parameters from the notation; for example, $P_{1}$ denotes the family of
polyhedra $P_{1}(a,b)$. Observe that $P_{3}$ is not connected to any other family. The circular 
arrow indicates the self-duality of the family (or polyhedron).
\beq
\label{displayone}
\matrix{
Q^* & \stackrel{\delta}{\longleftrightarrow} \mkern-30mu & Q & \mkern-30mu 
\stackrel{\p_2}{\longrightarrow}& P_2 \cr \cr
& &\;\;\downarrow\! {\scriptstyle\eta} & & & \cr \cr
P_3 & & P & \mkern-30mu \stackrel{\p_2}{\longrightarrow}& P_1 \cr
& & 
\begin{picture}(60,60)
\put(9,46){\oval(42,42)[b]}
\put(9,46){\oval(42,42)[tl]}
\put(9,67){\vector(1,0){2}}
\put(-19,34){$\scriptstyle\delta$}
\end{picture}
&&
} \eeq

The corresponding parameters can be read off the next diagram for (congruence classes of)
individual polyhedra; here, $Q(c,d)$ is the reference point. 
\beq
\label{displaytwo}
\matrix{
Q^*(c,d) \mkern20mu & \stackrel{\delta}{\longleftrightarrow} 
& Q(c,d) & 
\stackrel{\p_2}{\longrightarrow}& P_{2}(c,d) \cr \cr
& &\;\;\downarrow {\scriptstyle\eta} & & &\cr\cr
P_{3}(c',d')  \mkern20mu & & P(c-d,c+d) & 
\stackrel{\p_2}{\longrightarrow}& P_{1}(c-d,c+d) \cr
& & 
\begin{picture}(60,60)
\put(9,46){\oval(42,42)[br]}
\put(-18,46){\oval(42,42)[tl]}
\put(-18,46){\oval(42,42)[bl]}
\put(-18,25){\line(1,0){30}}
\put(-18,67){\vector(1,0){2}}
\put(-46,36){$\scriptstyle\delta$}
\end{picture}
&&
} \eeq

It is instructive to list the families of polyhedra by the structure of their special group, along
with the two regular polyhedra occurring in each family. 
\medskip
\begin{center}
\begin{tabular}{|c|c|c|c|c|c|}
\hline
$[3,3]^{*}$ & $[4,3]$ & $[3,4]$ & $[3,3]^{+}$ & $[4,3]^{+}$ & $[3,4]^{+}$ \\
\hline\hline
$P(a,b)$ & $Q(c,d)$ & $Q(c,d)^{*}$ & $P_{1}(a,b)$ & $P_{2}(c,d)$ & $P_{3}(c,d)$ \\
\hline
$\{6,6\}_{4}$ & $\{4,6\}_{6}$ & $\{6,4\}_{6}$ & $\{\infty,3\}^{(a)}$ &
$\{\infty,3\}^{(b)}$ & $\{\infty,4\}_{\cdot,*3}$ \\
$\{6,6 \hole 3\}$ & $\{4,6 \hole 4\}$ & $\{6,4 \hole 4\}$ & $\{3,3\}$ & $\{4,3\}$ &  $\{3,4\}$ \\
\hline
\end{tabular}
\end{center}
\medskip
In this table, the columns are indexed by the special groups to which the respective polyhedra
correspond. The second row contains the families. The last two rows comprise nine of the twelve discrete pure regular polyhedra in $\BE^{3}$, namely those listed in the table of \cite[p.225]{msarp} with dimension vectors 
\[ (\dim R_{0},\dim R_{1},\dim R_{2}) =  (1,2,1), (1,1,1) , (2,1,2) ,\]  
as well as the three (finite) ``crystallographic" Platonic polyhedra. The three remaining
pure regular polyhedra $\{\infty,6\}_{4,4}$, $\{\infty,4\}_{6,4}$ and $\{\infty,6\}_{6,3}$ have
dimension vector $(1,1,2)$ and do not occur in families alongside chiral polyhedra (see
Lemma~\ref{twi}). Note that $\p_{2}$, when applicable, maps a polyhedron to one in the same row.

A final comment about the regular polyhedra $P = \{\infty,3\}^{(a)}$, $\{\infty,3\}^{(b)}$ or
$\{\infty,4\}_{\cdot,*3}$ is appropriate. In each case the full symmetry group $G(P)$ is generated
by half-turns $R_{0},R_{1},R_{2}$, so the special group $G(P)_{0}$ consists only of rotations and
hence coincides with $G_{0}$. The generators $R_{0}',R_{1}',R_{2}'$ of $G(P)_{0}$ are also
half-turns. For example, if $P = \{\infty,4\}_{\cdot,*3}$, then $G(P)_{0} = G_{0} = [3,4]^{+}$
(viewed as the rotation group of an octahedron) and the rotation axes of $R_{0}'$, $R_{1}'$ or
$R_{2}'$ are the lines through $o$ and $(-1,1,0)$, $(-1,0,1)$, or $(0,0,1)$, respectively. It is not difficult to see that, in each case, $G(P)_{0}$ is a C-group with respect to these generators (see \cite[Section 2E]{msarp}), so there exists a finite abstract regular polyhedron $\CP^{0}$ with automorphism group $G(P)_{0}$ and distinguished generators $R_{0}',R_{1}',R_{2}'$.  (The abstract polyhedron $\CQ^{0}$ mentioned earlier arises in a similar way from $\{4,6\hole 4\} = Q(0,1)$.)   Inspection shows that $\CP^{0} = \{3,3\}$, $\{4,3\}_{3}$ or $\{3,4\}_{3}$, respectively. However, except in the first case, $\CP^{0}$ does not coincide with the regular polyhedron determined by $G_{0}$ and its generators $S_{1}',S_{2}'$; we know the latter to be $\{3,3\}$, $\{4,3\}$ and $\{3,4\}$, respectively. In fact, only when $\CP^{0}$ is directly regular should we expect isomorphism for these polyhedra (see \cite[p.510]{swc}). On the other hand, the very fact that conjugation by $R_{2}'$ yields an involutory group automorphism of $G_{0}$ with $S_{1}' \mapsto S_{1}'{S_{2}'}^{2}$ and $S_{2}'\mapsto {S_{2}'}^{-1}$, correctly predicts that the polyhedron associated with $G_{0}$ and $S_{1}',S_{2}'$ is directly regular and has group $G_{0}\rtimes C_{2}$.  (Since  $\{\infty,3\}^{(a)}$ and $\{\infty,3\}^{(b)}$ are Petrie-duals, the same is true for their respective polyhedra $\CP^{0}$, namely $\{3,3\}$ and $\{4,3\}_{3}$.  Moreover, since
$\{6,4 \hole 4\}^{\sg} = \{\infty,4\}_{\cdot,*3}$, we also have  $(\{4,6\hole \cdot, 2\}^{*})^{\sg} = \{3,4\}_{3}$, which is self-Petrie.) 

Notice, however, that $\CP^{0}$ does not correspond to a geometrically regular polyhedron; in fact, no point (distinct from $o$) is admisssible as initial vertex for Wythoff's construction applied to $G(P)_{0}$ and $R_{0}',R_{1}',R_{2}'$.

These considerations also shed some light on why there are exactly two regular polyhedra in each
of the six families. In fact, given a family with special group $G_{0}$, there are precisely two
involutory isometries $R$ such that conjugation by $R$ determines a group automorphism of $G_{0}$
with $S_{1}' \mapsto  S_{1}'{S_{2}'}^{2}$ and $S_{2}'\mapsto {S_{2}'}^{-1}$ (and hence $T'\mapsto
T'$). Each $R$ is the image $R_{2}'$, in $G_{0}$, of the generator $R_{2}$ for a regular polyhedron
in the same family. For example, if $P = \{\infty,4\}_{\cdot,*3}$, then $S_{1}',S_{2}' \,(=S_{2})$
are as in \eqref{helthi} and $R$ is either the half-turn about the $\xi_{3}$-axis or the reflection
in the $\xi_{1}\xi_{2}$-plane (the two isometries $R$ have perpendicular mirrors). The first choice
of $R$ leads to $P_{3}(0,1) = \{\infty,4\}_{\cdot,*3}$ (with $R_{0}',R_{1}',R_{2}'$ as above), and
the second to $P_{3}(1,0) = \{3,4\}$.  For the three families $P_{1}(a,b)$, $P_{2}(c,d)$ and $P_{3}(c,d)$ of helix-faced polyhedra, $R$ is either the reflection in the plane spanned by the rotation axes of $T'$ and $S_{2}'$, or the half-turn about the line perpendicular to this plane at $o$. The same remains true for the three families $P(a,b)$, $Q(c,d)$ and $Q(c,d)^{*}$ of finite-faced polyhedra, except that here $S_{2}'$ is replaced by $S_{2}'^{2}$.

\vfill
\end{document}